\documentclass[final]{siamltex}%

\pdfoutput=1
\usepackage{amstext}
\usepackage{amsfonts}
\usepackage{subfigure}
\usepackage{wrapfig}

\usepackage{longtable}
\usepackage{lscape}
\usepackage{multirow}
\usepackage[pdftex]{graphicx}
\usepackage{rotating}
\usepackage{hyperref}
\usepackage{float}

\usepackage{listings,url,verbatim}    

\newfloat{program}{thp}{lop}
\floatname{program}{Program}

\floatstyle{boxed}
\restylefloat{alg}
\newif\ifpdf \ifx\pdfoutput\undefined
\pdffalse 
\else
\pdfoutput=1 
\pdftrue \fi

\ifpdf \usepackage[pdftex]{graphicx} \else \usepackage{graphicx} \fi

\usepackage{color}
\usepackage{array}

\def\headerrule{\rule[.01in]{\textwidth}{.01in}}
\def\headerruleA{\rule[.0in]{0.in}{.3in}}

\newenvironment{INDENT}{\let\\=\item \parskip .2em
  \list{}{\topsep 0em \partopsep 0em
	  \parsep\parskip \itemsep 0em
	  \leftmargin 1.8em \rightmargin 0em}\item[]}{\endlist}
\newtheorem{algorithm}{Algorithm}[section]

\newtheorem{remark}{Remark}[section]
\newenvironment{REMARK}{\begin{remark}\rm}{\end{remark}}
\newtheorem{definition1}{Definition}[section]
\newenvironment{DEFINITION}{\begin{definition1}\rm}{\end{definition1}}

\def\IF/{{\bf if}}
\def\THEN/{{\bf then}}
\def\FOR/{{\bf for}}
\def\TO/{{\bf to}}
\def\STEP/{{\bf step}}
\def\DO/{{\bf do}\begin{INDENT}}
\def\ENDFOR/{\end{INDENT}{\bf endfor}}
\def\ENDIF/{\end{INDENT}{\bf endif}}
\def\WHILE/{{\bf while}}
\def\ENDWHILE/{\end{INDENT}{\bf endwhile}}
\def\REPEAT/{{\bf repeat}\begin{INDENT}}
\def\UNTIL/{\end{INDENT}{\bf until}}
\def\BEGIN/{\begin{INDENT}{\bf begin}\begin{INDENT}}
\def\END/{\end{INDENT}{\bf end}\end{INDENT}}

\newcommand{\ben}{\begin{enumerate}}
\newcommand{\een}{\end{enumerate}}
\newcommand{\be}{\begin{equation}}
\newcommand{\ee}{\end{equation}}

\newcommand{\Q}{{\bf Q.E.D.}}

\newtheorem{prop}[theorem]{Proposition}

\title{Subspace Iteration Randomization and Singular Value Problems}

\author{M.~Gu\thanks{This research was supported in part by NSF 
    Awards CCF-0830764 and CCF-1319312, and by the DOE Office of Advanced Scientific 
Computing Research under contract number DE-AC02-05CH11231. Email: mgu@math.berkeley.edu.}}

\begin{document}

\ifpdf \DeclareGraphicsExtensions{.pdf, .jpg} \else
\DeclareGraphicsExtensions{.eps, .jpg} \fi

\maketitle

\begin{abstract} A classical problem in matrix computations is the
efficient and reliable approximation of a given matrix by a matrix of lower
rank. The truncated singular value decomposition
(SVD) is known to provide the best such approximation for any given
fixed rank. However, the SVD is also known to be very costly to
compute. Among the different approaches in the literature for
computing low-rank approximations, randomized algorithms have
attracted researchers' recent attention due to their surprising
reliability and computational efficiency in different application
areas. Typically, such algorithms are shown to compute with
very high probability low-rank approximations that are within a constant factor
from optimal, and are known to perform even better in many practical
situations. In this paper, we present a novel error analysis that
considers randomized algorithms within the subspace iteration
framework and show with very high probability that highly accurate low-rank
approximations as well as singular values can indeed be computed
quickly for matrices with rapidly decaying singular values. Such
matrices appear frequently in diverse application areas such as data
analysis, fast structured matrix computations and fast direct methods
for large sparse linear systems of equations and are the driving
motivation for randomized methods. Furthermore, we show that the
low-rank approximations computed by these randomized algorithms are
actually rank-revealing approximations, and the special case of a
rank-$1$ approximation can also be used to correctly estimate matrix
$2$-norms with very high probability. Our numerical experiments are in full
support of our conclusions.  
\end{abstract}

{\bf key words:} low-rank approximation, randomized
algorithms, singular values, standard Gaussian matrix. 

\section{Introduction} \label{Sec:Intro} Randomized algorithms have
established themselves as some of the most competitive methods for
rapid low-rank matrix approximation, which is vital in many areas of
scientific computing, including principal component
analysis~\cite{Jolliffe,RokhlinSzlamTygert} and face
recognition~\cite{MullerMagaiaHerbst,turkpent}, large scale data
compression~\cite{DrineasKannanMahoneyII,DrineasMahoneyMuthukrishnan,Tropp,mahoney2011}
and fast approximate algorithms for PDEs and integral
equations~\cite{ChengGimbutasMartinssonRokhlin,Hackbusch1999,martinsson2010,SchmitzYing2D,SchmitzYing3D,XiaGu,XiaChandrasekaranGuLi}. In
this paper, we consider randomized algorithms for low-rank
approximations and singular value approximations within the subspace iteration
framework, leading to results that simultaneously retain the reliability of
randomized algorithms and the typical faster convergence of subspace iteration
methods.

Given any $m \times n$ matrix $A$ with $m
\geq n$, its singular value decomposition (SVD) is described by the
equation 
\begin{equation}\label{Eq:svd} 
A = U \Sigma V^{T},
\end{equation} 
where $U$ is an $m \times n$ column orthogonal matrix;
$V$ is an $n \times n$ orthogonal matrix; and $\Sigma =
\diag(\sigma_1, \cdots, \sigma_n)$ with $\sigma_1 \geq \sigma_2 \geq
\cdots \geq \sigma_n \geq 0$. Writing $U$ and $V$ in terms of their
columns, \[ U = \left(u_1, \cdots, u_n\right) \quad \mbox{and} \quad V
= \left(v_1, \cdots, v_n\right) ,\] then $u_j$ and $v_j$ are the left
and right singular vectors corresponding to $\sigma_j$, the $j$-th
largest singular value of $A$. For any $1 \leq k \leq n$, we let \[
A_k = \left(u_1, \cdots, u_k\right) \diag(\sigma_1, \cdots, \sigma_k)
\left(v_1, \cdots, v_k\right)^{T} \] be the (rank-$k$) truncated SVD
of $A$. The matrix $A_k$ is unique only if $\sigma_{k+1} <
\sigma_k$. The assumption that $m \geq n > \max{k,2} $ will be maintained
throughout this paper for ease of exposition. Our results
still hold for $m < n$ by applying all the algorithms on
$A^T$. Similarly, all our main results are derived under the
assumption that ${\bf rank}(A) = n$. But they remain unchanged even if
${\bf rank}(A) < n$, and hence remain valid by a continuity argument. 
All our analysis is done without consideration of round-off errors, and thus need not hold exactly true in finite precision, especially when the user tolerances for the low-rank approximation are close to machine precision levels. Additionally, we assume throughout this
paper that all matrices are real. In general, $A_k$ is an ideal
rank-$k$ approximation to $A$, due to the following celebrated
property of the SVD: 
\begin{theorem}\label{Thm:truncSVD} (Eckart and
  Young~\cite{EckartYoung}, Golub and van Loan~\cite{golubvanloan})
\begin{eqnarray}
 {\displaystyle \min_{{\bf rank}(B) \leq k} \| A - B\|_2} &=&  \| A - A_k\|_2 = \sigma_{k+1}. \label{Eqn:truncSVD2} \\
 { \min_{{\bf rank}(B) \leq k} \| A - B\|_F} &=&  \| A - A_k\|_F = \sqrt{\sum_{j=k+1}^n \sigma_j^2}. \label{Eqn:truncSVDF} 
\end{eqnarray}
\end{theorem}
\begin{REMARK} While there are results similar to Theorem~\ref{Thm:truncSVD} for all
unitarily invariant matrix norms, our work on low-rank matrix
approximation bounds will only focus on the two most popular of such norms:
the 2-norm and the Frobenius norm. 
\end{REMARK}

Theorem~\ref{Thm:truncSVD} states that the truncated SVD
provides a rank-$k$ approximation to $A$ with the smallest
possible 2-norm error and Frobenius-norm error. In the 2-norm,
any rank-$k$ approximation will result in an error no less
than $\sigma_{k+1}$, and in the Frobenius-norm, any rank-$k$
approximation will result in an error no less than
$\sqrt{\sum_{j=k+1}^n \sigma_j^2}$. Additionally, the
singular values of $A_k$ are exactly the first $k$
singular values of $A$, and the singular vectors of $A_k$
are the corresponding singular vectors of $A$. Note,
however, that while the solution to
problem~(\ref{Eqn:truncSVDF}) must be $A_k$, solutions to
problem~(\ref{Eqn:truncSVD2}) are not unique and include, for example,  
the rank-$k$ matrix $B$ defined below for any $0 \leq \theta
\leq 1$: 
\begin{equation}\label{Eqn:B}
B  = A_k -\theta \sigma_{k+1}\left(u_1, \cdots, u_k\right) \left(v_1, \cdots,
v_k\right)^{T}. 
\end{equation}
This subtle distinction between the 2-norm and Frobenius norm
will later on become very important in our analysis of
randomized algorithms (see Remark~\ref{Rem:Fvs2}.) In Theorem~\ref{Thm:Fto2} we prove 
an interesting result related to Theorem~\ref{Thm:truncSVD} for rank-$k$
approximations that only solve
problems~(\ref{Eqn:truncSVD2}) and~(\ref{Eqn:truncSVDF})
approximately. 

To compute a truncated SVD of a general $m \times n$ matrix
$A$, one of the most straightforward techniques is to compute the
full SVD and truncate it, with a standard linear algebra
software package like the LAPACK~\cite{LAPACK}. This
procedure is stable and accurate, but it requires $O(mn^2)$
floating point operations, or {\em flops}. This is
prohibitively expensive for applications such as data
mining, where the matrices involved are typically sparse
with huge dimensions. In other practical applications
involving the truncated SVD, often the very objective of
computing a rank-$k$ approximation is to avoid excessive
computation on $A$. Hence it is desirable to have schemes
that can compute a rank-$k$ approximation more
efficiently. Depending on the reliability requirements, a
good rank-$k$ approximation can be a matrix that is accurate
to within a constant factor from the optimal, such as a
rank-revealing factorization (more below), or it can be a
matrix that closely approximates the truncated SVD itself.

Many approaches have been taken in the literature for
computing low-rank approximations, including rank-revealing decompositions based on the QR, LU,
or two-sided orthogonal (aka UTV)
factorizations~\cite{ci,FierroHansen,ge7,HP,MG,Pan,VHZ}. Recently, there has been an explosion of randomized algorithms for computing
low-rank 
approximations~\cite{ChengGimbutasMartinssonRokhlin,DrineasKannanMahoneyII,DrineasMahoneyMuthukrishnan,FriezeKannanVempala,FriezeKannanVempalaII,LibertyAilonSinger,Liberty,LibertyWoolfeMartinssonRokhlinTygert,NguyenDoTran,WoolfeLibertyRokhlinTygert,sarlos}.
There is also software package available for computing interpolative
decompositions, a form of low-rank approximation, and for computing
the PCA, with randomized
sampling~\cite{MartinssonRokhlinShkolniskyTygert}. These algorithms
are attractive for two main reasons: they have been shown to be
surprisingly efficient computationally; and like subspace methods, the
main operations involved in many randomized algorithms can be
optimized for peak machine performance on modern architectures. For a
detailed analysis of randomized algorithms and an extended reference
list, see~\cite{Tropp}; for a survey of randomized algorithms in data
analysis, see~\cite{mahoney2011}.

The subspace iteration is a classical approach for computing
singular values. There is extensive convergence analysis on subspace
iteration methods~\cite{golubvanloan,demmelbk,BatheWilson,Template}
and a large literature on accelerated subspace iteration
methods~\cite{saad2}. In general, it is well-suited for fast
computations on modern computers because its main computations are in
terms of matrix-matrix products and QR factorizations that have been
highly optimized for maximum efficiency on modern serial and parallel
architectures~\cite{demmelbk,golubvanloan}. There are two well-known
weaknesses of subspace iteration, however, that limit its practical
use.  On one hand, subspace iteration typically requires very good
separation between the wanted and unwanted singular values for good
convergence. On the other hand, good convergence also often critically
depends on the choice of a good start
matrix~\cite{BatheWilson,Template}.

Another classical class of approximation methods for computing an
approximate SVD are the Krylov subspace methods, such as the Lanczos
algorithm (see, for
example~\cite{CalvettiReichelSorensen,CullumWilloughbyI,KokiopoulouBekasGallopoulos,LehoucqSorensenYang,Saad2011,WuSimon}.)
The computational cost of these methods depends heavily on several
factors, including the start vector, properties of the input matrix
and the need to stabilize the algorithm. One of the most important
part of the Krylov subspace methods, however, is the need to do a
matrix-vector product at each iteration. In contrast to matrix-matrix
products, matrix-vector products perform very poorly on modern
architectures due to the limited data reuse involved in such
operations, In fact, one focus of Krylov subspace research is on effective avoidance of matrix-vector operations in
Krylov subspace methods (see, for
example~\cite{GrimesLewisSimon,Ruhe79}.)

This work focuses on the surprisingly strong performance of
randomized algorithms in delivering highly accurate low-rank
approximations and singular values. To illustrate, we introduce
Algorithm~\ref{Alg:RandSamI}, one of the basic randomized
algorithms (see~\cite{Tropp}.)
\begin{algorithm}\label{Alg:RandSamI}{\bf Basic Randomized Algorithm} \\
\headerrule

\begin{tabular}{ll}
{\bf Input:} & $m \times n$ matrix $A$ with $ m \geq n$,
integers $k>0$ and $n > \ell >k$. \\
{\bf Output:} & a rank-$k$ approximation.  \\
\end{tabular}

\headerrule
\begin{INDENT}
    \begin{enumerate}
   \item Draw a random $n \times \ell $ test matrix $\Omega$. 
 \item     Compute $Y = A \, \Omega$.
 \item    Compute an orthogonal column basis $Q$ for $Y$. 
 \item     Compute $B =  Q^{T} \, A $. 
 \item    Compute $B_k$, the rank-$k$ truncated SVD of $B$. 
 \item     Return $QB_k$.
    \end{enumerate}
\end{INDENT}
\headerrule
\end{algorithm}
\begin{REMARK}\label{Rem:rand}
Throughout this paper, a random matrix, such as $\Omega$ in
Algorithm~\ref{Alg:RandSamI}, is a standard Gaussian matrix,
i.e., its entries are independent standard normal variables 
of zero mean and standard deviation $1$. 

While other random matrices might work equally well, the
choice of the Gaussian matrix provides two unique
advantages: First, the distribution of a standard Gaussian
matrix is rotationally invariant: If $V$ is an orthonormal
matrix, then $V^{T}\Omega$ is itself a standard Gaussian
matrix with the same statistical properties as
$\Omega$~\cite{Tropp}. Second, our analysis is much simplified by the
vast literature on the singular value probability density
functions of the Gaussian matrix.
\end{REMARK}

While Algorithm~\ref{Alg:RandSamI}
looks deceptively simple, its analysis is long,
arduous, and involves very strong doses of
statistics~\cite{Tropp}. The following theorem establishes
an error bound on the accuracy of $QQ^{T} A $ as a low-rank
approximation to $A$. There are similar results in the
Frobenius norm.
\begin{theorem} (Halko, Martinsson, Tropp~\cite[Corollary 10.9]{Tropp})
  \label{Thm:Troppbound2}
The column-orthonormal matrix $Q$ produced by Step 3 in 
Algorithm~\ref{Alg:RandSamI} satisfies
\[{\displaystyle \|\left(I-QQ^{T}\right)A\|_2 \leq \left(1 +
  17\sqrt{1+ \frac{k}{p}}\right) \sigma_{k+1} +
  \frac{8\sqrt{k+p}}{p+1} \sqrt{\sum_{j=k+1}^n \sigma_j^2},
  \quad \mbox{provided that} \quad p = \ell - k \geq 4, }\]
with failure probability at most $6e^{-p}$. 
\end{theorem}

\begin{REMARK}\label{Rem:tropp}
Comparing Theorem~\ref{Thm:Troppbound2} with
Theorem~\ref{Thm:truncSVD}, it is clear that
Algorithm~\ref{Alg:RandSamI} could provide a very good low rank
approximation to $A$ with probability at least $1-6e^{-p}$,
despite its simple operations, provided that $\sigma_{k+1}
\ll \|A\|_2$. While
algorithms~\cite{ChengGimbutasMartinssonRokhlin,DrineasKannanMahoneyII,DrineasMahoneyMuthukrishnan,FriezeKannanVempala,FriezeKannanVempalaII,LibertyAilonSinger,Liberty,WoolfeLibertyRokhlinTygert}
differ in their algorithm design, efficiency, and domain
applicability, they typically share the same advantages of
computational efficiency and approximation accuracy.
\end{REMARK}

Algorithm~\ref{Alg:RandSamI} is the combination of Stages A
and B of the Proto Algorithm in~\cite{Tropp}, where
the truncated SVD is considered separately from low-rank
approximation. In Section~\ref{Sec:Trunc} we will discuss
the pros and cons of SVD truncation vs. no
truncation. Algorithm~\ref{Alg:RandSamI} is a special case
of the randomized subspace iteration method (see Algorithm~\ref{Alg:RandSI}), for
which Halko, Martinsson, Tropp~\cite{Tropp} have developed similar results.

However, while the upper bound in Theorem~\ref{Thm:Troppbound2} can be
very satisfactory for many applications, there may be situations where
singular value approximations are also desirable. In addition, it is
well-known that in practical computations randomized algorithms often
far outperform their error
bounds~\cite{Tropp,MartinssonRokhlinShkolniskyTygert,RokhlinTygert},
whereas the results in ~\cite{Tropp} do not suggest
convergence of the computed rank-$k$ approximation to the truncated
SVD in either Algorithm~\ref{Alg:RandSamI} or the more general
randomized subspace iteration method.

Our entire work is based on novel analysis of the subspace iteration
method, and we consider randomized algorithms within the subspace
iteration framework. This allows us to take advantage of existing
theories and technical machinery in both fields.

Current analysis on randomized algorithms focuses on the
errors in the approximation of $A$ by a low rank matrix,
whereas classical analysis on subspace iteration methods
focuses on the accuracy in the approximate singular
values. Our analysis allows us to obtain both kinds of
results for both of these methods, leading to the stronger
rank-revealing approximations. In terms of randomized
algorithms, our matrix approximation bounds are in general
tighter and can be drastically better than existing ones; in
terms of singular values, our relative convergence lower
bounds can be interpreted as simultaneously convergence
error bounds and rank-revealing lower bounds.

Our analysis has lead us to some interesting conclusions,
all with high probability (more precise statements are in
Sections~\ref{Sec:StatAnal} through~\ref{Sec:CondEst}): 
\begin{itemize}
\item The leading $k$ singular values computed by randomized
algorithms are at least a good fraction of the true ones,
regardless of how the singular values are distributed, and
they converge quickly to the true singular values in case of
rapid singular value decay. In particular, this result
implies that randomized algorithms can also be used as
efficient and reliable condition number estimators.
\item The above results, together with the fact that
randomized algorithms compute low-rank approximations up to
a dimension dependent constant factor from optimal, mean
that these low-rank approximations are in fact
rank-revealing factorizations. In addition, for rapidly
decaying singular values, these approximations can be as
accurate as a truncated SVD.
\item The subspace iteration method in general and the power method in
particular is still slowly convergent without over-sampling in the
start matrix. We present an alternative choice of the start matrix
based on our analysis, and demonstrate its competitiveness.
\end{itemize}

The rest of this paper is organized as follows: In
Section~\ref{Sec:Alg} we discuss subspace iteration methods
and their randomized versions in more detail; in
Section~\ref{Sec:Setup} we list a number of preliminary as
well as key results needed for later analysis; in
Section~\ref{Sec:DetAnal} we derive deterministic lower
bounds on singular values and upper bounds on low-rank
approximations; in Section~\ref{Sec:StatAnal} we provide
both average case and large deviation bounds on singular
values and low-rank approximations; in
Section~\ref{Sec:RankRev} we compare these approximations
with other rank-revealing factorizations; in
Section~\ref{Sec:CondEst} we discuss how randomized
algorithms can be used as efficient and reliable
condition number estimators; in Section~\ref{Sec:Num} we
present supporting numerical experimental results; and in
Section~\ref{Sec:Conc} we draw some conclusions and point
out possible directions for future research.

Much of our analysis has its origin in the analysis of subspace
iteration~\cite{saad2} and randomized algorithms~\cite{Tropp}. It
relies both on linear algebra tools as well as statistical analysis to
do some of the needed heavy lifting to reach our conclusions. 
To limit the length of this paper, we have put the
more detailed parts of the analysis as well as some additional
numerical experimental results in the {\em Supplemental Material},
which is accessible at SIAM's on-line portal.

\section{Algorithms}\label{Sec:Alg}
In this section, we present the main algorithms that are 
discussed in the rest of this paper. We also discuss 
subtle differences between our presentation of randomized
algorithms and that in~\cite{Tropp}.

\subsection{Basic Algorithms} \label{Sec:BasicAlg}
We start with the classical subspace iteration method for
computing the largest few singular values of a given matrix.
\begin{algorithm} \label{Alg:BasicSI}{\bf Basic Subspace Iteration}\\
\headerrule

\begin{tabular}{ll}
{\noindent \bf Input:} & $m \times n$ matrix $A$ with $ n \leq
m$, integers $0 < k \leq \ell < n$,\\
& and $n \times \ell$ start matrix $\Omega$. \\
{\noindent \bf Output:} & a rank-$k$ approximation. \\
\end{tabular}

\headerrule
\begin{INDENT}
    \begin{enumerate}
   \item Compute $Y = \left(AA^{T}\right)^q A \, \Omega$.
   \item Compute an orthogonal column basis $Q$ for $Y$. 
   \item Compute $B =  Q^{T} \, A $. 
   \item Compute $B_k$, the rank-$k$ truncated SVD of $B$. 
   \item Return $QB_k$.
    \end{enumerate}
\end{INDENT}
\headerrule
\end{algorithm}
Given the availability of Lanczos-type algorithms for the
singular value computations, the classical subspace
iteration method is not widely used in practice except when
$k \ll n$.  We present it here for later comparisons with
its randomized version. We ignore the vast literature of
accelerated subspace iteration methods (see, for example~\cite{saad2}) in this paper since
our main goal here is to analyze the convergence behavior of
subspace iteration method with and without randomized start
matrix $\Omega$.

We have presented Algorithm~\ref{Alg:BasicSI} in an
over-simplified form above to convey the basic ideas
involved. In practice, the computation of $Y$ would be prone
to round-off errors. For better numerical accuracy,
Algorithm~\ref{Alg:B} in the Appendix should be
used numerically to compute the $Q$ matrix in 
Algorithm~\ref{Alg:BasicSI}. In practical computations,
however, Algorithm~\ref{Alg:B} is often performed once every
few iterations, to balance efficiency and numerical
stability (see Saad~\cite{saad2}.) In the rest of
Section~\ref{Sec:Alg}, any QR factorization of the matrix $Y
= \left(A A^{T}\right)^q A \Omega$ should be computed
numerically through periodic use of Algorithm~\ref{Alg:B}.

While there is little direct analysis of subspace iteration methods
for singular values (Algorithm~\ref{Sec:BasicAlg}) in the literature,
one can generalize results of subspace iteration methods for
symmetric matrices to the singular value case in a straightforward
fashion. The symmetric matrix version of Theorem~\ref{Thm:BatheWilson} can be
found in~\cite{BatheWilson}.  
\begin{theorem} (Bathe and Wilson) \label{Thm:BatheWilson}
Assume that
  Algorithm~\ref{Alg:BasicSI} converges as $q \rightarrow \infty$. Then
\[{\displaystyle |\sigma_j- \sigma_j(Q^T B_k)| \leq
O\left(\left(\frac{\sigma_{\ell+1}}{\sigma_k}\right)^{2q+1}\right)
.} \]
\end{theorem}
Thus convergence is governed by the ratio
${\displaystyle\frac{\sigma_{\ell+1}}{\sigma_k}}$. The
per-iteration cost of Algorithm~\ref{Alg:BasicSI} depends
linearly on $\ell \geq k$. A choice $\ell > k$ can be
economical if the more rapid convergence obtained through
the ratio ${\displaystyle\frac{\sigma_{\ell+1}}{\sigma_k}}$
can more than offset the extra cost per iteration. Another important issue with Algorithm~\ref{Alg:BasicSI} is the
constant hidden in the $O$ notation.
This constant can be exceedingly large for the unfortunate
choices of $\Omega$. In fact, an $\Omega$ matrix that is
almost orthogonal to any leading singular vectors will lead
to large number of iterations. Both issues will be made clearer with our relative convergence theory for
Algorithm~\ref{Alg:BasicSI} in Theorem~\ref{Thm:SVdetBound}. 

A special case of Algorithm~\ref{Alg:BasicSI} is when $k =
\ell = 1$. This is the classical power method for computing
the 2-norm of a given matrix. This method, along with its
randomized version, is included in Appendix~A 
for later discussion in our numerical experiments (see
Section~\ref{Sec:Num}.) The power method has the same
convergence properties of Algorithm~\ref{Alg:BasicSI}. More generally,
the subspace iteration method is typically run with $k = \ell$. 

\subsection{Randomized Algorithms} \label{Sec:Random}
In order to enhance the convergence of Algorithm~\ref{Sec:BasicAlg}
in the absence of any useful information about the leading singular
vectors, a sensible approach is to replace the deterministic
start matrix with a random one, leading to
\begin{algorithm} \label{Alg:RandSI}{\bf Randomized Subspace Iteration}\\ 
\headerrule

\begin{tabular}{ll}
{\noindent \bf Input:} & $m \times n$ matrix $A$ with $ n \leq
m$, integers $0 < k \leq \ell $,\\
{\noindent \bf Output:} & a rank-$k$ approximation. \\
\end{tabular}

\headerrule
\begin{INDENT}
    \begin{enumerate}
   \item Draw a random $n \times \ell$ start matrix $\Omega$. 
   \item Compute a rank-$k$ approximation with Algorithm~\ref{Alg:BasicSI}. 
    \end{enumerate}
\end{INDENT}
\headerrule
\end{algorithm}

\begin{REMARK}\label{Rem:basicrand} Since Algorithm~\ref{Alg:RandSI} is the special case of
Algorithm~\ref{Alg:BasicSI} with $\Omega$ being chosen as
random, all our results for Algorithm~\ref{Alg:BasicSI}
equally hold for Algorithm~\ref{Alg:RandSI}.
\end{REMARK}

The only difference between Algorithm~\ref{Alg:BasicSI} and
Algorithm~\ref{Alg:RandSI} is in the choice of $\Omega$,
yet this difference will lead to drastically different
convergence behavior. One of the main purposes of this paper
is to show that the slow or non-convergence of
Algorithm~\ref{Alg:BasicSI} due to bad choice of $\Omega$
vanishes with near certainty in
Algorithm~\ref{Alg:RandSI}.  In particular, 
a single iteration ($q = 0$ in Algorithm~\ref{Alg:RandSI}) 
in the randomized subspace iteration method is often
sufficient to return good enough singular values and
low-rank approximations (Section~\ref{Sec:StatAnal}).

Our analysis of deterministic and randomized subspace iteration method
was in large part motivated by the analysis and discussion of
randomized algorithms in~\cite{Tropp}. We have chosen to present the
algorithms in Section~\ref{Sec:Alg} in forms that are not identical to
those in~\cite{Tropp} for ease of stating our results in
Sections~\ref{Sec:DetAnal} through~\ref{Sec:Num}. Versions of
Algorithm~\ref{Alg:RandSI} have also appeared in~\cite{XiangZou} for
solving large-scale discrete inverse problems.

\subsection{To Truncate or not to Truncate}\label{Sec:Trunc} 

The randomized algorithms in Section~\ref{Sec:Alg} are presented in a
slight different form than those in~\cite{Tropp}. One key difference
is in the step of SVD truncation, which is considered an optional
postprocessing step there. In this section, we discuss the pros and
cons of SVD truncation. We start with the following simple lemma, 
versions of which appear in~\cite{BoutsidisDrineasMahoney,DrineasMahoneyMuthukrishnan2008,Tropp}.

\begin{lemma} \label{Lem:Trunc} Given an $m \times \ell$ matrix with orthonormal columns $Q$, with $\ell \leq n$, then
  for any $\ell \times n$ matrix $B$, 
\[ \|A - Q \left(Q^TA\right) \|_2 \leq \|A - Q B \|_2 
\quad \mbox{and} \quad 
\|A - Q \left(Q^TA\right) \|_F \leq \|A - Q B \|_F . \]
\end{lemma} 

Lemma~\ref{Lem:Trunc} makes it obvious that any SVD
truncation of $Q^T A$ will only result in a less accurate
approximation in the 2-norm and Frobenius norm. This is
strong motivation for no SVD truncation. The SVD truncation
of $Q^T A$ also involves the computation of the SVD of $Q^T
A$ in some form, which also results in extra computation.

On the other hand, since singular values of $Q^TA$ approximate their
corresponding singular values in $A$ at different rates, some singular
values of $Q^TA$ may be poor approximations of those of $A$, and
$QQ^TA$ need not be a good rank-$\ell$ approximation to $A$,
either. In contrast, for the right choices of $k$, the rank-$k$
truncated SVD of $Q^TA$ can contain $k$ excellent approximate singular values
and result in a good rank-$k$ approximation to $A$ as
well. So the choice of whether to truncate the SVD of $Q^TA$ depends on
practical considerations of computational efficiency and demands on
quality of singular value and low-rank approximations. This paper
focuses on a rank-$k$ approximations obtained from truncated SVD of
$Q^TA$.

\section{Setup}\label{Sec:Setup}
In this section we build some of the technical machinery needed for
our heavy analysis later on. We start by reciting two well-known
results in matrix analysis, and then develop a number of theoretical
tools that outline our approach in the low-rank approximation
analysis. Some of these results may be of interest in their own
right. For any matrix $X$, we use $\sigma_j(X)$ to denote its $j$-th
largest singular value.

The Cauchy interlacing theorem shows the limitations of any
approximation with an orthogonal projection. 
\begin{theorem}(Golub and van Loan~\cite[p. 411]{golubvanloan})
  \label{Thm:Cauchy} 
 Let $A$ be an $m\times n$ matrix and $Q$ be a matrix with
 orthonormal columns. Then $\sigma_j(A) \geq \sigma_j(Q^T A)$
 for $1 \leq j \leq \min(m,n)$. 
\end{theorem}

\begin{REMARK}\label{Rem:Interlace} A direct consequence of
 Theorem~\ref{Thm:Cauchy} is that $\sigma_j(A) \geq
 \sigma_j(\widehat{A})$, where $\widehat{A}$ is any
 submatrix of $A$.
\end{REMARK}

Weyl's monotonicity theorem relates singular values of
matrices $X$ and $Y$ to those of $X+Y$. 
\begin{theorem}\label{Thm:Weyl} (Weyl's monotonicity
  theorem~\cite[Thm. 3.3.16]{HORNJOHNSON})
Let $X$ and $Y$ be $m \times n$ matrices with $m \geq n$. Then
\[ \sigma_{i+j-1}\left(X+Y\right) \leq
\sigma_i(X)+\sigma_j(Y) \quad \mbox{for all $i, j \geq 1$
  such that $i+j-1\leq n$.} \]
\end{theorem}

The Hoffman-Wielandt theorem bounds the errors in the 
differences between the singular values of $X$ and those of $Y$ in terms of $\|X-Y\|_F$.
\begin{theorem}(Hoffman and Wielandt~\cite{HW}) \label{Thm:HW} 
Let $X$ and $Y$ be $m \times n$ matrices with $m \geq n$. Then
\[{\displaystyle \sqrt{\sum_{j=1}^n \left|\sigma_j\left(X\right) - \sigma_j\left(Y\right)\right|^2} \leq \|X-Y\|_F . } \]
\end{theorem}

Below we develop a number of theoretical results that will
form the basis for our later analysis on low-rank
approximations. Theorem~\ref{Thm:Fto2} below is of potentially broad independent interest. Let $B$ be a
rank-$k$ approximation to $A$. Theorem~\ref{Thm:Fto2} below
relates the approximation error in the Frobenius norm to that in
the 2-norm as well as the approximation errors in the
leading $k$ singular values. It will be called the {\em
Reverse Eckart and Young Theorem} due to its complimentary
nature with Theorem~\ref{Thm:truncSVD} in the Frobenius
norm. 

\begin{theorem}\label{Thm:Fto2} (Reverse Eckart and Young) 
Given any $m \times n$ matrix $A$, and let $B$ be a matrix with rank at most $k$ such that
\begin{equation}\label{Eqn:Fbound}
\| A - B \|_F \leq {\displaystyle \sqrt{\eta^2 +
\sum_{j=k+1}^n \sigma^2_{j} }}
\end{equation}
for some $\eta \geq 0$. Then we must have
\begin{eqnarray}
\| A - B \|_2 & \leq & {\displaystyle \sqrt{{\eta^2} + \sigma_{k+1}^2 
    },} \label{Eqn:F2bound} \\
\sqrt{\sum_{j=1}^k \left(\sigma_{j}- \sigma_{j}(B)\right)^2}
    & \leq &\eta. \label{Eqn:F2Hoffman} 
\end{eqnarray}
\end{theorem}
\begin{REMARK}\label{Rem:quad}
Notice that 
\[{\displaystyle \sqrt{{\eta^2}+ \sigma_{k+1}^2 } =
  \sigma_{k+1} + \frac{\eta^2}{\sqrt{{\eta^2}+\sigma_{k+1}^2 
      } + \sigma_{k+1}} }.\]
Equation~(\ref{Eqn:F2bound}) can be simplified to 
\begin{equation}\label{Eqn:eta}
\| A - B \|_2  \leq \sigma_{k+1} + \eta 
\end{equation}
when $\eta$ is larger than or close to $\sigma_{k+1}$. On
the other hand, if $\eta \ll \sigma_{k+1}$, then 
equation~(\ref{Eqn:F2bound}) simplifies to
\[\| A - B \|_2  \leq  \sigma_{k+1} +
\frac{\eta^2}{2\sigma_{k+1}}, \] 
where the last ratio can be
much smaller than $\eta$, implying a much better rank-$k$
approximation in $B$. Similar comments apply to
equation~(\ref{Eqn:Fbound}). This interesting feature of
Theorem~\ref{Thm:Fto2} is one of the reasons why our
eventual 2-norm and Frobenius norm upper bounds are much
better than those in Theorem~\ref{Thm:Troppbound2} in the
event that $\eta \ll \sigma_{k+1}$. This also has made our
proofs in Appendix B somewhat involved in places.  
\end{REMARK}

\begin{REMARK} Equation~(\ref{Eqn:F2Hoffman}) asserts that a small $\eta$ in
equation~(\ref{Eqn:Fbound}) necessarily means good
approximations to all the $k$ leading singular values of
$A$. In particular, $\eta = 0$ means the leading $k$
singular values of $A$ and $B$ must be the same. However,
our singular value analysis will not be based on
Equation~(\ref{Eqn:F2Hoffman}), as our approach in
Section~\ref{Sec:DetAnal} provides us with much
better results.
\end{REMARK}

{\noindent \bf Proof of Theorem~\ref{Thm:Fto2}:} Write $A =
\left(A-B\right) + B$. It follows from
Theorem~\ref{Thm:Weyl} that for any $1 \leq i \leq n-k$:
\[{\displaystyle  \sigma_{i+k}(A) \leq \sigma_{i}(A-B) +
  \sigma_{k+1}(B) = \sigma_{i}(A-B), } \]
since $B$ is a rank-$k$ matrix. It follows that
\[{\displaystyle \|A-B\|_F^2 = \sum_{i=1}^n \sigma_i^2(A-B)
  \geq \sigma_1^2(A-B)+ \sum_{i=2}^{n-k} \sigma_i^2(A-B)
  \geq \sigma_1^2(A-B)+ \sum_{i=2}^{n-k} \sigma_{i+k}^2.} \]
Plugging this into equation~(\ref{Eqn:Fbound})
  yields~(\ref{Eqn:F2bound}). 

As to equation~(\ref{Eqn:F2Hoffman}), we observe that the $(k+1)-st$ through the
last singular values of $B$ are all zero, given that $B$ has rank
$k$. Hence the result trivially follows from Theorem~\ref{Thm:HW},
\[\sum_{j=1}^k \left(\sigma_j - \sigma_j(B)\right)^2 + 
\sum_{j=k+1}^n \sigma_j^2 \leq\| A - B \|_F^2 \leq
\eta^2 +  \sum_{j=k+1}^n \sigma_j^2. \quad \mbox{\Q}\]

Our next theorem is a generalization of Theorem~\ref{Thm:truncSVD}. 
\begin{theorem} \label{Thm:SVD2}
Let $Q$ be an $m \times \ell$ matrix with orthonormal columns, let $1 \leq k
\leq \ell$, and let $B_k$ be the rank-$k$ truncated SVD of $Q^{T}
A$. Then $B_k$ is an optimal solution to the following problem
\begin{equation}\label{Eqn:SVD2}
\min_{{\bf rank}(B) \leq k,  B \;\; \mbox{$B$ is $\ell \times n$}} \| A - Q B \|_F = \| A - Q B_k \|_F .
\end{equation}
In addition, we also have
\begin{equation}\label{Eqn:Bk2AkF}
\| A - Q B_k \|_F^2 \leq  \| \left(I- Q Q^{T}\right)A_k \|_F^2  + \sum_{j=k+1}^n \sigma^2_{j}. 
\end{equation}
\end{theorem}

\begin{REMARK}\label{Rem:Fvs2} Problem~(\ref{Eqn:SVD2}) in Theorem~\ref{Thm:SVD2} is a type of
restricted SVD problem. Oddly enough, this problem becomes much harder
to solve for the 2-norm. In fact, $B_k$ might not even be the solution
to the corresponding restricted SVD problem in
2-norm. Combining Theorems~\ref{Thm:Fto2}
and~\ref{Thm:SVD2}, we obtain
\begin{equation}\label{Eqn:Bk2Ak2}
\| A - Q B_k \|_2^2 \leq  \| \left(I- Q Q^{T}\right)A_k \|_F^2  + \sigma^2_{k+1}. 
\end{equation}
Our low-rank approximation analysis in the $2$-norm will be
based on equation~(\ref{Eqn:Bk2Ak2}). While this is
sufficient, it also makes our $2$-norm results perhaps
weaker than they should be due to the mixture of the
$2$-norm and the Frobenius norm.
\end{REMARK}

By Theorem~\ref{Thm:truncSVD}, $A_k$ is the best Frobenius
norm approximation to $A$, whereas by Theorem~\ref{Thm:SVD2}
$QB_k$ is the best restricted Frobenius norm approximation
to $A$. This leads to the following interesting consequence
\begin{equation}\label{A2B2A}
\|A-A_k \|_F \leq \|A-QB_k \|_F \leq \|A-QQ^{T} A_k \|_F . 
\end{equation}
Thus we can expect $QB_k$ to also be an excellent rank-$k$
approximation to $A$ as long as $Q$ points to the principle
singular vector directions.

{\noindent \bf Proof of Theorem~\ref{Thm:SVD2}:} We first rewrite 
\[\| A - Q B \|_F^2 = \| \left(I - Q Q^T\right) A + Q \left(Q^TA - B\right) \|_F^2 
= \| \left(I - Q Q^T\right) A\|_F^2 + \|\left(Q^TA - B\right)
\|_F^2. \] 
Result~(\ref{Eqn:SVD2}) is now an immediate consequence
of Theorem~\ref{Thm:truncSVD}. To prove~(\ref{Eqn:Bk2AkF}), we observe that 
\begin{eqnarray*}
\|A - Q Q^T A_k\|_F^2 & = & {\bf trace}\left(\left(A - Q Q^T
A_k\right)^T\left(A - Q Q^T A_k\right)\right)  \\
&=& {\bf trace}\left(\left(A - A_k + A_k - Q Q^T
A_k\right)^T\left(A - A_k + A_k - Q Q^T A_k\right)\right) \\
&=& \|A-A_k\|_F^2 + \|A_k-QQ^T A_k\|_F^2 
+ 2 {\bf trace}\left(\left(A - A_k\right)^T\left(A_k - Q Q^T
A_k\right)\right)  \\
&=& \sum_{j=k+1}^n \sigma_j^2 + \|A_k-QQ^T A_k\|_F^2 
+ 2 {\bf trace}\left(\left(\left(I - Q
Q^T\right)A_k\right)\left(A - A_k\right)^T\right) . 
\end{eqnarray*}
The third term in the last equation is zero because $A_k\left(A -
A_k\right)^T = 0.$ Combining this last relation with
equation~(\ref{A2B2A}) gives us relation~(\ref{Eqn:Bk2AkF}). \Q

\section{Deterministic Analysis}\label{Sec:DetAnal} In this
section we perform deterministic convergence analysis on
Algorithm~\ref{Alg:BasicSI}. Theorem~\ref{Thm:SVdetBound} is
a relative convergence lower bound, and
Theorem~\ref{Thm:DetBound2} is an upper bound on the matrix
approximation error. Both appear to be new for subspace
iteration. Our approach, while quite novel, was motivated in
part by the analysis of subspace iteration methods by
Saad~\cite{saad2} and randomized algorithms
in~\cite{Tropp}. Since Algorithm~\ref{Alg:RandSamI} is a
special case of Algorithm~\ref{Alg:RandSI} with $q = 0$,
which in turn is a special case of
Algorithm~\ref{Alg:BasicSI} with an initial random matrix,
our analysis applies to them as well and will form the basis
for additional probabilistic analysis in
Section~\ref{Sec:StatAnal}.

\subsection{A Special Orthonormal Basis}\label{Sec:DetPrelim}
We begin by noticing that the output $QB_k$ in
Algorithm~\ref{Alg:BasicSI} is also the rank-$k$ truncated
SVD of the matrix $QQ^TA$, due to the fact that $Q$ is
column orthonormal. In fact, columns of $Q$ are nothing but
an orthonormal basis for the column space of matrix
$\left(AA^{T}\right)^q A \, \Omega$. This is the reason why
Algorithm~\ref{Alg:BasicSI} is called subspace
iteration. Lemma~\ref{Lem:Basis} below shows how to obtain
alternative orthonormal bases for the same column space. We
omit the proof.
\begin{lemma}\label{Lem:Basis} In the notation of
  Algorithm~\ref{Alg:BasicSI}, assume that $X$ is a
  non-singular $\ell \times \ell$ matrix and that $\Omega$
  has full column rank. Let $\widehat{Q} \widehat{R}$ be the
  QR factorization of the 
  matrix $\left(AA^T\right)^q A \Omega X$, then 
\[ Q Q^T  = \widehat{Q} \widehat{Q}^T. \]
\end{lemma} 

Since 
\[{\displaystyle \left(AA^T\right)^q A \Omega  = U \Sigma^{2q+1}
  V^T \Omega ,}\]
define and partition 
\begin{equation}\label{Eqn:partition}
{\displaystyle \widehat{\Omega} \stackrel{def}= V^T \Omega = \begin{array}{rl}
\begin{array}{cc}
\ell - p & \left\{\right. \cr
n-\ell + p & \left\{\right. \end{array} 
& \hspace{-0.15in} \left(\begin{array}{c}
  \widehat{\Omega}_1  \cr \widehat{\Omega}_2\end{array}
  \right) \end{array},} 
\end{equation}
where $0 \leq p \leq \ell - k$. The introduction of the additional parameter $p$
is to balance the need for oversampling for reliability (see Theorem~\ref{Thm:Troppbound2})
and oversampling for faster convergence (see Theorem~\ref{Thm:BatheWilson}). We also partition 
$\Sigma = \diag\left(\Sigma_1,\Sigma_2,\Sigma_3\right)$,
where $\Sigma_1$, $\Sigma_2$, and $\Sigma_3$ are $k \times
k$, $(\ell-p -k)\times (\ell-p-k)$, and $(n-\ell+p)\times (n-\ell+p)$.
This partition allows us to further write
\begin{equation}\label{Eqn:AO}
{\displaystyle \left(AA^T\right)^q A \Omega  = U 
\left(\begin{array}{c}
\left(\begin{array}{cc}
\Sigma_1 & \cr
& \Sigma_2 \end{array}\right)^{2q+1} 
\widehat{\Omega}_1
\cr
 \cr \Sigma_3^{2q+1}\widehat{\Omega}_2\end{array}
\right).} 
\end{equation}
The matrix ${\displaystyle \widehat{\Omega}_1}$ has at least
as many columns as rows. Assume it is of full row rank so
that its pseudo-inverse satisfies
\[{\displaystyle  \widehat{\Omega}_1 
\widehat{\Omega}_1^{\dagger} = I.} \]

Below we present a special choice of $X$ that will reveal 
the manner in which convergence to singular values and
low-rank approximations takes place. Ideally, such an $X$
would orient the first $k$ columns of ${\displaystyle U 
\left(\begin{array}{c}
\left(\begin{array}{cc}
\Sigma_1 & \cr
& \Sigma_2 \end{array}\right)^{2q+1} 
\widehat{\Omega}_1
\cr
 \cr \Sigma_3^{2q+1}\widehat{\Omega}_2\end{array}
\right)X}$ in the directions of the leading $k$ singular
vectors in $U$. We choose
\begin{equation}\label{Eqn:X}
{\displaystyle  X = \left( 
\widehat{\Omega}_1^{\dagger} \left(\begin{array}{cc}
\Sigma_1& \cr 
& \Sigma_2 \end{array}
  \right)^{-(2q+1)}, \quad \quad \widehat{X}
 \right)},  
\end{equation}
where the $\ell \times p$ matrix $\widehat{X}$ is chosen so that 
$X$ is non-singular and $\widehat{\Omega}_1 \widehat{X} =
0$. Recalling equation~(\ref{Eqn:AO}), this choice of $X$ allows
  us to write
\begin{equation}\label{Eqn:OX}
{\displaystyle \left(AA^T\right)^q A \Omega  X = 
U \left(\begin{array}{ccc}
I & 0 & 0 \cr
0 & I & 0 \cr
H_1 & H_2 & H_3\end{array}
  \right) ,} 
\end{equation}
where
\[{\displaystyle 
H_1 = \Sigma_3^{2q+1} \widehat{\Omega}_2
\widehat{\Omega}_1^{\dagger} \left(\begin{array}{c}
\Sigma_1^{-(2q+1)} 
\cr
0 \end{array}
  \right), \quad 
H_2 = \Sigma_3^{2q+1} \widehat{\Omega}_2 \widehat{\Omega}_1^{\dagger} \left(\begin{array}{c}
0
\cr
\Sigma_2^{-(2q+1)} 
\end{array}
  \right), \quad 
H_3 = \Sigma_3^{2q+1} \widehat{\Omega}_2 \widehat{X}.} \]

Notice that we have created a ``gap'' in $H_1$: the largest
singular value in $\Sigma_3$ is $\sigma_{\ell-p+1}$, which
is potentially much smaller than $\sigma_k$, the smallest
singular value in $\Sigma_1$. We can expect $H_1$ to
converge to $0$ rather quickly when $q \rightarrow \infty$, if $\sigma_{\ell-p+1} \ll
\sigma_{k}$ and if the matrix ${\displaystyle
\widehat{\Omega}_1^{\dagger}}$ is not too large in norm. Our
convergence analysis of Algorithms~\ref{Alg:BasicSI}
and~\ref{Alg:RandSI} will mainly involve deriving upper
bounds on various functions related to $H_1$. Our gap
disappears when we choose $p = \ell - k$, in which case our
results in Section~\ref{Sec:StatAnal} will be more in line
with Theorem~\ref{Thm:Troppbound2}. 

By equation~(\ref{Eqn:OX}), the QR factorization of
$\left(AA^T\right)^q A \Omega X$ can now be written in the
following $3\times 3$ partition:
\begin{equation}\label{Eqn:QR}
{\displaystyle U \left(\begin{array}{ccc}
I & 0 & 0 \cr
0 & I & 0 \cr
H_1 & H_2 & H_3\end{array}
  \right) = \widehat{Q} \widehat{R} 
= \left(\begin{array}{ccc}
\widehat{Q}_1 & \widehat{Q}_2 & \widehat{Q}_3\end{array}
  \right) \left(\begin{array}{ccc}
\widehat{R}_{11}  & \widehat{R}_{12}  & \widehat{R}_{13}  \cr
 & \widehat{R}_{22}  & \widehat{R}_{23} \cr
 & & \widehat{R}_{33}\end{array}   \right). } 
\end{equation}
We will use this representation to derive convergence upper
bounds for singular value and rank-k approximations. In
particular, we will make use of the fact that the above QR
factorization also embeds another one
\begin{equation}\label{Eqn:Q1}
{\displaystyle U \left(\begin{array}{c}
I  \cr
0  \cr
H_1 \end{array}
  \right) = \widehat{Q}_1  \widehat{R}_{11} . } 
\end{equation}

We are now ready to derive a lower bound on $\sigma_k(B_k)$.
\begin{lemma}\label{Lem:Bk} Let $H_1$ be defined in
  equation~(\ref{Eqn:OX}), and assume that the matrix
  $\widehat{\Omega}_1$ has full row rank, then the matrix
  $B_k$ computed in Algorithm~\ref{Alg:BasicSI} must satisfy 
\begin{equation}\label{Eqn:RSigma}
{\displaystyle \sigma_k\left(B_k\right) \geq
  \frac{\sigma_k}{\sqrt{1 + \left\|H_{1}\right\|_2^2}}} .
\end{equation}
\end{lemma} 

\begin{REMARK} It might seem more
  intuitive in equation~(\ref{Eqn:X}) to choose $X = \left(X_1 \quad X_2\right)$ where
  $X_1$ solves the following least squares problem
\[ {\displaystyle  \min_{X_1} \left\| 
\left(\begin{array}{c}
\left(\begin{array}{cc}
\Sigma_1 & \cr
& \Sigma_2 \end{array}\right)^{2q+1} 
\widehat{\Omega}_1
\cr
 \cr \Sigma_3^{2q+1}\widehat{\Omega}_2\end{array}
\right)X_1 - \left(\begin{array}{c}
I \cr
0 \cr
0 \end{array}\right)\right\|_2.} \]
Our choice of $X$ seems as
effective and allows simpler analysis.
\end{REMARK}

{\noindent \bf Proof of Lemma~\ref{Lem:Bk}:} We note by Lemma~\ref{Lem:Basis} that 
\begin{equation}\label{Eqn:tmp1}
{\displaystyle Q Q^T A = \widehat{Q} \widehat{Q}^T A 
= \widehat{Q} \left(\begin{array}{c|c} 
\widehat{Q}_1^T U \left(\begin{array}{c}
\Sigma_1  \cr
0  \cr
0 \end{array}
  \right) & \widehat{Q}_1^T U\left(\begin{array}{cc}
0 &0    \cr
\Sigma_2 &   \cr
& \Sigma_3 \end{array}
  \right) \cr
& \cr
 \hline
& \cr
\left(\begin{array}{c}
\widehat{Q}^T_2  \cr
\widehat{Q}^T_3  \cr
 \end{array}
  \right) U\left(\begin{array}{c}
\Sigma_1  \cr
0  \cr
0 \end{array}
  \right) & \left(\begin{array}{c}
\widehat{Q}^T_2  \cr
\widehat{Q}^T_3  \cr
 \end{array}
  \right) U \left(\begin{array}{cc}
0 &0    \cr
\Sigma_2 &   \cr
& \Sigma_3 \end{array}
  \right) \end{array}
  \right) V^T.}
\end{equation}
From equations~(\ref{Eqn:tmp1}) and~(\ref{Eqn:Q1}), we see that the matrix 
\[{\displaystyle  \widehat{Q}_1^T U \left(\begin{array}{c}
\Sigma_1  \cr
0  \cr
0 \end{array}
  \right) = \left(U \left(\begin{array}{c}
I  \cr
0  \cr
H_1 \end{array}
  \right)\widehat{R}_{11}^{-1}\right)^T  U \left(\begin{array}{c}
\Sigma_1  \cr
0  \cr
0 \end{array}
  \right) = \widehat{R}_{11}^{-T}\Sigma_1} \]
is simply a submatrix of the middle matrix 
on the right hand side of
equation~(\ref{Eqn:tmp1}). By Remark~\ref{Rem:Interlace}, it follows immediately that 
\[{\displaystyle \sigma_k\left(B_k\right) =
\sigma_k\left(\widehat{Q}\widehat{Q}^TA\right)  \geq
\sigma_k\left(\widehat{R}_{11}^{-T}\Sigma_1\right) } . \]
On the other hand, we also have
\[{\displaystyle \sigma_k =
  \sigma_k\left(\widehat{R}_{11}^T\left(\widehat{R}_{11}^{-T}\Sigma_1\right)\right)
  \leq \left\|\widehat{R}_{11}^T\right\|_2
  \sigma_k\left(\widehat{R}_{11}^{-T}\Sigma_1\right).} \]
Combining these two relations, and together with the fact
  that ${\displaystyle \left\|\widehat{R}_{11}^T\right\|_2 =
  \sqrt{1 + \left\|H_{1}\right\|_2^2}}$,  we
  obtain~(\ref{Eqn:RSigma}).  \Q   
\subsection{Deterministic Bounds} \label{Sec:DetMain}
In this section we develop the analysis in Section~\ref{Sec:DetPrelim}
into deterministic lower bounds for singular values and upper bounds
for rank-k approximations. 

Since the interlacing theorem~\ref{Thm:Cauchy} asserts an upper bound
$\sigma_k\left(B_k\right) \leq \sigma_k$, equation~(\ref{Eqn:RSigma})
provides a nice lower bound on $\sigma_k\left(B_k\right)$. These
bounds mean that $\sigma_k\left(B_k\right)$ is a good approximation to
$\sigma_k$ as long as $\left\|H_{1}\right\|_2$ is small. This
consideration is formalized in the theorem below.  
\begin{theorem}
\label{Thm:SVdetBound} Let $A = U \Sigma V^{T}$ be the SVD of $A$, let
$0 \leq p \leq \ell - k$, and let $V^T \Omega$ be partitioned in
equation~(\ref{Eqn:partition}). Assume that the matrix
$\widehat{\Omega}_1$ has full row rank, then
Algorithm~\ref{Alg:BasicSI} must satisfy for $j = 1, \cdots, k$: 
\[{\displaystyle \sigma_j \geq \sigma_j\left(B_k\right)
  \geq \frac{\sigma_j}{\displaystyle \sqrt{1 +
    \left\|\widehat{\Omega}_2\right\|_2^2 \left\|\widehat{\Omega}_1^{\dagger}\right\|_2^2 
    \left(\frac{\sigma_{\ell - p +
    1}}{\sigma_j}\right)^{4q+2}}}. } \]
\end{theorem}

{\noindent \bf Proof of Theorem~\ref{Thm:SVdetBound}:} By the definition of the matrix $H_1$ in
  equation~(\ref{Eqn:OX}), it is straightforward to get
\begin{equation}\label{Eqn:H1}
{\displaystyle \| H_1\|_2 \leq \left\|\widehat{\Omega}_2\right\|_2\left\|\widehat{\Omega}_1^{\dagger}\right\|_2 
    \left(\frac{\sigma_{\ell - p +
    1}}{\sigma_k}\right)^{2q+1}.} 
\end{equation}
This, together with lower bound~(\ref{Eqn:RSigma}), gives
the result in Theorem~\ref{Thm:SVdetBound} for $j =
k$. To prove Theorem~\ref{Thm:SVdetBound} for any $1 \leq j
< k$, we observe that since $\sigma_j\left(B_k\right)
    =\sigma_j\left(B_j\right)$, all that is needed is to 
repeat all previous arguments for a rank $j$ truncated
    SVD. \Q

\begin{REMARK}\label{Rem:SI} Theorem~\ref{Thm:SVdetBound}
makes explicit the two key factors governing the convergence
of Algorithm~\ref{Alg:BasicSI}. On one hand, we can expect
fast convergence for $\sigma_j(B_k)$ if $\sigma_{\ell-p+1} \ll
\sigma_j$. On the other hand, an unfortunate choice of
$\Omega$ could potentially make $\|\Omega_1^{\dagger}\|_2$
very large, leading to slow converge even if the singular
values do decay quickly. The main effect of randomization in
Algorithm~\ref{Alg:RandSI} is to ensure a reasonably sized
$\|\widehat{\Omega}_2\|_2\|\widehat{\Omega}_1^{\dagger}\|_2$
with near certainty. See Theorem~\ref{Thm:HDBoundsIII} for a
precise statement and more details.
\end{REMARK}

Now we consider rank-k approximation upper bounds. Toward
this end and considering Theorem~\ref{Thm:SVD2}, we would
like to start with an upper bound on $ \| \left(I- Q
Q^{T}\right)A_k \|_F$. By Lemma~\ref{Lem:Basis} and
equation~(\ref{Eqn:QR}), we have 
\[ {\displaystyle \left\| \left(I- Q
Q^{T}\right)A_k \right\|_F}  =  {\displaystyle 
\left\| \left(I- \widehat{Q} \widehat{Q}^T\right)A_k \right\|_F \leq 
\left\| \left(I- \widehat{Q}_1 \widehat{Q}_1^T\right)A_k \right\|_F . }
\]
Since $A_k = U \diag\left(\Sigma_1,  0,  0\right)V^T$, and
since $\widehat{Q}_1 = {\displaystyle U \left(\begin{array}{c}
I  \cr
0  \cr
H_1 \end{array}
  \right) \widehat{R}_{11}^{-1}}$ according to
equation~(\ref{Eqn:Q1}), the above right hand side becomes
\[ {\displaystyle \left\| \left(I- \widehat{Q}_1
  \widehat{Q}_1^T\right)A_k \right\|_F = \left\| \left(I- \left(\begin{array}{c}
I  \cr
0  \cr
H_1 \end{array}
  \right)\left(I + H_1^T H_1\right)^{-1} 
\left(\begin{array}{c}
I  \cr
0  \cr
H_1 \end{array}
  \right)^T\right) \left(\begin{array}{c}
\Sigma_1  \cr
0  \cr
0 \end{array}
  \right) \right\|_F  ,}\]
 where we have used the fact that (see~(\ref{Eqn:Q1}))
\[ {\displaystyle  \widehat{R}_{11}^{-1}
  \widehat{R}_{11}^{-T} =
  \left(\widehat{R}_{11}^T\widehat{R}_{11}\right)^{-1} =
  \left(I + H_1^T H_1\right)^{-1}. }\] 
Continuing, 
\begin{eqnarray}
\left\| \left(I- \widehat{Q}_1 \widehat{Q}_1^T\right)A_k
\right\|_F & = &  {\displaystyle \left\| \left(\begin{array}{ccc}
I - \left(I + H_1^T H_1\right)^{-1} & 0 & - \left(I + H_1^T
H_1\right)^{-1} H_1^T \cr
0 & I & 0 \cr
-H_1 \left(I + H_1^T H_1\right)^{-1} & 0 & 
I - H_1 \left(I + H_1^T H_1\right)^{-1} H_1^T
\end{array}
  \right) \left(\begin{array}{c}
\Sigma_1  \cr
0  \cr
0 \end{array}
  \right) \right\|_F  } \nonumber \\
& = & {\displaystyle \left\| \left(\begin{array}{c}
H_1^T \left(I + H_1 H_1^T\right)^{-1} H_1\cr
 \cr
-\left(I + H_1 H_1^T\right)^{-1} H_1
\end{array}
  \right) \Sigma_1  \right\|_F  } \nonumber \\
& = & {\displaystyle \sqrt{{\bf trace}\left(\Sigma_1H_1^T
\left(I + H_1 H_1^T\right)^{-1}
H_1\Sigma_1\right)}}\label{Eqn:trace}
\end{eqnarray}
We are now ready to prove 
\begin{theorem} \label{Thm:DetBound2} With the notation of
  Section~\ref{Sec:DetAnal}, we have  
\begin{eqnarray*}
\| \left(I-QQ^{T}\right)A \|_F & \leq & {\displaystyle \| A - Q B_k \|_F
  \leq 
  \sqrt{\left(\sum_{j=k+1}^n \sigma^2_{j}\right) + 
\frac{\alpha^2
  \left\|\widehat{\Omega}_2\right\|_2^2\left\|\widehat{\Omega}_1^{\dagger}\right\|_2^2}
{{1 + \gamma^2
\left\|\widehat{\Omega}_2\right\|_2^2\left\|\widehat{\Omega}_1^{\dagger}\right\|_2^2}}
, }} \\ 
\| \left(I-QQ^{T}\right)A \|_2 & \leq & {\displaystyle \| A
  - Q B_k \|_2 \leq 
  \sqrt{\sigma^2_{k+1} + 
\frac{\alpha^2
  \left\|\widehat{\Omega}_2\right\|_2^2\left\|\widehat{\Omega}_1^{\dagger}\right\|_2^2}
{{1 + \gamma^2
\left\|\widehat{\Omega}_2\right\|_2^2\left\|\widehat{\Omega}_1^{\dagger}\right\|_2^2}}
, }}
\end{eqnarray*}
where 
\[{\displaystyle  \alpha = \sqrt{k}\sigma_{\ell-p+1} 
\left(\frac{\sigma_{\ell-p+1}}{\sigma_{k}}\right)^{2q} \quad
\mbox{and} \quad \gamma = {\displaystyle
  \left(\frac{\sigma_{\ell-p+1}}{\sigma_{1}}\right)\left(\frac{\sigma_{\ell-p+1}}{\sigma_{k}}\right)^{2q} .}}
\]
\end{theorem} 
\begin{REMARK}\label{Rem:SI3} 
  Theorem~\ref{Thm:DetBound2} trivially simplifies to 
\begin{eqnarray*}
\| \left(I-QQ^{T}\right)A \|_F & \leq & {\displaystyle \| A - Q B_k \|_F
  \leq 
  \sqrt{\left(\sum_{j=k+1}^n \sigma^2_{j}\right) + \alpha^2
  \left\|\widehat{\Omega}_2\right\|_2^2\left\|\widehat{\Omega}_1^{\dagger}\right\|_2^2},} \\ 
\| \left(I-QQ^{T}\right)A \|_2 & \leq & {\displaystyle \| A
  - Q B_k \|_2 \leq 
  \sqrt{\sigma^2_{k+1} + {\alpha^2
  \left\|\widehat{\Omega}_2\right\|_2^2\left\|\widehat{\Omega}_1^{\dagger}\right\|_2^2}}.}
\end{eqnarray*}
However, when $\Omega$ is taken to be Gaussian, only
the bounds in Theorem~\ref{Thm:DetBound2} allow average case
analysis for all values of $p$ (See
Section~\ref{Sec:StatExp}.)
\end{REMARK}

\begin{REMARK}\label{Rem:SI2} Not surprisingly, the two key
  factors governing the singular value convergence of
  Algorithm~\ref{Alg:BasicSI} also govern the convergence of
  the low-rank approximation. Hence Remark~\ref{Rem:SI}
  applies equally well to Theorem~\ref{Thm:DetBound2}.
\end{REMARK}

{\bf \noindent Proof of Theorem~\ref{Thm:DetBound2}:} We
first assume that the matrix $H_1$ 
in equation~(\ref{Eqn:OX}) has full column rank. Rewrite
\begin{eqnarray*} {\displaystyle \Sigma_1H_1^T
\left(I + H_1 H_1^T\right)^{-1}
H_1\Sigma_1} & = &
  \left(\left(\left(H_1\Sigma_1\right)^T\left(H_1\Sigma_1\right)\right)^{-1} + \Sigma_1^{-2}\right)^{-1} \\
& = & \left(\left\|H_1\Sigma_1\right\|_2^{-2} I  +
  \Sigma_1^{-2}\right)^{-1} - \\
&& \left(\left(\left\|H_1\Sigma_1\right\|_2^{-2} I  +
  \Sigma_1^{-2}\right)^{-1}-
  \left(\left(\left(H_1\Sigma_1\right)^T\left(H_1\Sigma_1\right)\right)^{-1} + \Sigma_1^{-2}\right)^{-1}\right).
\end{eqnarray*}
The last expression is a symmetric positive semi-definite matrix. This allows us to write
\begin{eqnarray}  
{\displaystyle \left\| \left(I- Q
Q^{T}\right)A_k \right\|_F} & \leq &  {\displaystyle \sqrt{{\bf trace}\left(\Sigma_1H_1^T
\left(I + H_1 H_1^T\right)^{-1}
H_1\Sigma_1\right)}} \nonumber \\
&\leq&{\displaystyle \sqrt{ {\bf trace}\left(\left(
\left\|H_1\Sigma_1\right\|_2^{-2} I + \Sigma_1^{-2} \right)^{-1}
\right)}} =  {\displaystyle \left\|H_1\Sigma_1\right\|_2 \sqrt{
{\bf trace}\left(\Sigma_1\left(
\left\|H_1\Sigma_1\right\|_2^2 I + \Sigma_1^2
\right)^{-1}\Sigma_1\right)}} \nonumber \\
&\leq & {\displaystyle \frac{\sqrt{k}\left\|H_1\Sigma_1\right\|_2
    \sigma_1}{\sqrt{\sigma_1^2 +
      \left\|H_1\Sigma_1\right\|_2^2}}.}\label{Eqn:H1S1bound}
\end{eqnarray}
By a continuity argument, the last relation remains valid
even if $H_1$ does not have full column rank. 

Due to the special form of $H_1$ in equation~(\ref{Eqn:OX}),
we can write $H_1\Sigma_1$ as
\[ {\displaystyle H_1\Sigma_1 = \Sigma_3^{2q+1} \widehat{\Omega}_2
\widehat{\Omega}_1^{\dagger} \left(\begin{array}{c}
\Sigma_1^{-(2q)} 
\cr
0 \end{array}
  \right).} \]
Hence
\[ {\displaystyle \|H_1\Sigma_1\|_2 \leq \sigma_{\ell-p+1}
  \left(\frac{\sigma_{\ell-p+1}}{\sigma_{k}}\right)^{2q} \|
  \widehat{\Omega}_2\|_2 \|\widehat{\Omega}_1^{\dagger}
  \|_2.} \]
Plugging this into equation~(\ref{Eqn:H1S1bound}) and
  dividing both the nomerator and denominator by $\sigma_1$,
\[{\displaystyle \left\| \left(I- Q
Q^{T}\right)A_k \right\|_F \leq \frac{\alpha
  \left\|\widehat{\Omega}_2\right\|_2\left\|\widehat{\Omega}_1^{\dagger}\right\|_2}
  {\sqrt{1 + \gamma^2
  \left\|\widehat{\Omega}_2\right\|_2^2\left\|\widehat{\Omega}_1^{\dagger}\right\|_2^2}}.}
  \]
Comparing this with Theorems~\ref{Thm:Fto2}
and~\ref{Thm:SVD2} proves Theorem~\ref{Thm:DetBound2}. \Q   

\section{Statistical Analysis}\label{Sec:StatAnal} 
This section carries out the needed statistical analysis to
reach our approximation error bounds. In
Section~\ref{Sec:StatTool} we make a list of the statistical
tools used in this analysis; in Section~\ref{Sec:StatExp} we
perform average value analysis on our error bounds; and in
Section~\ref{Sec:StatDev} we provide large deviation bounds.

\subsection{Statistical Tools}\label{Sec:StatTool} 
The simplest of needed statistical results necessary for our
analysis is the following proposition from~\cite{Tropp}.
\begin{prop}\label{Lem:FnormExpV} For fix matrices $S, T$ and
  standard Gaussian matrix $G$, we have 
\[ \mathbb{E} \|S G T\|_2 \leq \|S\|_2 \|T\|_F +
\|S\|_F \|T\|_2. \] 
\end{prop}

The following large deviation bound for the pseudo-inverse
of a Gaussian matrix is also from~\cite{Tropp}.

\begin{lemma}\label{Lem:LargDev} Let $G$ be an $(\ell-p) \times \ell$
Gaussian matrix where $ p \geq 0$ and $ \ell -p \geq 2$. Then ${\bf
  rank}(G) = \ell-p$ with probability $1$. For all $t \geq
1$, 
\[ {\displaystyle {\mathbb{P}} \left\{\left\|G^{\dag}\right\|_2
  \geq \frac{e t\sqrt{\ell}}{p+1}\right\} \leq
  t^{-(p+1)}. } \]
\end{lemma}

The following theorem provides classical tail bounds for
functions of Gaussian matrices. It was taken from~\cite{Bogdanov}[Thm. 4.5.7]. 
\begin{theorem}\label{Thm:Gaussfun}  Suppose that $h$ is a
  real valued Lipschitz function on matrices:  
\[|h(X)-h(Y)| \leq {\cal L} \|X-Y\|_F \quad \mbox{for all
  $X, Y$ and a constant ${\cal L}>0.$} \]
Draw a standard Gaussian matrix $G$. Then
\[{\displaystyle \mathbb{P} \left\{h(G) \geq \mathbb{E} h(G)
  + {\cal L} u\right\} \leq e^{-u^2/2}.} \] 
\end{theorem}

The two propositions below will be used in our average case error bounds analysis, both for
singular values and rank-k approximations. Their proofs are
lengthy and can be found in the {\em Supplemental Material}.

\begin{prop}\label{Lem:Omega2Bound} 
Let $\alpha > 0$, $\beta > 0$, $\gamma > 0$ and $\delta >
0$, and let $G$ be an $m \times n $ Gaussian matrix. Then
\begin{eqnarray}
{\displaystyle \mathbb{E} \left( \frac{1}{\sqrt{1+\alpha^2
\|G\|_2^2}}\right) } &\geq & {\displaystyle
\frac{1}{\sqrt{1+\alpha^2 {\cal C} ^2}}} \label{Eqn:SVAB1} \\
{\displaystyle \mathbb{E}\left( {\sqrt{\delta^2+
\frac{\alpha^2 \|G\|_2^2}{\beta^2 + \gamma^2
\|G\|_2^2}}}\right)} &\leq & {\displaystyle
{\sqrt{\delta^2+\frac{\alpha^2 {\cal C} ^2}{\beta^2 + \gamma^2
{\cal C} ^2}}}} , \label{Eqn:SVAB2}
\end{eqnarray}
where ${\cal C}  = \sqrt{m} + \sqrt{n} + 7$.
\end{prop}

There are lower and upper bounds similar to
Proposition~\ref{Lem:Omega2Bound} for the pseudo-inverse of
a Gaussian, with a significant complication. When $G$ is a
square Gaussian matrix, it is non-singular with probability
$1$. However, the probability density function for its
pseudo-inverse could have a very long tail according to
Lemma~\ref{Lem:LargDev}. A similar argument could also be made
when $G$ is almost a square matrix. This complication will
have important implications for parameter choices in
Algorithm~\ref{Alg:RandSI} (see Sections~\ref{Sec:StatExp}
and~\ref{Sec:StatDev}.) Function $\log(\cdot)$ below is base-$e$.

\begin{prop}\label{Lem:Omega1Bound} 
Let $\alpha > 0$, $\beta > 0$, $\gamma > 0$ and $\delta >
0$, and let $G$ be an $(\ell - p) \times \ell $ Gaussian
matrix. Then ${\bf rank}(G) = \ell - p$ with probability
$1$, and 
\begin{eqnarray}
{\displaystyle \mathbb{E} \left( \frac{1}{\sqrt{1+\alpha^2
      \|G^{\dagger}\|_2^2}}\right) } &\geq & 
\left\{\begin{array}{ll} {\displaystyle 
\frac{1}{\sqrt{1+\alpha^2 {\cal C} ^2}}}& \mbox{for $p \geq
  2$,} \cr
& \cr
 {\displaystyle \frac{1}{ 1+ \alpha^2 {\cal C} ^2 \log 
      \frac{2\sqrt{1 + \alpha^2{\cal C}^2 }}{ \alpha{\cal C}}}} &  \mbox{for $p = 1$,} \cr
& \cr
{\displaystyle \frac{1}{1 + \alpha {\cal C}}} & \mbox{for $p = 0.$} 
\end{array} \right. \label{Eqn:SVAB3} \\
{\displaystyle \mathbb{E}\left( {\sqrt{\delta^2+
      \frac{\alpha^2 \|G^{\dagger}\|_2^2}{\beta^2 + \gamma^2
      \|G^{\dagger}\|_2^2}}}\right)} &\leq &
\left\{\begin{array}{ll}
 {\displaystyle  
{\sqrt{\delta^2+\frac{\alpha^2 {\cal C} ^2}{\beta^2 + \gamma^2
      {\cal C} ^2}}}} & \mbox{for $p \geq
  2$,} \cr
& \cr
{\displaystyle  \delta+
  \frac{\alpha^2\left(\ell-1\right)}{\delta \beta^2} \left(2 +
  \frac{1}{2} \log \left(1 +
  \frac{\delta^2\beta^2}{\alpha^2}\right)\right)} 
& \mbox{for $p = 1$,} \cr
& \cr
{\displaystyle  
\delta + \frac{4\sqrt{\ell} \sqrt{\delta^2\gamma^2 +
      \alpha^2}}{\beta} \log  \left(1 +
      \left(\frac{\alpha}{\delta\gamma}\right)^2\right)} & \mbox{for $p = 0$,} 
\end{array} \right.  \label{Eqn:SVAB4} 
\end{eqnarray}
where ${\displaystyle
  {\cal C}  = \frac{ 4e \sqrt{\ell}}{p+1} }$.
\end{prop}

\subsection{Average Case Error Bounds}\label{Sec:StatExp}
This section is devoted to the average case analysis of
Algorithm~\ref{Alg:RandSI}. This
work requires us to study the average case behavior on the
upper and lower bounds in Theorems~\ref{Thm:SVdetBound} and~\ref{Thm:DetBound2}. As observed in
Section~\ref{Sec:Random}, the distribution of a standard
Gaussian matrix is rotationally invariant, and hence 
the matrices $\widehat{\Omega}_1$ and
$\widehat{\Omega}_{2}$ are themselves independent standard
Gaussian matrices. With the tools established in
Section~\ref{Sec:StatTool}, our analysis here consists
mostly of stitching together the right pieces from there.

We first analyze the singular value lower bounds in
Theorems~\ref{Thm:SVAvBounds}. This will require separate analysis for
$ p \geq 2$, $p =1$, and $p = 0$, as suggested in
Section~\ref{Sec:StatTool}. We then analyze the low-rank approximation
bounds in Theorem~\ref{Thm:DetBound2}, which also requires separate
analysis for the same three cases of $ p$. Throughout
Section~\ref{Sec:StatExp}, we will need the following definition for
any $0 \leq p \leq \ell$: 
\[{\displaystyle {\cal C}_1 =  \sqrt{n-\ell+p}+\sqrt{\ell} + 7, \quad 
{\cal C}_2 =  \frac{4 e
  \sqrt{\ell}}{p+1}, \quad {\cal C} = {\cal C}_1 {\cal C}_2, 
\quad \mbox{and} \quad \tau_j = \frac{\sigma_{\ell-p+1}}
{\sigma_j}.} \]
\begin{theorem} \label{Thm:SVAvBounds} Let $A = U \Sigma
  V^{T}$ be the SVD of $A$, and let $QB_k$ be a rank-$k$
  approximation computed by Algorithm~\ref{Alg:RandSI}.
  Then for $j =1, \cdots, k,$
\begin{equation}\label{Eqn:proof3}
{\displaystyle \mathbb{E} \left( \sigma_j\left(QB_k\right) \right) \geq 
\left\{\begin{array}{ll} {\displaystyle \frac{\sigma_j}{\sqrt{1+{\cal
      C}^2\tau_j^{4q+2}}}} &  \mbox{for $p
  \geq 2$,} \cr 
& \cr
 {\displaystyle \frac{\sigma_j}{1+{\cal C}^2\tau_j^{4q+2}
     \log \sqrt{{\cal C}^2 + \tau_j^{-(4q+2)}}}} &  \mbox{for $p
  = 1$,} \cr 
& \cr
{\displaystyle \frac{\sigma_j}{1+{\cal C} \tau_j^{2q+1}}}  &  \mbox{for $p
  = 0$.} \end{array} \right. }
\end{equation}
\end{theorem}

\begin{REMARK}\label{Rem:LargeDev1} The value of $p$ is not
part of Algorithm~\ref{Alg:RandSI} and can thus be
arbitrarily chosen within $[0, \ell -k].$ Since our bounds
for $p\leq 1$ are worse than that for $p \geq 2$, they
should probably not be used unless $\ell-k \leq 1$ or unless
there is a large singular value gap at $\sigma_{\ell}$ or
$\sigma_{\ell+1}$.
\end{REMARK}

\begin{REMARK}\label{Rem:LargeDev2} Theorem~\ref{Thm:SVAvBounds} strongly suggests that 
in general some over-sampling in the number of columns can
significantly improve convergence in the singular value
approximation. This is consistent with the
literature~\cite{ChengGimbutasMartinssonRokhlin,DrineasKannanMahoneyII,DrineasMahoneyMuthukrishnan,FriezeKannanVempala,FriezeKannanVempalaII,LibertyAilonSinger,Liberty,NguyenDoTran,WoolfeLibertyRokhlinTygert,sarlos,Tropp} 
and is very significant for practical implementations. 
\end{REMARK}

\begin{REMARK}\label{Rem:BadPower} A typical implementation of the
classical subspace iteration method in general and the classical power method
in particular chooses $\ell = k$, which leads to $p =
0$. Theorem~\ref{Thm:SVAvBounds} implies that this choice in general
leads to slower convergence than $p > 0$ and thus should be avoided. We will
elaborate this point in more detail in Section~\ref{Sec:StatDev} and
provide numerical evidence to support this conclusion in
Section~\ref{Sec:Num}.  
\end{REMARK}

\begin{REMARK}\label{Rem:lowerbound}
Since $\tau_j \leq 1$ for all $j$, 
Theorem~\ref{Thm:SVAvBounds} implies that for $p \geq 2$
and for all $j \leq k$, 
\[{\displaystyle \mathbb{E} \left( \sigma_j\left(QB_k\right)
  \right) \geq \frac{\sigma_j}{\sqrt{1+{\cal C}^2}}. } \] In
other words, Algorithm~\ref{Alg:RandSI} approximates the
leading $k$ singular values by a good fraction on average,
regardless of how the singular values are distributed, even
for $q = 0$. This result is surprising and yet valuable. It
will have applications in condition number estimation (see
Sections~\ref{Sec:StatDev} and~\ref{Sec:CondEst} for more
discussion.)
\end{REMARK}

\begin{REMARK}\label{Rem:smallp}
For matrices with rapidly decaying singular values,
convergence could be so rapid that one could even set $q =
0$ in some cases (Section~\ref{Sec:StatDev}.) This is the basis of the excitement
about Algorithm~\ref{Alg:RandSI} in that very little work is
typically sufficient to realize an excellent low-rank
approximation. The faster the singular values decay, the
faster Algorithm~\ref{Alg:RandSI} converges.
\end{REMARK}

\begin{REMARK}
Kuczy\'{n}ski and Wo\'{z}niakowski~\cite{Kuczynski2}
developed probabilistic error bounds for computing the
largest eigenvalue of an SPD matrix by the power method
for a unit start vector under the uniform distribution. 
Their results correspond to the case of $\ell = k = 1$ and
$p = 0$ in Theorem~\ref{Thm:SVAvBounds}. However, our
results appear to be much stronger.
\end{REMARK}

{\noindent \bf Proof of Theorem~\ref{Thm:SVAvBounds}:} As in Theorem~\ref{Thm:SVdetBound},
we will only prove Theorem~\ref{Thm:SVAvBounds} for $j =
k$. All other values of $j$ can be proved by simply citing
Theorem~\ref{Thm:SVAvBounds} for a rank-$j$ SVD
truncation. Since $\widehat{\Omega}_2$ and
$\widehat{\Omega}_1$ are independent of each other, we will
take expectations over $\widehat{\Omega}_2$ and
$\widehat{\Omega}_1$ in turn, based on Propositions~\ref{Lem:Omega2Bound}
and~\ref{Lem:Omega1Bound}.

Let ${\displaystyle \alpha  
= \left\|\widehat{\Omega}_1^{\dagger}\right\|_2 \tau_k^{2q+1}}$. By
Theorem~\ref{Thm:SVdetBound} and Proposition~\ref{Lem:Omega2Bound}, 
\begin{equation}\label{Eqn:proof1}
{\displaystyle \mathbb{E} \left(\sigma_k\left(QB_k\right)     \; \middle\vert\; 
\widehat{\Omega}_1 \right) \geq \frac{\sigma_k}{\sqrt{
1 + \alpha^2{\cal C}_1^2}}.}
\end{equation}

For $p \geq 2$, we further take expectation over $\widehat{\Omega}_1$
according to Proposition~\ref{Lem:Omega1Bound}. By equation~(\ref{Eqn:proof1}), 
\[\mathbb{E}\left(\sigma_k\left(QB_k\right) \right)  = 
\mathbb{E}\left( \mathbb{E}\left(\sigma_k \left(QB_k\right)
    \; \middle\vert\; \widehat{\Omega}_1 \right)\right) \geq  \mathbb{E}\left( \frac{\sigma_k}{\sqrt{1 +
    \left(\left\|\widehat{\Omega}_1^{\dagger}\right\|_2
    \tau_k^{2q+1}\right)^2{\cal C}_1^2}}\right)  \geq
  \frac{\sigma_k}{\sqrt{1 +{\cal C}^2\tau_k^{4q+2}}}. \]
To complete the proof, we note that the results for $p = 1$
  and $p = 0$ can be obtained similarly by taking expectation of
$\widehat{\Omega}_1$ over equation~(\ref{Eqn:proof1}) and simplifying.  \Q

It is now time for average case analysis of low-rank matrix
approximations. Again, we base our arguments on
Propositions~\ref{Lem:Omega2Bound} and~\ref{Lem:Omega1Bound}. For ease of notation, let 
\[ {\displaystyle \widehat{\delta}_{k+1} =
  \sqrt{\sum_{j=k+1}^n\sigma_j^2}.} \]
For the sake of simplicity, in Theorem~\ref{Thm:SVAvBounds}
below we have omitted 
\[\mathbb{E} \| \left(I-QQ^{T}\right)A \|_F  \leq
	 {\displaystyle \mathbb{E} \| A - Q B_k \|_F .} \]
\begin{theorem} \label{Thm:LRKAvBounds} Let $QB_k$ be a rank-$k$
  approximation computed by Algorithm~\ref{Alg:RandSI}. Then
\begin{eqnarray*}
 {\displaystyle \mathbb{E} \| A - Q B_k \|_F}
&  \leq &  \left\{\begin{array}{ll} 
{\displaystyle 
  \sqrt{\widehat{\delta}_{k+1}^2 + k {\cal C}^2
    \sigma_{\ell-p+1}^2\tau_k^{4q}}} & \mbox{for $p \geq
  2$,} \cr
& \cr
 {\displaystyle \widehat{\delta}_{k+1} +
 \frac{k{\cal
 C}^2\sigma_{\ell-p+1}^2\tau_k^{4q}}{\widehat{\delta}_{k+1}}
 \log \sqrt{{\cal C}^2 + \frac{1}{k}
 \left(\frac{\widehat{\delta}_{k+1}}{\sigma_{\ell-p+1}}\right)^2\tau_k^{-4q}}}&
 \mbox{for $p = 1$,} \cr
& \cr
{\displaystyle \widehat{\delta}_{k+1} + \sqrt{n}{\cal C}\sigma_{\ell-p+1}\tau_k^{2q}
  \log \left(1 +  k
  \left(\frac{\sigma_1}{\widehat{\delta}_{k+1}}\right)^2\right)}& \mbox{for $p = 0.$} 
\end{array} \right. \\
& \\
 {\displaystyle \mathbb{E} \| A - Q B_k \|_2}
&  \leq &  \left\{\begin{array}{ll} 
{\displaystyle 
  \sqrt{{\sigma}^2_{k+1} + k {\cal C}^2
    \sigma_{\ell-p+1}^2\tau_k^{4q}}} & \mbox{for $p \geq
  2$,} \cr
& \cr
 {\displaystyle {\sigma}_{k+1} +
 \frac{k{\cal
 C}^2\sigma_{\ell-p+1}^2\tau_k^{4q}}{{\sigma}_{k+1}}
 \log \sqrt{{\cal C}^2 + \frac{1}{k}
 \left(\frac{{\sigma}_{k+1}}{\sigma_{\ell-p+1}}\right)^2\tau_k^{-4q}}}&
 \mbox{for $p = 1$,} \cr
& \cr
{\displaystyle {\sigma}_{k+1} + \sqrt{(k+1)}{\cal C}\sigma_{\ell-p+1}\tau_k^{2q}
  \log \left(1 +  k
  \left(\frac{\sigma_1}{{\sigma}_{k+1}}\right)^2\right)}& \mbox{for $p = 0.$} 
\end{array} \right. 
\end{eqnarray*}
\end{theorem}

\begin{REMARK} Remark~\ref{Rem:quad} applies to
 Theorem~\ref{Thm:LRKAvBounds} as well.
\end{REMARK} 

{\noindent \bf Proof of Theorem~\ref{Thm:LRKAvBounds}:} We only prove
Theorem~\ref{Thm:LRKAvBounds} for the Frobenius norm. The
case for the 2-norm is completely analogous. As in the proof for
Theorem~\ref{Thm:SVAvBounds}, this one involves taking
expectations over $\widehat{\Omega}_2$ first and
$\widehat{\Omega}_1$ next. Let $\displaystyle \widehat{\delta}_{k+1} =
\sqrt{\sum_{j=k+1}^n \sigma^2_{j}}$. Fixing $\widehat{\Omega}_1$ 
in Theorem~\ref{Thm:DetBound2} and taking expectation on
$\widehat{\Omega}_2$ according to Proposition~\ref{Lem:Omega2Bound},
we obtain immediately
\begin{eqnarray}
\mathbb{E} {\displaystyle \| A - Q B_k \|_F}
&  \leq  & {\displaystyle 
  \sqrt{\widehat{\delta}_{k+1}^2 + 
\frac{\alpha^2 {\cal C}_1^2 \left\|\widehat{\Omega}_1^{\dagger}\right\|_2^2} 
{{1 + \gamma^2 {\cal C}_1^2 \left\|\widehat{\Omega}_1^{\dagger}\right\|_2^2}}
,}} \label{proof2a}
\end{eqnarray}
with ${\displaystyle  \alpha = \sqrt{k}\sigma_{\ell-p+1} \tau_k^{2q}}$ 
and ${\displaystyle  \gamma = {\displaystyle 
  \left(\frac{\sigma_{\ell-p+1}}{\sigma_{1}}\right)\tau_k^{2q} .}}$ 

For $p \geq 2$, we further take expectation over $\widehat{\Omega}_1$ according
to Proposition~\ref{Lem:Omega1Bound}. By
equation~(\ref{Eqn:SVAB4}),
\begin{eqnarray*}
{\displaystyle \mathbb{E}\left( \mathbb{E}\| A - Q B_k \|_F  \; \middle\vert\; \widehat{\Omega}_1 \right)} 
& \leq & {\displaystyle \mathbb{E} 
 \left(
  \sqrt{\widehat{\delta}^2_{k+1} + 
\frac{\alpha^2 {\cal C}_1^2 \left\|\widehat{\Omega}_1^{\dagger}\right\|_2^2} 
{{1 + \gamma^2 {\cal C}_1^2 \left\|\widehat{\Omega}_1^{\dagger}\right\|_2^2}}}
\right) } \\
& \leq & {\displaystyle  \sqrt{\widehat{\delta}^2_{k+1} + \frac{\alpha^2 {\cal C}^2} 
{{1 + \gamma^2 {\cal C}^2}}} \leq \sqrt{\widehat{\delta}^2_{k+1} + \alpha^2 {\cal C}^2},} 
\end{eqnarray*}
which is the Frobenius norm upper bound in
Theorem~\ref{Thm:LRKAvBounds}.

For $p = 1$, we again take expectation over
$\widehat{\Omega}_1$ in equation~(\ref{proof2a}) according to Proposition~\ref{Lem:Omega1Bound}:
\begin{eqnarray*}
{\displaystyle \mathbb{E}\| A - Q B_k \|_F} & \leq & {\displaystyle \mathbb{E} 
 \left(\sqrt{\widehat{\delta}^2_{k+1} + 
\frac{\alpha^2 {\cal C}_1^2 \left\|\widehat{\Omega}_1^{\dagger}\right\|_2^2} 
{{1 + \gamma^2 {\cal C}_1^2 \left\|\widehat{\Omega}_1^{\dagger}\right\|_2^2}}}
\right) } \\
& \leq & {\displaystyle \widehat{\delta}_{k+1} + \frac{\alpha^2
 {\cal C}_1^2\left(\ell-1\right)}{\widehat{\delta}_{k+1}} \left( 2
 + \frac{1}{2} \log \left(1 +  \frac{ \widehat{\delta}_{k+1}^2}{\alpha^2 {\cal C}_1^2}\right)\right),}
\end{eqnarray*}
which is bounded above by the corresponding expression in
Theorem~\ref{Thm:LRKAvBounds}.

Now, we turn our attention to the case $p = 0$. Taking expectations as before, 
\begin{eqnarray*}
{\displaystyle \mathbb{E}\| A - Q B_k \|_F} & \leq & {\displaystyle \mathbb{E} 
 \left(\sqrt{\widehat{\delta}^2_{k+1} + 
\frac{\alpha^2 {\cal C}_1^2 \left\|\widehat{\Omega}_1^{\dagger}\right\|_2^2} 
{{1 + \gamma^2 {\cal C}_1^2 \left\|\widehat{\Omega}_1^{\dagger}\right\|_2^2}}}
\right) } \\
& \leq & {\displaystyle \widehat{\delta}_{k+1} + 4\sqrt{\ell}
 \sqrt{\widehat{\delta}^2_{k+1} \gamma^2{\cal C}_1^2 + \alpha^2
 {\cal C}_1^2 } \log \left(1 +  \left(\frac{\alpha{\cal
 C}_1}{\widehat{\delta}^2_{k+1}\gamma{\cal
 C}_1}\right)^2\right)} \\
& = & {\displaystyle \widehat{\delta}_{k+1} + 4\sqrt{\ell}{\cal C}_1
 \sqrt{\widehat{\delta}^2_{k+1} \gamma^2 + \alpha^2
  } \log \left(1 +  \left(\frac{\alpha}{\widehat{\delta}_{k+1}\gamma}\right)^2\right).}
\end{eqnarray*}
Plugging in the expressions for $\alpha$ and $\gamma$ in
equation~(\ref{proof2a}), 
\[{\displaystyle \mathbb{E}\| A - Q B_k \|_F}  \leq 
{\displaystyle \widehat{\delta}_{k+1} + 4\sqrt{\ell}{\cal C}_1
 \sqrt{\left(\frac{\widehat{\delta}_{k+1}}{\sigma_1}\right)^2
 + k}\,\, \sigma_{\ell-p+1}\tau_k^{2q}
  \log \left(1 +  k \left(\frac{\sigma_1}{\widehat{\delta}_{k+1}}\right)^2\right),}
\]
which is bounded above by the corresponding expression in Theorem~\ref{Thm:LRKAvBounds} since
$\widehat{\delta}_{k+1} \leq \sqrt{n-k} \,\,{\sigma_1}$. \Q
 
\subsection{Large Deviation Bounds}\label{Sec:StatDev} In
this section we develop approximation error tail bounds. Theorems~\ref{Thm:SVdetBound}
and~\ref{Thm:DetBound2} dictate that our main focus will be in
developing probabilistic upper bounds on $
\left\|\widehat{\Omega}_{2} \right\|_2
\left\|\widehat{\Omega}_1^{\dag}\right \|_2$. 

\begin{theorem} \label{Thm:HDBoundsIII} 
Let $A = U \Sigma V^{T}$ be the SVD of $A$, and $0 \leq p
\leq \ell - k$.  Further let $QB_k$ be a rank-$k$ 
approximation computed by
Algorithm~\ref{Alg:RandSI}. Given any $0 < \Delta \ll 1$,
define 
\[ {\displaystyle {\cal C}_{\Delta} = \frac{e \sqrt{\ell}}{p+1}
  \left(\frac{2}{\Delta}\right)^{\frac{1}{p+1}} 
\left(\sqrt{n-\ell+p}+\sqrt{\ell} + \sqrt{2 \log
  \frac{2}{\Delta}}\right). } \]
We must have for $j =1, \cdots, k$,
\begin{eqnarray*}
{\displaystyle \sigma_{j}\left(QB_k)\right)} &\geq
&{\displaystyle \frac{ \sigma_{j}}{\sqrt{1 + {\displaystyle {\cal C}_{\Delta}^2 
\left(\frac{\sigma_{\ell - p +
    1}}{\sigma_j}\right)^{4q+2}}}}}, 
\end{eqnarray*}
and
\begin{eqnarray*}
\| \left(I-QQ^{T}\right)A \|_F & \leq & {\displaystyle \| A - QB_k \|_F
\leq  \sqrt{\left(\sum_{j=k+1}^n \sigma^2_{j}\right) + k {\cal
    C}_{\Delta}^2 \sigma_{\ell-p+1}^2 \left(\frac{\sigma_{\ell-p+1}}
{\sigma_{k}}\right)^{4q}} }, \\
\| \left(I-QQ^{T}\right)A \|_2 & \leq & {\displaystyle \| A - QB_k \|_2
\leq  \sqrt{ \sigma^2_{k+1} + k {\cal
    C}_{\Delta}^2 \sigma_{\ell-p+1}^2 \left(\frac{\sigma_{\ell-p+1}}
{\sigma_{k}}\right)^{4q}} }.
\end{eqnarray*}
with exception probability at most ${\Delta} $.
\end{theorem}

\begin{REMARK}\label{Rem:Conv}
Like the average case, the factor $\sigma_{\ell-p+1}^2$
shows up in all three bounds, for all $q \geq 0$. Hence Algorithm~\ref{Alg:RandSamI} can make
significant progress toward convergence in case of rapidly
decaying singular values in $A$, with probability
$1-\Delta$. This is clearly a much stronger result than
Theorem~\ref{Thm:Troppbound2}. 
\end{REMARK}

\begin{REMARK}\label{Rem:Fail}
While the value of $\Delta$ could be set arbitrarily tiny,
it can never be set to $0$. This implies that there is a
chance, however arbitrarily small, that
Algorithm~\ref{Alg:RandSI} might not converge according to
the bounds in Theorem~\ref{Thm:HDBoundsIII}. This small
exception chance probably has less to do with
Algorithm~\ref{Alg:RandSI} and more to do with the inherent
complexity of efficiently computing accurate matrix
norms. Since Algorithm~\ref{Alg:RandSI} accesses $A$ only
through the $2q+2$ matrix-matrix products of the form $A X$
or $A^TY$ for different $X$ and $Y$ matrices, it can be used
to efficiently compute $\|A^{-1}\|_2$ (setting $k = 1$)
provided that a factorization of $A$ is available or if $A$
is itself a non-singular triangular matrix. On the other
hand, it is generally expected that even {\rm estimating}
$\|A^{-1}\|_2$ to within a constant factor independent of
the matrix $A$ must cost as much, asymptotically, as
computing $A^{-1}$. Demmel, Diament, and
Malajovich~\cite{ddm2001} show that the cost of computing an
estimate of $\|A^{-1}\|$ of guaranteed quality is at least
the cost of testing whether the product of two $n \times n$
matrices is zero, and performing this test is conjectured to
cost as much as actually computing the
product~\cite[p. 288]{Highambk}. Since
Algorithm~\ref{Alg:RandSI} costs only $O(n^2 q \ell)$
operations to provide a good estimate for $\|A^{-1}\|_2$, it
probably can not be expected to work without {\em any}
failure. See Section~\ref{Sec:CondEst} for more comments.
\end{REMARK}

{\bf \noindent Proof of Theorem~\ref{Thm:HDBoundsIII}:} Since
$\widehat{\Omega}_2$ and $\widehat{\Omega}_1$ are independent from
each other, we can study how the error depends on the matrix
$\widehat{\Omega}_2$ when $\widehat{\Omega}_1$ is reasonably bounded. To this end, we define an
event as follows:
\[ {\displaystyle  {\bf E}_t = \left\{ \widehat{\Omega}_1:
\left\|\widehat{\Omega}_1^{\dag} \right\|_2 \leq t {\cal L}\right\},
  \quad \mbox{where} \quad {\cal L} = \frac{e\sqrt{\ell}}{p+1}. }\]
Invoking the conclusion of Lemma~\ref{Lem:LargDev}, we find
  that 
\begin{equation}\label{Eqn:rare}
\mathbb{P}\left({\bf E}_t^c\right) \leq t^{-(p+1)}. 
\end{equation}
In other words, we have just shown
that ${\displaystyle \left\|\widehat{\Omega}_1^{\dag}
  \right\|_2 \leq t{\cal L}}$ with probability at least $1 -
t^{-(p+1)}$. 

Below we consider the function  
\[ {\displaystyle h(X)  = \left\| X \right\|_2
\left\|\widehat{\Omega}_1^{\dag}\right \|_2 , } \]
where $X$ has the same dimensions as
$\widehat{\Omega}_2$. It is straightforward to show that  
\[{\displaystyle  |h(X) - h(Y)| \leq \left\|\widehat{\Omega}_1^{\dag} \right\|_2 \left\|X-Y\right\|_F \leq t
{\cal L} \left\|X-Y\right\|_F ,} \]
under event ${\bf E}_t$. Also under event ${\bf E}_t$ and by
Proposition~\ref{Lem:FnormExpV}, we have   
\[{\displaystyle  \mathbb{E} h(X)} \leq {\displaystyle 
  \left(\sqrt{n-\ell+p}+\sqrt{\ell}\right)
  \left\|\widehat{\Omega}_1^{\dag}\right \|_2}  
\leq  {\displaystyle \frac{e
  t\sqrt{\ell}}{p+1}
  \left(\sqrt{n-\ell+p}+\sqrt{\ell}\right) \stackrel{def}= t
  {\cal E} . } \]
Applying the concentration of measure equation,
Theorem~\ref{Thm:Gaussfun}, conditionally to
$\widehat{\Omega}_{2}$ under event ${\bf E}_t$,  
\[{\displaystyle \mathbb{P} \left\{
\left\|\widehat{\Omega}_{2} \right\|_2
\left\|\widehat{\Omega}_1^{\dag}\right \|_2
 \geq t {\cal E} + t {\cal L} u  \; \middle\vert\;  {\bf E}_t  \right\} \leq  e^{-u^2/2}  . } \] 
Use the equation~(\ref{Eqn:rare}) to remove the
restriction on $\widehat{\Omega}_1$, therefore,  
\[{\displaystyle \mathbb{P} 
\left\{
\left\| \widehat{\Omega}_{2} \right\|_2
\left\|\widehat{\Omega}_1^{\dag}\right \|_2
\geq t {\cal E} + t {\cal L} u  \right\} \leq
 t^{-(p+1)} + e^{-u^2/2}  ,} \]
Now we choose 
\[ {\displaystyle t = \left(\frac{2}{\Delta}\right)^{1/(p+1)}
  \quad \mbox{and} \quad u = \sqrt{2 \log
  \frac{2}{\Delta }}} \]
so that $t^{-(p+1)} + e^{-u^2/2} = \Delta$. With this choice
  of $t$ and $u$, 
\[{\displaystyle t {\cal E} + t {\cal L} u   =  {\cal
  C}_{\Delta } \quad \mbox{or} \quad 
\mathbb{P} 
\left\{
\left\| \widehat{\Omega}_{2} \right\|_2
\left\|\widehat{\Omega}_1^{\dag}\right \|_2
\geq {\cal C}_{\Delta} \right\} \leq \Delta .} \]
Plugging this bound into the formulas in
Theorem~\ref{Thm:SVdetBound} and Remark~\ref{Rem:SI3} proves
Theorem~\ref{Thm:HDBoundsIII}. \Q

While the value of oversampling size $p$ does not look so important in
the average case error bounds as long as $p \geq 2$, it makes an
oversized difference in large deviation bounds. Consider the case $p =
2$ with a tiny $\Delta >0$. In this case, $ {\displaystyle
{\cal C}_{\Delta} }$ may still be quite large, and quite a few extra
number of iterations might be necessary to ensure satisfactory
convergence with small exception probability. 

For $p \leq 1$, the large deviation bound is brutal. For very small
values of $k$, such as $1$ in the case of the randomized
power method (see Algorithm~\ref{Alg:RandPM}), it seems
unreasonable to require a relatively large value of $p$. On the other
hand, a small $p$ value would significantly impact convergence. We
will address this conflicting issue of choosing $p$ further in
Section~\ref{Sec:Num}.

But for any large enough values of $k$ (such as $k = 20$ or more, for
example,) a reasonable choice would be to choose $p$ so
${\displaystyle \left(\frac{2}{{\Delta} }\right)^{1/(p+1)}}$ is a
modest number. We will now choose 
\begin{equation}\label{Eqn:p}
{\displaystyle p = \lceil \log_{10}\left( \frac{2}{{\Delta} } \right)
\rceil - 1}. 
\end{equation} 
This choice gives ${\displaystyle
\left(\frac{2}{{\Delta} }\right)^{1/(p+1)} \leq 10}$.  For a typical
choice of ${\Delta} = 10^{-16}$, equation~(\ref{Eqn:p}) gives $ p =
16$. For this value of ${\Delta} $, the exception probability is
smaller than that of matching DNA
fingerprints~\cite{RischDevlin}. Given that the "random numbers"
generated on modern computers are really only {\em pseudo random
numbers} that may have quite different upper tail
distributions than the true Gaussian (see, for
example~\cite{ThomasLuk,WichmannHill}), and given that only finite precision computations are
typically done in practice, it is probably meaningless to require
${\Delta} $ to be much less than $10^{-16}$, the double precision.
Additionally, with this choice of $p$, the large deviation bounds are
very similar to the average case error bounds, suggesting that the
typical behavior is also the worst case behavior, with
probability $1- \Delta$.

Our final observation on Theorem~\ref{Thm:HDBoundsIII} is so
important that we present it in the form of a Corollary. We
will not prove it 
because it is a direct consequence.
\begin{corollary} \label{Cor:HDBoundsIIIv2} In the notation
  of Theorem~\ref{Thm:HDBoundsIII}, we must have for $j =1,
\cdots, k$, 
\begin{equation}\label{Eqn:CondCor}
{\displaystyle \sigma_{j}\left(QB_k)\right)} \geq
{\displaystyle \frac{ \sigma_{j}}{\sqrt{1 + {\displaystyle
	{\cal C}_{\Delta}^2 }}}} \quad \mbox{and} \quad 
{\displaystyle \| A - QB_k \|_2 }  \leq   {\displaystyle \sigma_{k+1} \sqrt{1 + k {\cal
    C}_{\Delta}^2 }}
\end{equation}
with exception probability at most ${\Delta} $. 
\end{corollary} 

This is a surprisingly strong result. We will discuss its
implications in terms of rank-revealing factorizations in
Section~\ref{Sec:RankRev} and condition number estimation 
in Section~\ref{Sec:CondEst}. 

\section{Rank-revealing Factorizations}\label{Sec:RankRev}
Rank-revealing factorizations were first discussed in
Chan~\cite{chan}. Generally speaking, there are
rank-revealing UTV
factorizations~\cite{FierroHansen,stewart93}, QR
factorizations~\cite{chanhan2,ci,ge7}, and LU
factorizations~\cite{MG,pan96}. While there is no uniform
definition of {\em the} rank-revealing factorization, a
comparison of different forms of rank-revealing
factorizations has appeared in Foster and
Liu~\cite{fosterliu}. For the discussions in this section,
we make the following definition, which is loosely
consistent with those in~\cite{fosterliu}. 
\begin{DEFINITION}
Given $m \times n$ matrices $A$ and $B$ and integer $k < 
\min(n,m)$, we call $B$ a {\em rank-revealing rank-$k$
  approximation} to $A$ if ${\bf rank}(B) \leq k$ and if 
there exist polynomials $c_1(m,n)$, and $c_2(m,n)$ such that 
\begin{eqnarray}
\sigma_j(B) & \geq & \frac{\sigma_j(A)}{c_2(m,n)}, \quad j =
1, \cdots, k, \label{Eqn:c2} \\
\|A-B\|_2 &\leq & c_1(m,n) \sigma_{k+1}(A). \label{Eqn:c1}
\end{eqnarray}
\end{DEFINITION}

A {\em rank-revealing rank-$k$ approximation} differs from
an ordinary rank-$k$ approximation in the extra
condition~(\ref{Eqn:c2}), which requires some accuracy in
{\em all} $k$ leading singular values. Therefore a
rank-revealing rank-$k$ approximation is likely a stronger
approximation than a simple low rank approximation. To see
why~(\ref{Eqn:c2}) is so important, we consider for an
example the case where the leading $k+1$ singular values of
$A$ are identical: $\sigma_1(A) = \cdots =
\sigma_{k+1}(A)$. This includes the $n \times n$ identity
matrix as a special case. Now choose $\theta = 1$ in
equation~(\ref{Eqn:B}). It follows that $B = 0$ is an {\em
optimal} rank-$k$ approximation to $A$, which is likely
unacceptable to most users. On the
other hand, $B = 0$ obviously does not satisfy
condition~(\ref{Eqn:c2}) for any polynomial $c_2(m,n)$, and
therefore is not a rank-revealing rank-$k$ approximation to
$A$. Similarly, {\em any} orthogonal matrix $Q$ would
satisfy the bound in Theorem~\ref{Thm:Troppbound2} for such
an $A$ matrix, and only the matrix $Q$ from
Algorithm~\ref{Alg:RandSI} would satisfy
Theorem~\ref{Thm:HDBoundsIII}. 

By definition, Algorithm~\ref{Alg:RandSI} produces 
a rank-revealing rank-$k$ approximation with probability at
least $1-\Delta$. In this section, we compare this
approximation with the strong
RRQR factorization developed in Gu and Eisenstat~\cite{ge7}.
\begin{theorem} \label{Thm:GuEisenstat_SISC96} (Gu and Eisenstat~\cite{ge7}) Let
$A$ be an $m \times n$ matrix and let $1 \leq k \leq
  \min(m,n)$. For any given parameter $f>1$, there exists a
  permutation $\Pi$ such that 
\[ A \Pi = Q  \left(
\begin{array}{cc} 	
R_{11} & R_{12} \\
& R_{22} \\
\end{array}\right) , \]
where for any $1 \leq i \leq k$ and $1 \leq j \leq n-k$,
\begin{equation}\label{eqn:RRQR}
1 \leq \frac{\sigma_i(A)}{\sigma_i(R_{11})},
  \frac{\sigma_j(R_{22})}{\sigma_{k+j}} \leq
  \sqrt{1+f^2k(n-k)}. 
\end{equation}
\end{theorem}

Let ${\displaystyle \widehat{A}_k = Q  \left(
\begin{array}{cc} 	
R_{11} & R_{12} \\
& 0\end{array}\right)\Pi^T}$. Then $\widehat{A}_k$ is a
  rank-$k$ matrix. It follows from equation~(\ref{eqn:RRQR})
  that 
\begin{eqnarray*}
{\displaystyle \sigma_j\left(\widehat{A}_k\right) } & \geq &{\displaystyle  \frac{\sigma_j}{
  \sqrt{1+f^2k(n-k)}}, \quad j = 1, \cdots, k,} \\
{\displaystyle \left\| A - \widehat{A}_k\right\|_2 } & \leq &
  {\displaystyle  \sigma_{k+1} \sqrt{1+f^2k(n-k)}} .
\end{eqnarray*}
These properties are compatible with the inequalities in
Theorem~\ref{Thm:HDBoundsIII}. The strong RRQR factorization
in Theorem~\ref{Thm:GuEisenstat_SISC96} also includes a
permutation $\Pi$ that selects $k$ linearly independent
columns of $A$ such that $\left\|R_{11}^{-1} R_{12}\right\|_2
\leq f$. Such information could be useful in some
applications~\cite{MartinssonRokhlinShkolniskyTygert}. 

But the matrix $QB_k$, being a two-sided orthogonal
approximation, does not contain any information about such
permutation. On the other hand, it is likely to be cheaper
to compute due to the matrix-matrix product operations
involved, and for rapidly decaying singular values or by potentially increasing the value of $q$, it
could make a much better approximation than $\widehat{A}_k$.

\section{Condition Number Estimation}\label{Sec:CondEst}
For any given square non-singular matrix $A$, define
\[\kappa(A) = \|A\| \|A^{-1}\| , \]
as its condition number. Here $\| \cdot \|$ is any matrix norm, such as the
matrix $1$-norm, $2$-norm, $\infty$-norm, Frobenius norm, or 
$\max$-norm. Condition numbers are of central importance in solving
many matrix computation problems, such as linear equations, least
squares problems, eigenvalue/eigenvector problems, and sparse matrix
problems. For a detailed discussion of condition number estimation,
see the survey paper by Higham~\cite{HIGHAM1} and the references
therein. More recent work includes Laub and Xia~\cite{LaubXia}.

A typical condition estimator uses a matrix norm estimator to estimate
$\|A\|$ and $\|A^{-1}\|$ separately, and multiply them together to get
an estimate for $\kappa(A)$. A typical matrix norm estimator, in turn,
only accesses the matrix $A$ through matrix-matrix or matrix-vector
multiplications, without the need to directly access entries of
$A$. Thus the costs of estimating $\|A\|$ and $\|A^{-1}\|$ are similar
if a factorization for $A$ is available. The goal in matrix norm
estimation is to compute a reliable estimate of $\|A\|$ up to a factor
that does not grow too fast with the dimension of $A$, perhaps without
direct access to entries of $A$, at a cost that is considerably less
than that of matrix factorization or inversion, something that is believed
to be impossible (see Remark~\ref{Rem:Fail}.)

However, by Corollary~\ref{Cor:HDBoundsIIIv2}, we know
Algorithm~\ref{Alg:RandSI} does compute a reliable estimate
for $\|A\|_2$ with $k=1$ and a reasonable choice of $\ell >
1$, due to the randomization of the start matrix. Below we
concentrate on estimating $\|A\|_1$. Currently, Hager's
method is one of the most popular 
estimators for $\|A\|_1$, is the default $1$-norm estimator of
LAPACK~\cite{LAPACK,HAGER,HIGHAM1,HIGHAM90A}. Hager's method is based
on a variant of the gradient descent method to find a local
maximizer for the following optimization problem:
\begin{equation}\label{Eqn:Hag}
{\displaystyle \left\| A \right\|_1 =
\max_{x\in {\cal S}} \left\|Ax \right\|_1, \quad \mbox{where} \quad
{\cal S} = \left\{ x \in {\bf R}^n: \|x\|_1 \leq 1. \right\}
} 
\end{equation}
\begin{algorithm}\label{Alg:Hager}{\bf Hager's Method} \\
\headerrule

\begin{tabular}{ll}
{\bf Input:} & $m \times n$ matrix $A$, and initial $1$-norm unit vector $x$. \\
{\bf Output:} & An estimate for $\|A\|_1$.  \\
\end{tabular}

\headerrule \\
{\noindent \bf repeat} 
\begin{INDENT}
    \begin{enumerate}
   \item Compute $y = A x, \quad z = A^T {\bf sign}(y)$. 
   \item {\bf if} $\|z\|_{\infty} \leq z^T x \;\; $ {\bf then} 
\begin{INDENT}
{\bf return} $\gamma = \|y\|_{1}$.
\end{INDENT}
   \item  $x = e_j, \quad \mbox{where} \quad j = {\bf argmax}_k |z_k|.$
   \end{enumerate}
\end{INDENT}
\headerrule
\end{algorithm}

The $e_j$ is the $j$-th unit vector. While it could
occasionally take much longer, Hager's method typically
takes very few (less than $5$) iterations to converge to a
local maximum that is within a reasonable factor (like $10$
or less) of $\|A\|_1$. As Algorithm~\ref{Alg:RandSI} already
computes a reliable estimate for $\|A\|_2$, it is
straightforward to combine Algorithms~\ref{Alg:RandSI}
and~\ref{Alg:Hager} to obtain a reliable estimate for
$\|A\|_1$, which satisfies $\|A\|_1 \geq \|A\|_2/\sqrt{n}$. 

\begin{algorithm}\label{Alg:Hager2}{\bf Randomized Hager's Method} \\
\headerrule

\begin{tabular}{ll}
{\bf Input:} & $m \times n$ matrix $A$, and integer $\ell > 1$. \\
{\bf Output:} & An estimate for $\|A\|_1$.  \\
\end{tabular}

\headerrule
\begin{enumerate}
\item Compute rank-$1$ approximation $QB_1$ to $A$ using
Algorithm~\ref{Alg:RandSI}
\item Set $\widehat{u}$ to be the right singular vector of
  $QB_1 $.  
\item Run Algorithm~\ref{Alg:Hager} on $A$ with
  initial vector $x = \widehat{u}/\|\widehat{u}\|_1$.
\item Return $\gamma$ from Algorithm~\ref{Alg:Hager}.
\end{enumerate}
\headerrule
\end{algorithm}

Since $QB_1$ is a rank-1 matrix, $\widehat{u}$ is
straightforward to compute. The number of iterations in
Algorithm~\ref{Alg:Hager} can be restricted to as few as
$1$ or $2$. This is because Algorithm~\ref{Alg:Hager} is
only used to find a
column whose vector $1$-norm provides the estimate for
$\|A\|_1$, no local maximum to problem~(\ref{Eqn:Hag}) is
necessary. Corollary~\ref{Cor:Hager2} directly follows from 
Corollary~\ref{Cor:HDBoundsIIIv2}.
\begin{corollary} \label{Cor:Hager2} For any $0 < \Delta
  \ll 1$, the output $\gamma$ from
  Algorithm~\ref{Alg:Hager2} must satisfy
\[ {\displaystyle \gamma } \geq
{\displaystyle \frac{ \|A\|_1}{\sqrt{n}\sqrt{1 + {\displaystyle
	\widehat{\cal C}_{\Delta}^2 }}}} \quad \mbox{where}
	\quad 
{\displaystyle \widehat{\cal C}_{\Delta} =
  \frac{e}{\sqrt{\ell}} \left(\frac{2}{\Delta}\right)^{\frac{1}{\ell}} 
\left(\sqrt{n}+\sqrt{\ell} + \sqrt{2 \log
  \frac{2}{\Delta}}\right), } \]
with exception probability at most ${\Delta} $.
\end{corollary} 

\begin{REMARK}\label{Rem:Hag2} One probably does not need to
choose a very tiny $\Delta$ for matrix norm estimation. In
our numerical experiments, $\ell = 5$ worked very well. For
matrices of dimension up to $200$,
Algorithm~\ref{Alg:Hager2} never under-estimated the true
norm by a factor over $10$. In general, we can choose 
${\displaystyle \ell = \lceil \log_{2}\left( \frac{2}{{\Delta} } \right) \rceil}$, in which case 
 the constants ${\displaystyle \widehat{\cal C}_{\Delta}}$ and $\gamma$ above satisfy
\[{\displaystyle \widehat{\cal C}_{\Delta} < 2 e \left(\sqrt{\frac{n}{\ell}}+3\right)}
\quad \mbox{and} \quad 
{\displaystyle \gamma } \geq
{\displaystyle \frac{ \|A\|_1}{2 e \sqrt{n}{\displaystyle \left(\sqrt{\frac{n}{\ell}}+4\right)}}}. \]
\end{REMARK}

\begin{REMARK}\label{Rem:Hag} 
Hager's method has been
generalized by Higham~\cite{highampnorm} to estimate the
matrix $p$-norm for any $p \geq 1$ and the mixed matrix norm
$\|A\|_{\alpha,\beta}$ for $\alpha \geq 1$ and $\beta \geq
1$. In particular, the $\max$-norm is the special case with $\alpha =
\infty$ and $\beta = 1$. Algorithm~\ref{Alg:Hager2} can be
trivially generalized to those cases as well, by replacing
Hager's method in Algorithm~\ref{Alg:Hager2} with its
generalized version, leading to a
Corollary~\ref{Cor:Hager2}-like conclusion for
reliability. We omit the details. 
\end{REMARK}

\begin{REMARK}
Kuczy\'{n}ski and Wo\'{z}niakowski~\cite{Kuczynski}
developed probabilistic error bounds for estimating the
condition number using the Lanczos algorithm for unit
start vectors under the uniform distribution. However, our
results appear to be much stronger.
\end{REMARK}

Below, we demonstrate the robustness of
Algorithm~\ref{Alg:Hager2} through the following
example. Let 
\[{\displaystyle A = \left(\begin{array}{cc} \alpha & b^T
    \cr b & \rho E\widehat{A} E\end{array}\right), \quad
    \mbox{for} \quad {\displaystyle E = I - 
\frac{1}{n-1} 
\left(\begin{array}{c} 1
    \cr \vdots \cr 1\end{array}\right)
\left(\begin{array}{c} 1
    \cr \vdots \cr 1\end{array}\right)^T}, } \]
where $\alpha>0, \rho > 0$ are scalars, $b > 0$ is an $n-1$
dimensional vector, and $\widehat{A}$ is
an $(n-1)\times (n-1)$ matrix. If we take the 
initial vector $x$ in Algorithm~\ref{Alg:Hager} to be the vector
of all $1$'s (the default choice in LAPACK), then
Algorithm~\ref{Alg:Hager} will always return $\alpha +
\|b\|_1$ as the $1$-norm estimate, regardless of
$\rho \widehat{A}$. 

In our numerical experiment, we set $n = 100$, $\rho =
10^{10}$ and chose $\alpha$, $b$ and $\widehat{A}$ to be
random, with $\|A\|_1 \approx 8.35 \times 10^{11}$. 
For $\ell = 5$, we obtained $\|A\|_1
\approx 2.46 \times 10^{11}$ from
Algorithm~\ref{Alg:Hager2}. On the other hand, 
Algorithm~\ref{Alg:Hager} returned $\|A\|_1 \approx 4.72
\times 10^{1}$, which was completely wrong.

\section{Numerical Experiments}\label{Sec:Num}
In this section we perform numerical experiments to shed
more light on randomized algorithms. Our main purpose of
these experiments is to provide numerical support to our
probabilistic analysis and to demonstrate that different
applications can lead to different singular value
distributions in the matrix and impose different accuracy
requirements, and thus demand different levels of
computational effort on the randomized algorithms.

\subsection{Improved Randomized Power Iteration} 

In the case of a small $k$, it seems unreasonable to require
a potentially large value of $p$ as suggested in
equation~(\ref{Eqn:p}). However, for a truely small value of
$p$, going random is still not enough to overcome the
potential problem of slow convergence associated with a poor
start matrix in Algorithm~\ref{Alg:RandSI}, and some
additional work maybe needed (see
Sections~\ref{Sec:StatAnal}.)

This discussion is particularly relevant for $k = 1$, which
corresponds to the classical power method,
Algorithm~\ref{Alg:BasicPM}, and its randomized version, 
Algorithm~\ref{Alg:RandPM}, in Appendix A. Any value of $p
>0$ seems to be too much work, but $p =0$ does not lead to
fast enough convergence. 

According to Corollary~\ref{Cor:HDBoundsIIIv2}, Algorithm~\ref{Alg:RandSI} can already compute
order of magnitude approximations to all the leading
singular values with $q = 0$. Thus, an obvious improvement of
Algorithm~\ref{Alg:RandSI} for small values of $k$ would be
to compute $\Omega$ with Algorithm~\ref{Alg:RandSamI} and then
compute a subspace approximation with
Algorithm~\ref{Alg:BasicSI}. Algorithm~\ref{Alg:RandSIAd2}
below is designed for subspace computations where
${\displaystyle k =O\left(\lceil \log_{10}\left(
\frac{2}{{\Delta} }\right) \rceil\right)}$ or smaller.

\begin{algorithm} \label{Alg:RandSIAd2}{\bf Improved
    Randomized Subspace Iteration for small $k$}\\
\headerrule

\begin{tabular}{ll}
{\noindent \bf Input:} & $m \times n$ matrix $A$ with $ n \leq m$, \\
& integers $q$ and $\ell_1 > \ell_2 \geq k$.  \\
{\noindent \bf Output:} & a rank-$k$ approximation.  \\
\end{tabular}

\headerrule
\begin{INDENT}
    \begin{enumerate} 
    \item Run Algorithm~\ref{Alg:RandSamI} with 
    $\ell = \ell_1$ for a rank-$\ell_2$ approximation. 
    \item Set $\Omega$ to be approximate right singular
    vector matrix. 
    \item Run Algorithm~\ref{Alg:BasicSI} with $\Omega$ and
    $\ell = \ell_2$ for a rank-$k$ approximation. 
    \end{enumerate}
\end{INDENT}
\headerrule
\end{algorithm}

We perform our experiments with $4000 \times 4000$ matrices 
of the form
\[ A =  \left(\log\left\|X_i-Y_j\right\|_2\right) , \]
where $\{X_i\}$ are $n$-dimensional Gaussian random
variables with mean $0$ and standard deviation $1$, and
where $\{Y_j\}$ are $n$-dimensional Gaussian random
variables with mean $\mu$ and standard deviation $1$. We
choose different $\mu$ values to control the ratio of the
two leading singular values of $A$. 

\begin{figure}[th!]
\vspace{-1.8in}
\begin{center}
{\includegraphics[width=.7\textwidth]{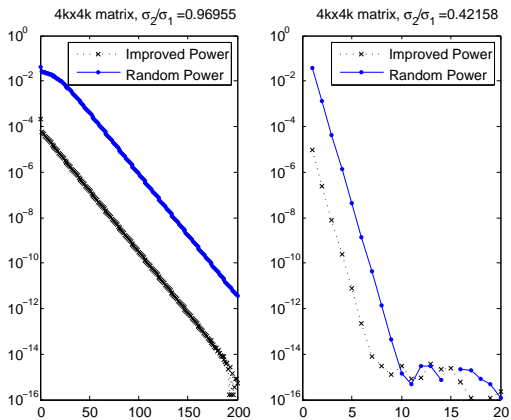}}
\label{fig:pi}
\vspace{-1.9in}
\caption{Faster Convergence of
  Algorithm~\ref{Alg:RandSIAd2} due to Better Choice of
  Start Vector}
\end{center}
\end{figure}

We ran Algorithm~\ref{Alg:RandSIAd2} with $\ell_1 = 5$ and
$\ell_2 = k = 1$. We also ran Randomized Power Method,
Algorithm~\ref{Alg:RandPM}, to compute $\|A\|_2$. We choose
$\mu = 1$ for a large $\sigma_2/\sigma_1$ ratio and  $\mu =
2.5$ for a small ratio. The results are summarized in
Figure~\ref{fig:pi}.  

For the case of large $\sigma_2/\sigma_1$ ratio,
Algorithm~\ref{Alg:RandSIAd2} converged to $\|A\|_2$ in
about $250$ steps, as opposed to about $350$ steps for
Algorithm~\ref{Alg:RandPM}. For the case of a small 
$\sigma_2/\sigma_1$ ratio, both algorithms performed
equally well. Algorithm~\ref{Alg:RandSIAd2} converged
slightly more quickly, but that is offset by the extra
work needed to compute the initial $\Omega$. 

Figure~\ref{fig:pi} confirms our analysis. At the cost of
the initial step to obtain a good start vector,
Algorithm~\ref{Alg:RandSIAd2} can converge significantly
faster than Algorithm~\ref{Alg:RandPM}. 

\subsection{low-rank approximation}
\begin{figure}[th!]
\begin{center}
{\includegraphics[width=0.55\textwidth]{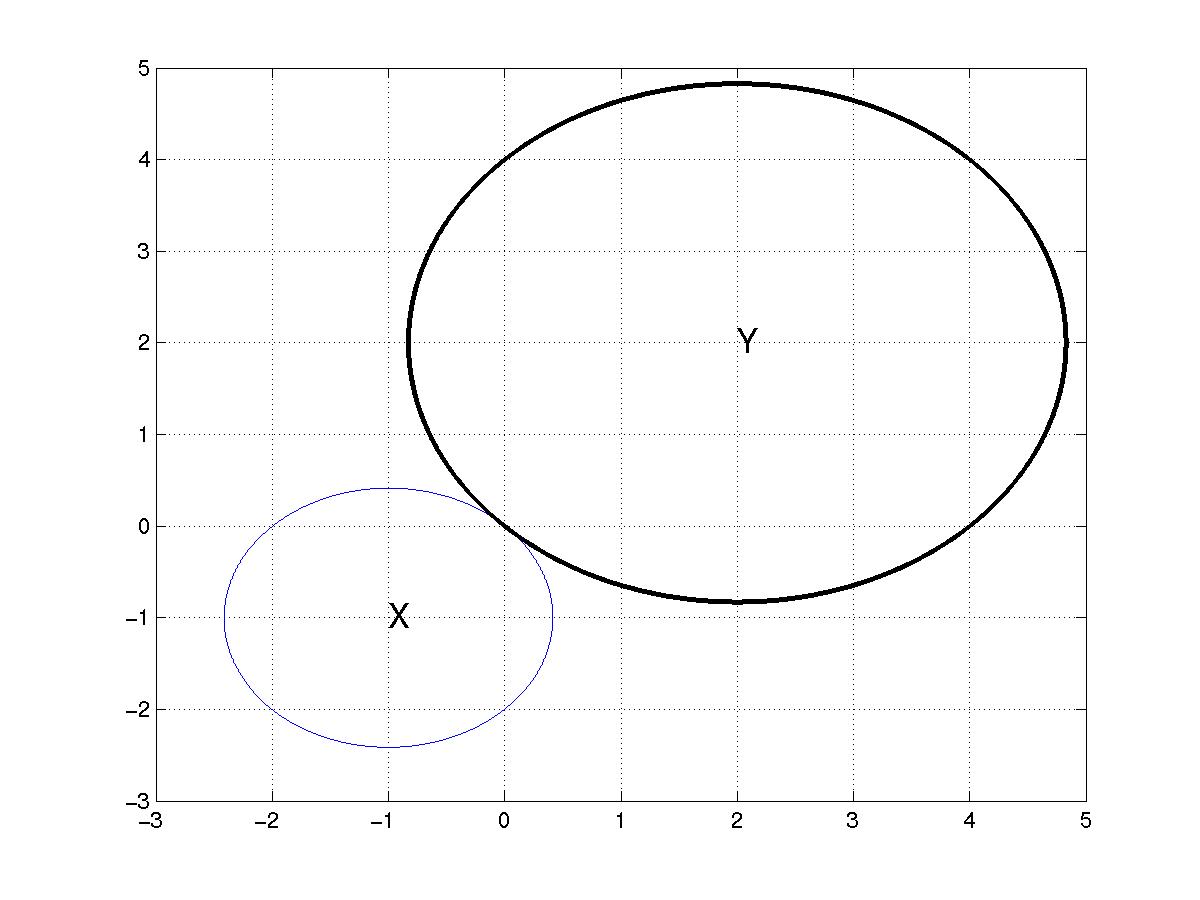}}
\caption{$X$ and $Y$ points in $A$}\label{fig:Dist}
\end{center}
\end{figure}

In this experiment, we consider a $4000 \times 4000$ matrix of
the form
\[ A =  \left(\log\left\|X_i-Y_j\right\|_2\right) , \]
where $\{X_i\}$ are equi-spaced points on the edge of the
disc $\|X-\pmatrix{-1\cr -1}\|_2 = \sqrt{2}$ and $\{Y_j\}$
equi-spaced points on the edge of the disc
$\|Y-\pmatrix{2\cr 2}\|_2 = 2\sqrt{2}$ (see
Figure~\ref{fig:Dist}.) We compare the performance of
Algorithms~\ref{Alg:RandSamI} and~\ref{Alg:RandSI}
against that of {\tt svds}, the matlab version of ARPACK~\cite{LehoucqSorensenYang} for finding a few selected 
singular values of large matrices. We choose $k = 50$. The
results are summarized in Table~\ref{Tab:t1}.

\begin{centering}
\begin{table}[th!]
\caption{\label{Tab:t1}
Numbers of Matrix-Vector Multiplies}
\begin{center}
\begin{tabular}{|l| l| c| c| c| c|}
\hline\hline
\headerruleA Tolerance & $q = 0$& $q = 2$ & $q = 4$ & svds \\
\hline
$10^{-6}$ & $143$ & $5\times 96$ & $9\times 79$ & $500$ \\
$10^{-8}$ & $180$ & $5\times 96$ & $9 \times 87$ & $600$ \\
$10^{-10}$ & $190$ & $5\times 96$ & $9\times 93$ & $600$ \\
\hline
\end{tabular}
\end{center}
\end{table}
\end{centering}

Since the singular values of this matrix decay relatively
quickly, Algorithm~\ref{Alg:RandSamI} seems to out-perform
Algorithm~\ref{Alg:RandSI} for any values of $q >
0$. Algorithm~\ref{Alg:RandSamI} also outperforms {\tt
svds}. As Algorithm~\ref{Alg:RandSamI} mostly computes
matrix-matrix products whereas each step of svds involves a
matrix-vector product, we would expect
Algorithm~\ref{Alg:RandSamI} to have even better
performance than {\tt svds} on modern serial and parallel
architectures.  This example demonstrates that for matrices
with fast decaying singular values, randomized algorithms
can be as competitive as the best methods for computing highly
accurate low-rank approximations.

\subsection{Structured Matrix Computations}

\begin{figure}[th!]
\begin{center}
{\includegraphics[width=.7\textwidth]{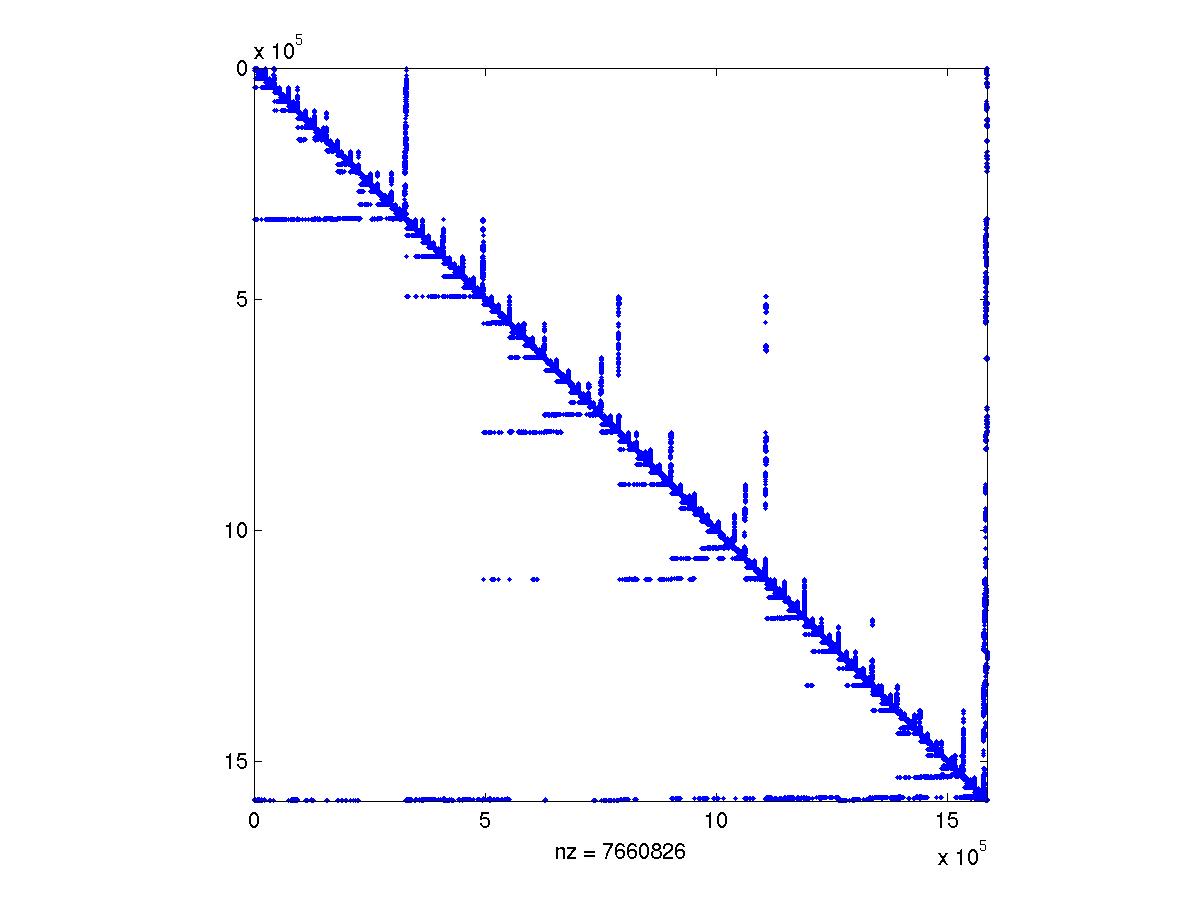}}
\label{fig:G3circuit}
\caption{A Sparse Matrix in Symmetric Minimum Degree Ordering}
\end{center}
\end{figure}

In this example, we demonstrate the effectiveness of
randomized algorithms for low-rank approximation in
the context of structured matrix computations. ${\tt G}3_{\tt circuit}$ is a $1585478 \times 1585478$ sparse SPD matrix
arising from circuit simulations. It is publicly available
in the {\em University of Flordia Sparse Matrix
Collection}~\cite{tdavis}. Figure~\ref{fig:G3circuit}
depicts its sparsity pattern in the symmetric minimum degree
ordering~\cite{georgeliu89}. A direct factorization of this
matrix creates a large amount of fill-in. In particular, the
Schur complement of the leading $1582178 \times 1582178$
principal submatrix, to be called $A$, is a 
$3300\times 3300$ dense submatrix. Here we compute hierarchical semiseparable (HSS) preconditioners
to $A$ with the techniques in~\cite{LiGuWuXia} and report
the numbers of preconditioned conjugate gradient (PCG) steps
to iteratively solve for a linear system of equations $A x=
b$ for a random right hand side $b$. The PCG is a very popular technique for solving
large SPD systems of equations~\cite{hest54,pcg}. We refer
the reader to~\cite{LiGuWuXia,martinsson2010} for 
details about the HSS matrix structure and its numerical
construction, but emphasize that the key and most
time-consuming step for computing HSS preconditioners is to
approximate various off-diagonal blocks of the matrix $A$ by
matrices of rank $k$ or less. We choose convergence
tolerance $\delta = 10^{-12}$. The conjugate gradient method
(CG) without any preconditioning takes $878$ iterations to reduce
the residual below this tolerance.

\begin{centering}
\begin{table}[th!]
\caption{\label{Tab:HSS}
Numbers of PCG Iterations}
\begin{center}
\begin{tabular}{|c| c| c| c| c|}
\hline\hline
\headerruleA Maximum off-diagonal rank $k$ & $p= 10$& $p = 20$ & $p = 40$  \\
\hline
$20$ & $75$ & $77$ & $72$   \\
$40$ & $69$ & $69$ & $69$   \\
$60$ & $64$ & $61$ & $61$   \\
\hline
\end{tabular}
\end{center}
\end{table}
\end{centering}

Table~\ref{Tab:HSS} summarizes our results. We can see that
all choices of $p$ drastically decrease the number of CG
iterations. Howver, the additional reduction in the number
of CG iterations is typically small for higher values of
$p$. Considering the extra cost involved in higher $p$
values in the construction of HSS preconditioners, it seems
that higher $p$ values are ineffective for this
application. This example suggests that for the purpose of
constructing preconditioners in structured matrix
computations, a small $p$ value is typically sufficient to
develop highly effective preconditioners. This is consistent
with the rule of thumb that typically randomized algorithms
require very little oversampling and a value of $p$ in
between $10$ to $20$
suffices~\cite{martinsson2010,MartinssonRokhlinShkolniskyTygert,RokhlinTygert}. 
In fact, the SVD truncation in
Algorithm~\ref{Alg:RandSI} and Algorithm~\ref{Alg:BasicSI} is unnecessary for this example. 

\subsection{Eigenfaces}

\begin{centering}
\begin{table}[th!]
\caption{\label{Tab:Eigenfacesmatch}
Comparison of Numbers of Incorrect Matches}
\begin{center}
\begin{tabular}{|l|  c| c| c| c|}
\hline\hline
\headerruleA Rank $k$ & $p = 10$& $p = 20$ & $p = 40$ & Truncated SVD \\
\hline
$10$ & $32$ & $25$ & $23$ & $24$   \\
$20$ & $25$ & $26$ & $25$ & $21$   \\
$30$ & $21$ & $20$ & $18$ & $17$   \\
$40$ & $20$ & $17$ & $17$ & $16$   \\
\hline
\end{tabular}
\end{center}
\end{table}
\end{centering}

Eigenfaces is a well studied method of face recognition
based on principal component analysis (PCA), popularised by
the seminal work of Turk and Pentland~\cite{turkpent}. For
more recent work and survey,
see~\cite{BrunelliPoggio,KirbySirovich,SinhaBalas,SirovichKirby,SirovichMeytlis}
and the references therein. In this experiment we
demonstrate the effects of randomized algorithms on face
recognition. 

\begin{figure}[th!]
\vspace{-2.5in}
{\includegraphics[width=1.\textwidth]{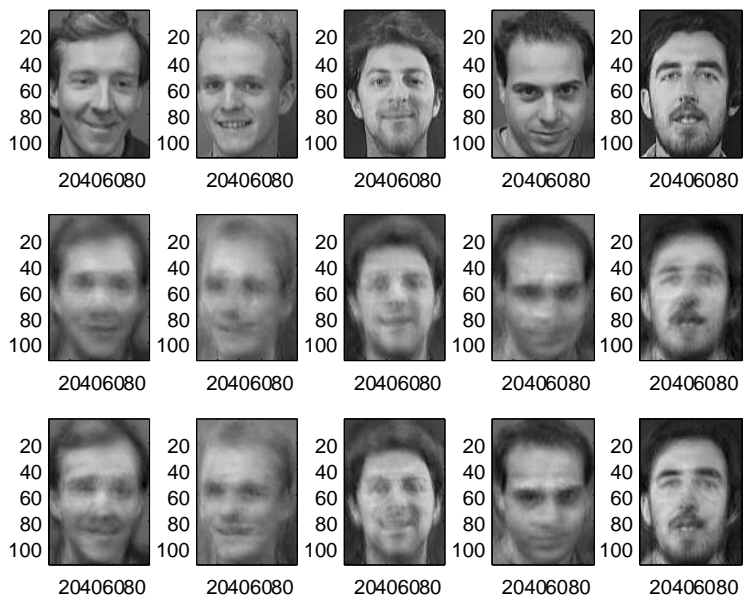}}
\label{fig:Eigenfaces}
\vspace{-2in}
\caption{Original Faces and their Eigenfaces}
\end{figure}

Typical face recognition starts with a data base of training
images, which are then processed as follows:
\begin{enumerate}
\item Calculate the mean of the training images.
\item Subtract the mean from the training images, obtaining the
mean-shifted images.
\item Calculate a truncated SVD of the mean-shifted images. 
\item Project the mean-shifted images into the singular
  vector space using the retained singular vectors,
  obtaining feature vectors. 
\end{enumerate}

To classify a new face, one does the following calculations:
\begin{enumerate}
\item Subtract the mean from the new image, obtaining the
mean-shifted image.
\item Project the mean-shifted image into the singular
  vector space, obtaining a new feature vector. 
\item Find the feature vector in the data base that best
  matches the new feature vector. 
\end{enumerate}

Our face data are obtained from the Database of Faces
maintained at the AT\&T Laboratories
Cambridge~\cite{attfaces}. All faces are greyscale images
with a consistent resolution.  There are ten different
images of each of 40 distinct subjects. The size of each
image is $92\times 112$ pixels, with $256$ grey levels per
pixel. We use 200 of these images, 5 from each individual,
as training images, and the remaining ones for
classification.

In Figure~\ref{fig:Eigenfaces}, the first row are the
original face images; the second row are eigenfaces with a
rank-$10$ truncated SVD, and the third row eigenfaces with a
rank-$20$. 

In addition to the exact truncated SVD, we also perform
image training and classification using
Algorithm~\ref{Alg:RandSamI} with different $p$ values. The
results are summarized in
Table~\ref{Tab:Eigenfacesmatch}. It is clear that smaller
$p$ values give worse results than truncated SVD, but $p =
40$ gives results that are very similar to truncated SVD,
even though some of the singular values are accurate to only
within $1$ to $2$ digits. This example demonstrates that
limited accuracy that goes beyond being correct to within a
constant factor is sufficient for some applications. 
 
\section{Conclusions and Future Work} \label{Sec:Conc}
We have presented some interesting results on randomized
algorithms within the framework of the subspace iteration
method for singular value and low-rank matrix
approximations. While randomized algorithms have been
primarily considered as an efficient tool to compute
low-rank approximations, our results further suggest that
they actually compute the much stronger rank-revealing
factorizations, and can be used to reliably estimate
condition numbers. We have also presented numerical
experimental results that support our analysis.

This work opens up many directions for future research. Most immediate
is the convergence analysis on singular vectors. We expect results
compatible to those for singular values. Variations of subspace
iteration methods exist for computing eigenvalues of symmetric and
non-symmetric matrices. It would be interesting to extend our results
to these methods. Additionally, for a given rank $k$, an interesting
and important issue is how to choose the parameters $q$ and $\ell$ to
minimize the overall cost in Algorithm~\ref{Alg:RandSI} and
Algorithm~\ref{Alg:BasicSI}.
 
{\noindent \bf Acknowledgments.} The author would like to thank
Shengguo Li, Michael Mahoney, Vladimir Rokhlin, Mark Tygert, Jianlin Xia and Chao Yang
for many helpful discussions on this subject. He would especially like
to thank Joel Tropp, whose interesting talk at UC Berkeley in the
Spring of 2010 sparked the author's interest on the subject that
eventually led to this work, and Chris Melgaard, with whom he had
extensive discussions about the material presented in this
work. Finally, the author would like to thank the anonymous referees
who go out of their ways to provide numerous helpful suggestions that
greatly improved the presentation of this paper, including a
shorter proof for Theorem~\ref{Thm:Fto2}.

\appendix 
\renewcommand{\theequation}{A.\arabic{equation}}
\renewcommand{\thetheorem}{A.\arabic{theorem}}
\renewcommand{\thesection}{A}

{\bf Appendix.} For numerical stability, Algorithm~\ref{Alg:B} below is
often performed once every few iterations in subspace
iteration methods, to balance efficiency and numerical
stability (see Saad~\cite{saad2}.) 

\begin{algorithm} \label{Alg:B}{\bf Orthorgonalization with QR}\\
\headerrule

\begin{tabular}{ll}
{\bf Input:} & $m \times n$ matrix $A$, $n \times
\ell$ start matrix $\Omega$, and integer $q \geq 0$. \\
{\bf Output:} & $Y = \left(A A^{T}\right)^q A \Omega$.  \\
\end{tabular}

\headerrule
\begin{INDENT}
    Compute $Y = A \Omega$, and QR factorize $Q R = Y$. \\
    \FOR/ $i = 1, \cdots, q$ \DO/
        $Y = A^{T} \, Q$; QR factorize $Q R = Y$; \\
        $Y = A \, Q$;  QR factorize $Q R = Y$.    
    \ENDFOR/
\end{INDENT}
\headerrule
\end{algorithm}

Below is the classical power method for computing the
$2$-norm of a given matrix.

\begin{algorithm} \label{Alg:BasicPM}{\bf Basic Power Method}\\
\headerrule

\begin{tabular}{ll}
{\noindent \bf Input:} & $m \times n$ matrix $A$ with $ n \leq
m$, \\
& and $n \times 1$ start vector $\Omega$. \\
{\noindent \bf Output:} & approximation to $\|A\|_2$. \\
\end{tabular}

\headerrule
\begin{INDENT}
    \begin{enumerate}
   \item Compute $Y = \left(AA^{T}\right)^q A \, \Omega$.
   \item Compute an orthogonal column basis $Q$ for $Y$. 
   \item Compute $B =  Q^{T} \, A $. 
   \item Return $\|B\|_2$.
    \end{enumerate}
\end{INDENT}
\headerrule
\end{algorithm}

In situations where no useful information about the leading
right singular vector is available, the vector $\Omega$ in
Algorithm~\ref{Alg:BasicPM} can also be chosen to  be
random, to enhance convergence, leading to 
\begin{algorithm} \label{Alg:RandPM}{\bf Randomized Power Method}\\
\headerrule

\begin{tabular}{ll}
{\noindent \bf Input:} & $m \times n$ matrix $A$ with $ n \leq
m$, \\
{\noindent \bf Output:} & approximation to $\|A\|_2$. \\
\end{tabular}

\headerrule
\begin{INDENT}
    \begin{enumerate}
   \item Draw a random $n \times 1$ vector $\Omega$. 
   \item Compute $Y = \left(AA^{T}\right)^q A \, \Omega$.
   \item Compute an orthogonal column basis $Q$ for $Y$. 
   \item Compute $B =  Q^{T} \, A $. 
   \item Return $\|B\|_2$.
    \end{enumerate}
\end{INDENT}
\headerrule
\end{algorithm}

\newpage

\begin{center}
{\bf Supplemental Material}
\end{center}

\setcounter{page}{1}
\pagenumbering{roman}

\renewcommand{\theequation}{S1.\arabic{equation}}
\renewcommand{\thetheorem}{S1.\arabic{theorem}}
\renewcommand{\thesection}{S1}

\section{Introduction} In the interest of reducing the length of the original paper,
we have put some of the non-essential material here. This {\em Supplemental Material} 
is organized as follows: In Section~\ref{Sec:Conv} we
discuss how the decaying rates of the singular values can
affect parameter choices in the randomized algorithms; in Section~\ref{Sec:Num} we
present additional supporting numerical experimental results; in
Section~\ref{Sec:Anal} we provide proofs for the two propositions in the original paper; and
in Section~\ref{Sec:ApPrelim} we list the facts we have used from calculus. 

\renewcommand{\theequation}{S2.\arabic{equation}}
\renewcommand{\thetheorem}{S2.\arabic{theorem}}
\renewcommand{\thesection}{S2}

\section{Further Convergence Considerations}\label{Sec:Conv}
In this section we discuss in more detail on the convergence
rate of Algorithm~\ref{Alg:RandSI}. We identify a singular
value distribution where Algorithm~\ref{Alg:RandSamI} is
likely to perform better, and another distribution where a
proper choice of $\ell$ leads to significant convergence
speedup. 

In the remainder of this section, we will always assume $p
\geq 2$ has been chosen according to equation~(\ref{Eqn:p}).
Furthermore, we will only consider the costs of computing
the matrix-vector products $A u$ and $A^{T} v$, which are
$O(mn)$ flops apiece. This is justified as the truncated SVD
only costs $O(n\ell^2)$ flops, smaller than the total
matrix-vector product cost, which is $O(mn\ell)$ flops, and
$\ell \ll n$. As per Section~\ref{Sec:StatAnal}, we will
concentrate on the expression
\[{\displaystyle \left(\frac{\sigma_{\ell - p
      +1}}{\sigma_k}\right)^{2q+1}} \]
as the key factor that controls singular value convergence. 

\subsection{Rapidly decaying singular value distributions}
First we consider the model where the singular values (except for the
first few) decay and satisfy the following equation
\begin{equation}\label{Eqn:DecayModel1}
\sigma_{s+t}  \leq \alpha \sigma_{s} \sigma_{t}, 
\end{equation}
for some constant $\alpha > 0$ and any $s, t > 1$. This model is
satisfied when the singular values of $A$ decay exponentially or
faster. We wish to show that Algorithm~\ref{Alg:RandSamI}
performs better than Algorithm~\ref{Alg:RandSI} with $q > 0$.  

Since Algorithm~\ref{Alg:RandSI} performs $\widehat{\ell} =
(2q+1) \ell$ matrix-vector products, we will thus allow 
Algorithm~\ref{Alg:RandSamI} to performs $\widehat{\ell}$ 
matrix-vector products as well. In this setting, both
algorithms cost roughly the same, and we will compare their
singular value convergence rates. For Algorithm~\ref{Alg:RandSI}, the ratio is
\[{\displaystyle \left(\frac{\sigma_{\ell - p +1}}{\sigma_k}\right)^{2q+1} \leq 
\left(\alpha \sigma_{\ell - p - k+1} \right)^{2q+1} },  \]
according to the singular value decay model~(\ref{Eqn:DecayModel1}).
On the other hand, for Algorithm~\ref{Alg:RandSamI}, the ratio is
\begin{eqnarray*}
{\displaystyle \left(\frac{\sigma_{\widehat{\ell} - p +1}}{\sigma_k}\right) }&\leq &
{\displaystyle \alpha \sigma_{ (2q+1) \ell - p - k+1} = 
\alpha \sigma_{ (2q+1) (\ell - p - k+1) + 2q(p+k-1)}
\leq  \alpha^2 \sigma_{ (2q+1) (\ell - p - k+1)}
\sigma_{2q(p+k-1)}} \\
&\leq & {\displaystyle \alpha^2 \left(\alpha \sigma_{\ell - p - k+1} \right)^{2q+1} \sigma_{2q(p+k-1)} 
},  
\end{eqnarray*}
which is a tighter upper bound. This comparison suggests that in
general there is little convergence advantage of
Algorithm~\ref{Alg:RandSI} over Algorithm~\ref{Alg:RandSamI} when
the singular values decay exponentially or faster.

\subsection{Slowly decaying singular value distributions}
Below we consider the model where the singular values (except for the
first few) decay and satisfy the following equation
\begin{equation}\label{Eqn:DecayModel2}
{\sigma_{s} \over \sigma_{k}}  \leq \alpha \left(k \over s\right)^T, 
\end{equation}
for some constants $\alpha \geq 1,  T \geq 1$ and any $s > k$. This model is
satisfied when the singular value $\sigma_{s}$ decays like $O(1/s^T)$.
With this model, we now analyze the performance of
Algorithm~\ref{Alg:RandSI} and derive a different set of choices
of the parameters $q$ and $\ell$. Again we allow
Algorithm~\ref{Alg:RandSI} to do a total of $\widehat{\ell}
= (2q+1) \ell$ matrix-vector products, and we will choose
$q$ and $\ell$ to optimize convergence. The ratio becomes 
\[{\displaystyle \left(\frac{\sigma_{\ell - p +1}}{\sigma_k}\right)^{2q+1} 
= \left(\frac{\sigma_{\ell - p +1}}{\sigma_k}\right)^{\widehat{\ell}/\ell} 
\leq \left(\alpha \left({ k \over \ell - p +1}\right)^T  \right){^{\widehat{\ell}/\ell}} 
\leq \left(\alpha \left({ k \over \ell - p +1}\right)^{1/\ell} \right)^{T\widehat{\ell}}} . \]
The last expression allows us to study the optimal choice of
$\ell$ that minimizes it. Define
\[f(\ell) = \left({ k \over \ell - p +1}\right)^{1/\ell}. \]
Then
\[f'(\ell) = - f(\ell) {g(\ell) \over \ell^2} , \quad
\mbox{where} \quad g(\ell) = { \ell \over \ell - p +1} +
\log {k \over \ell - p + 1} . \]
The optimal $\ell$, denoted $\ell_{\bf opt}$, is the unique solution
of $g(\ell) = 0$ and satisfies
\[ e k \leq \ell_{\bf opt} - p + 1 \leq \eta ((p-1)+ k) ,\]
where $\eta = 3.59 \cdots $ satisfies $1 + 1/\eta = \log \eta.$ This
range of $\ell_{\bf opt}$ suggests that Algorithm~\ref{Alg:RandSI}
should be used only when $\widehat{\ell} \gg \ell_{\bf opt}$. The benefit of
choosing $\ell = \ell_{\bf opt}$ is that the convergence rate for 
Algorithm~\ref{Alg:RandSI} now becomes exponential:
\[{\displaystyle \left(\frac{\sigma_{\ell - p +1}}{\sigma_k}\right)^{2q+1} 
= \left(\frac{\sigma_{\ell_{\bf opt} - p
+1}}{\sigma_k}\right)^{\widehat{\ell}/\ell_{\bf opt}} =
\left(e^{-1/\left(\ell_{\bf opt} - p
+1\right)}\right)^{\widehat{\ell}} \leq
\left(e^{-1/\left(\eta(p -1 + k)\right)}\right)^{\widehat{\ell}}}. \] 
The last expression is ${\displaystyle
e^{-O\left(\widehat{\ell}/k\right)}}$. In contrast, a naive
choice of $\ell = k + p$, would lead to a ratio of
${\displaystyle \left(\frac{\sigma_{\ell - p
+1}}{\sigma_k}\right)^{2q+1}} \leq {\displaystyle
e^{-O\left(\widehat{\ell}/k^2\right)}}$, slower than the
optimal one by a factor of $O(k)$ in the exponent.  Hence
the best choice of $\ell_{\bf opt}$ can lead to significantly
accelerated rate of overall convergence. 

\subsection{Adaptive Randomized Algorithms} \label{Sec:Adap}
From the two different singular value distributions
discussed above, it is clear that much
research is needed to design an efficient algorithm that can
automatically choose the right set of parameters for different singular
value distributions within the framework of
Algorithm~\ref{Alg:RandSI}.

In this section, we will limit our scope and present
an adaptive version of Algorithm~\ref{Alg:RandSI}, with the
assumption that the singular values decay slowly. Our goal
is to quickly compute rank-$k$ approximations up to the
tolerance provided. This algorithm is motivated by similar
work in~\cite{Tropp} and will be used later on in our
numerical experiments in Section~\ref{Sec:LSI}.

\begin{algorithm} \label{Alg:RandSIAd}{\bf : Adaptive Randomized Subspace Iteration Method}\\
\headerrule

\begin{tabular}{ll}
{\noindent \bf Input:} & $m \times n$ matrix $A$ with $ n \leq m$, 
accuracy tolerance $\tau > 0$, failure
tolerance $\Delta > 0$, \\
& integers $b, c$, $k$, $q>0$,  and ${\bf Cmax} > 0$.  \\
{\noindent \bf Output:} & a rank-$k$ approximation.  \\
\end{tabular}

\headerrule
\begin{INDENT}
    \begin{enumerate} 
    \item Compute ${\displaystyle p = \left\lceil
          \log_{10}\left( \frac{2}{\Delta} \right)\right \rceil}, \ell = c + k + p$. 
    \item Draw a random $n \times \ell $ test matrix $\Omega$. 
   \item Compute $Y = \left(A A^{T}\right)^q A \, \Omega$.
   \item Compute an orthogonal column basis $Q$ for $Y$. 
   \item Compute $B =  Q^{T} \, A $. 
   \item Compute the SVD of $B$ and  the rank-$k$ truncated SVD $B_k$.
   \item ${\displaystyle {\cal E} = \left(\frac{\sigma_{\ell - p +1}}{\sigma_k}\right)^{2q+1}}$.\\
   \WHILE/ ${\cal E} > \sqrt{\tau}$  \DO/
        \begin{enumerate} 
            \item Compute $\delta \ell = \max\left(b,\left \lceil
            \left(
            \left({\displaystyle\frac{\ell-p+1}{k}}\right)^{\displaystyle
            \frac{\log\sqrt{\tau}/{\cal E}}{\log \cal E}}-1\right)\left(\ell - p +1\right) \right\rceil\right)$. 
            \item \IF/ $(\delta \ell + \ell) > {\bf Cmax}$ \DO/
                       {\bf Quit}. {\sc Sampling size exceeding limit for the given tolerance} 
                  \ENDIF/
           \item Draw a random $n \times \delta \ell $ test matrix $\Omega$. 
           \item Update $Y = [Y \quad  \left(A A^{T}\right)^{q} A \, \Omega]$.
           \item Update the orthogonal column basis $Q$ for $Y$. 
           \item Update $B =  Q^{T} \, A $. 
           \item Update the SVD of $B$ and the rank-$k$ truncated SVD $B_k$.
           \item ${\cal E} = \left(\frac{\sigma_{\ell - p +1}}{\sigma_r}\right)^{2q}\sigma_{\ell-p + 1}$.
       \end{enumerate}
   \ENDWHILE/
    \item Return $QB_k$.
    \end{enumerate}
\end{INDENT}
\headerrule
\end{algorithm}
\begin{REMARK}
In Algorithm~\ref{Alg:RandSIAd}, an integer $c \geq 0$ was introduced to allow additional
initial convergence. Its value should be dependent on accuracy
tolerance $\tau$.  The integer $b > 0$ was introduced so at least $b$
columns will be sampled for each iteration. This is to avoid sampling
too few columns per iteration, as computations with too few columns
often incur additional data movement costs that slow down the
execution of the whole algorithm. The formula for $\delta
\ell$ was derived under the
assumption~(\ref{Eqn:DecayModel2}) with $\alpha = 1$. 
\end{REMARK}

\begin{REMARK}
The huristic choice of $\delta \ell$ in
Algorithm~\ref{Alg:RandSIAd} is aimed at reducing the
singular value error to $\tau$, under the assumption that
the singular values decay at least as fast as
model~(\ref{Eqn:DecayModel2}). 
\end{REMARK}

\renewcommand{\theequation}{S3.\arabic{equation}}
\renewcommand{\thetheorem}{S3.\arabic{theorem}}
\renewcommand{\thesection}{S3}

\section{Numerical Experiments}\label{Sec:LSI}.
In this section we report more numerical experimental results to shed
more light on randomized algorithms. 

Latent Semantic Indexing (LSI) is a massive data
processing application based on low-rank
approximations~\cite{BerryDumaisOBrien}. A data
base of terms and documents is processed to generate a
term-document matrix, where each column is a document with
each non-zero in the column represents the weighted number
of matches to a particular term.

\begin{centering}
\begin{table}[th!]
\caption{\label{Tab:LSI}
Number of Agreements with Truncated SVD}
\begin{center}
\begin{tabular}{|l|  c| c| c| }
\hline\hline
\headerruleA Tolerance $\tau$ & $q = 0$& $q = 2$ & $q = 4$ \\
\hline
$10^{-3}$  & $460$ & $709$ & $769$    \\
$10^{-7}$  & $511$ & $889$ & $910$    \\
$10^{-11}$ & $504$ & $909$ & $921$    \\
\hline
\end{tabular}
\end{center}
\end{table}
\end{centering}

Given a set of terms (a query), LSI attempts to find the
document that best matches it in some semantical sense. To
do so, LSI computes a rank-$k$ truncated SVD of the
term-document matrix so that $A \approx U_kS_k V_k^{T}$.

For any query vector $q$, compute the feature vector $d =
\left(q^{T} U\right)S_k^{-1}$. The document that most
matches $q$ is the row of $V_k$ that is the most parallel to
$d$.

We use the TDT2 text data~\cite{TDT2}. The TDT2 corpus 
consists of data collected during the first half of 1998 and
taken from 6 sources, including 2 newswires (APW, NYT), 2
radio programs (VOA, PRI) and 2 television programs (CNN,
ABC). It consists of 11201 on-topic documents which are
classified into 96 semantic categories. What is available
at~\cite{TDT2} is a subset of this corpus, with a total of
9,394 documents and over $36000$ terms.

\begin{figure}[th!]
\begin{center}
{\includegraphics[width=.75\textwidth]{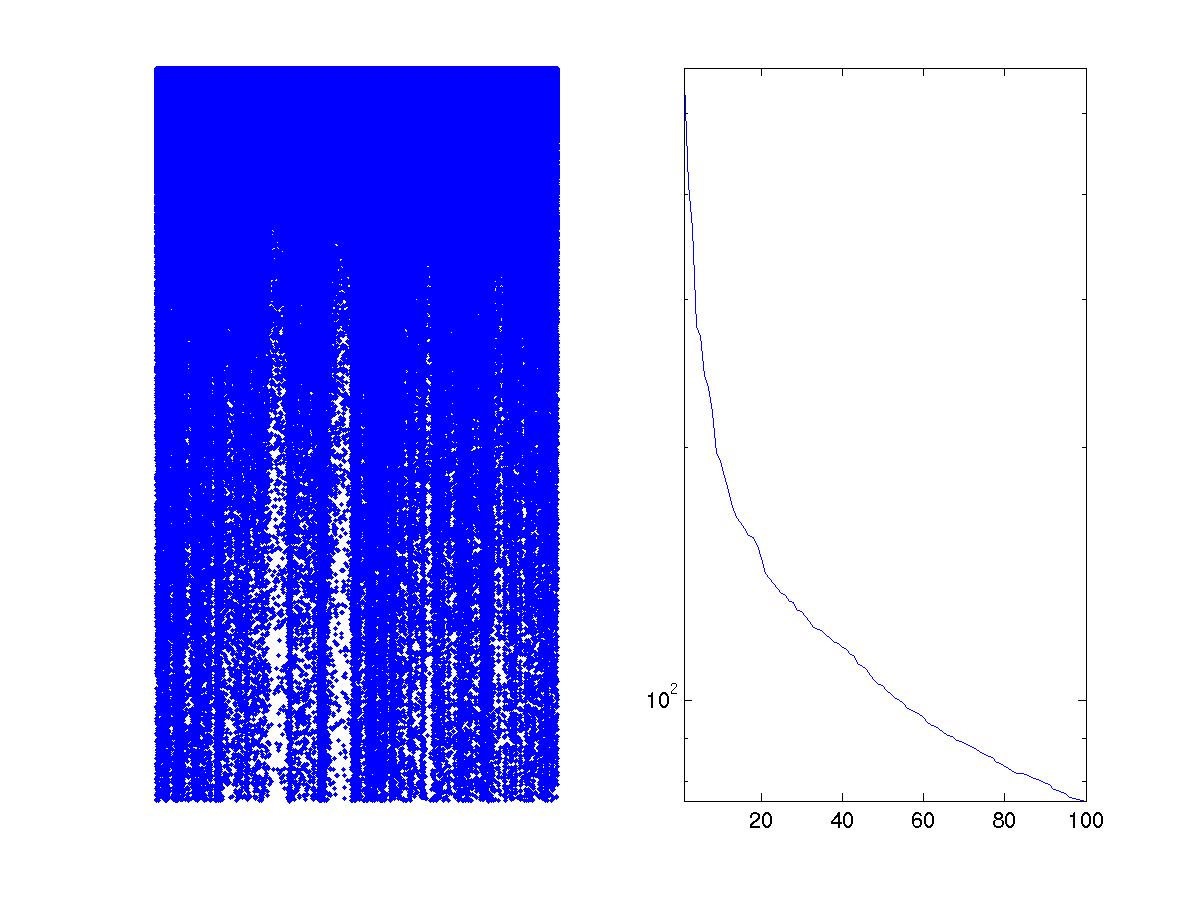}}
\label{fig:lsi}
\vspace{-.25in}
\caption{From left to right: Sparsity pattern of
  term-document matrix; Its leading Singular Values.}
\end{center}
\end{figure}

We performed $1000$ random queries with the truncated SVD
for different values of $k$. Then we repeat the same queries
with the low-rank approximation computed by
Algorithm~\ref{Alg:RandSIAd} for $q = 0, 2, 4$ and 
a decreasing set of $\tau$ values. For each $q$ and $\tau$, 
Algorithm~\ref{Alg:RandSIAd} automatically stops once
$500$ column samples have been reached in computing the
low-rank approximation.

Table~\ref{Tab:LSI} clearly indicates that better accuracy
in randomized algorithms leads to more agreement with the
truncated SVD in terms of query matches. Due to the nature
of this experiment, an agreement does not always mean a
better match. However, Table~\ref{Tab:LSI} does give some
indication that better accuracy in the low-rank
approximation is probably better for LSI. Since $q = 2$
looks significantly better than $q = 0$, this example
indicates that for LSI, it may be necessary to use
Algorithm~\ref{Alg:RandSIAd} with a small but positive
$q$ value for best performance.

\renewcommand{\theequation}{S4.\arabic{equation}}
\renewcommand{\thetheorem}{S4.\arabic{theorem}}
\renewcommand{\thesection}{S4}

\section{Proofs of Propositions~\ref{Lem:Omega2Bound}
and~\ref{Lem:Omega1Bound}}\label{Sec:Anal}

We begin with the following probability tool.  
\begin{lemma} (Chen and Dongarra~\cite{ChenDongarra})
  \label{Lem:ChenDon} Let $G$ be an $m
\times n$ standard Gaussian random matrix with $m \leq n$,
and let $f(x)$ denote the probability density function of
$\|G^{\dagger}\|_2^{-2}$, then $f(x)$ satisfies:
\[{\displaystyle  f(x) \leq L_{m,n} e^{-\frac{x}{2}}
  x^{\frac{1}{2}(n-m-1)}, \quad \mbox{where} \quad 
 L_{m,n} = \frac{{\displaystyle 2^{\frac{n-m-1}{2}}
    \Gamma\left(\frac{n+1}{2}\right)}}{{\displaystyle \Gamma\left(\frac{m}{2}\right)}}
    \Gamma\left(n - m + 1\right). }\]
\end{lemma}

The following classical result, the {\em law of the
  unconscious statistician}, will be very helpful to our analysis.
\begin{prop}\label{prop:law} Let $g(\cdot)$ be a
  non-negative continuously differentiable function with $g(0) = 0$, and let
  $G$ be a random matrix, we have   
\[{\displaystyle  \mathbb{E} g\left(\|G \|_2 \right) = \int_{0}^{\infty}
  g'\left(x \right) \mathbb{P} \left\{\|G \|_2 \geq x
  \right\} dx . }\]
\end{prop}

We also need to define the following functions
\begin{equation}\label{Eqn:gghat} 
g(x)  =  {\displaystyle 1 - \frac{1}{\sqrt{1 + \alpha^2
      x^2}}} \quad \mbox{and} \quad 
\widehat{g}(x) = \sqrt{\delta^2 + \frac{\alpha^2
      x^2}{\beta^2 + \gamma^2x^2}}-\delta,
\end{equation}
where $\alpha > 0$ and $\delta > 0$ are constants to be specified later on. 
It is easy to see tht $g(0) = 0$ and $\widehat{g}(0) = 0$, and 
\[ {\displaystyle g'(x) = \frac{\alpha^2 x}{\left(\sqrt{1 + \alpha^2
x^2}\right)^3} \quad \mbox{and} \quad 
\widehat{g}'(x) = \frac{\alpha^2 \beta^2 x}{\left(\beta^2 + \gamma^2x^2\right)^2
\sqrt{\delta^2 + \frac{\alpha^2  x^2}{\beta^2 +
    \gamma^2x^2}}}}. \]

{\noindent \bf Proof of Proposition~\ref{Lem:Omega2Bound}:}  
Define a function $h(G) = \|G\|_2$. Then by Proposition~\ref{Lem:FnormExpV}, we have 
\begin{equation}\label{Eqn:EhG}
 \mathbb{E} \left(h(G)\right) \leq \sqrt{m}+\sqrt{n} < \sqrt{m}+\sqrt{n} + 3 \stackrel{def}= {\cal E}. 
\end{equation}
$h$ is a Lipschitz function on matrices with Lipschitz constant ${\cal L} = 1$ (see
Theorem~\ref{Thm:Gaussfun}):
\[|h(X)-h(Y)| \leq \|X-Y\|_F \quad \mbox{for all} \quad X, Y.  \]
For equation~(\ref{Eqn:SVAB1}), we can rewrite, by way of
function $g(x)$ in~(\ref{Eqn:gghat}) and Proposition~\ref{prop:law}, 
\[ {\displaystyle \mathbb{E} \left( \frac{1}{\sqrt{1+\alpha^2 \|G\|_2^2}}\right) } =  {\displaystyle  
1 - \mathbb{E} \left( g\left(\|G\|\right)\right) = 1 - \int^{\infty}_{0} g'(x) 
\mathbb{P} \left\{\|G \|_2 \geq x
  \right\} dx . }\]
By Theorem~\ref{Thm:Gaussfun}, we have ${\displaystyle
\mathbb{P} \left\{\|G \|_2 \geq x \right\} \leq e^{-u^2/2}}$
for $u = x - {\cal E}$. Putting it all together,
\begin{eqnarray}  
 {\displaystyle \mathbb{E} \left( \frac{1}{\sqrt{1+\alpha^2 \|G\|_2^2}}\right) } & \geq &  
 {\displaystyle 1 - \int^{{\cal E}}_{0} g'(x) dx - \int^{\infty}_{{\cal E}} g'(x) \mathbb{P} \left\{\|G \|_2 \geq x
  \right\} dx } \nonumber \\
& \geq & {\displaystyle  \frac{1}{\sqrt{1 + \alpha^2 {\cal E}^2}} - 
\int^{\infty}_{{\cal E}} \frac{\alpha^2 x }{\left(\sqrt{1 + \alpha^2 x^2}\right)^3} e^{-(x - {\cal E})^2/2} dx }  \nonumber \\
& \geq & {\displaystyle  \frac{1}{\sqrt{1 + \alpha^2 {\cal E}^2}} - 
\frac{\alpha^2}{\left(\sqrt{1 + \alpha^2 {\cal E}^2}\right)^3} \int^{\infty}_{0} \left({\cal E} + u\right) e^{-u^2/2} d u  }  \nonumber \\
& = & {\displaystyle  \frac{1}{\sqrt{1 + \alpha^2 {\cal E}^2}} - 
\frac{\alpha^2 \left({\cal E}\sqrt{\pi/2} + 1\right)}{\left(\sqrt{1 + \alpha^2 {\cal E}^2}\right)^3},} \label{Eqn:SVB1a}
\end{eqnarray}
where in the last equation we have used the fact that ${\displaystyle
 \int^{\infty}_{0} e^{-u^2/2} d u = \sqrt{\pi/2}}$ and that
 ${\displaystyle \int^{\infty}_{0} u e^{-u^2/2} d u = 1}$. 

Comparing equations~(\ref{Eqn:SVAB1}) and~(\ref{Eqn:SVB1a}),
it is clear that we need to seek a ${\cal C} >0$ so that
\begin{equation}\label{Eqn:ConDelta}
{\displaystyle  \frac{1}{\sqrt{1 + \alpha^2 {\cal E}^2}} - 
\frac{\alpha^2 \left({\cal E}\sqrt{\pi/2} + 1\right)}{\left(\sqrt{1 + \alpha^2 {\cal E}^2}\right)^3} 
\geq \frac{1}{\sqrt{1 + \alpha^2 {{\cal C} }^2}}} 
\end{equation}
for all values of $\alpha > 0$. This is equivalent to 
\[{\displaystyle  \frac{1}{\sqrt{1 + \alpha^2 {\cal E}^2}} - \frac{1}{\sqrt{1 + \alpha^2 {{\cal C} }^2}} \geq 
\frac{\alpha^2 \left({\cal E}\sqrt{\pi/2} + 1\right)}{\left(\sqrt{1 + \alpha^2 {\cal E}^2}\right)^3}, } \]
or 
\[{\displaystyle  \frac{{\cal C} ^2 - {\cal E}^2}{\sqrt{1 + \alpha^2 {{\cal C} }^2} \left(
{\sqrt{1 + \alpha^2 {\cal E}^2}} +{\sqrt{1 + \alpha^2 {{\cal C} }^2}}\right)} \geq 
\frac{\left({\cal E}\sqrt{\pi/2} + 1\right)}{\left(\sqrt{1 + \alpha^2 {\cal E}^2}\right)^2}}, \]
which becomes
\[{\displaystyle  {\cal C} ^2 - {\cal E}^2  \geq \left({\cal E}\sqrt{\pi/2} + 1\right) 
\sqrt{\frac{1 + \alpha^2 {{\cal C} }^2}{1 + \alpha^2 {\cal E}^2}} \left(1 + 
\sqrt{\frac{1 + \alpha^2 {{\cal C} }^2}{1 + \alpha^2 {\cal E}^2}} \right).} \]
For ${\cal C} > {\cal E}$, the right hand side reaches its maximum as $\alpha$
approaches $\infty$. Hence it suffices to choose
${\cal C} $ such that   
\begin{equation}\label{Eqn:Delta}
 {\displaystyle  {\cal C} ^2 - {\cal E}^2  \geq  \left({\cal E}\sqrt{\pi/2} + 1\right) 
{\frac{{{\cal C} }}{{\cal E}}} \left(1 + \frac{{{\cal C} }}{{\cal
    E}} \right)}, 
\end{equation}
which solves to 
\[{\displaystyle  {\cal C}  \geq \frac{{\cal E}^3}{{\cal E}^2 - \left({\cal E}\sqrt{\pi/2} + 1\right)}}. \]
For $m \geq 1$ and $n \geq 1$, we have ${\cal E} \geq 5$. The last equation for ${\cal C} $ is easily satisfied when we
choose ${\cal C}  = {\cal E} + 4 = \sqrt{m}+\sqrt{n} + 7$.

We will now take a similar approach to prove
equation~(\ref{Eqn:SVAB2}). We rewrite, by way of function $\widehat{g}(x)$ in~(\ref{Eqn:gghat}), 
\[ {\displaystyle \mathbb{E} \left( {\sqrt{\delta^2+
      \frac{\alpha^2 \|G\|_2^2}{\beta^2 + \gamma^2
      \|G\|_2^2}}}\right)} =  {\displaystyle  
\delta + \mathbb{E} \left( \widehat{g}\left(\|G\|\right)\right) = \delta + \int^{\infty}_{0} \widehat{g}'(x) 
\mathbb{P} \left\{\|G \|_2 \geq x
  \right\} dx . }\]
Since ${\displaystyle \mathbb{P} \left\{\|G \|_2 \geq x \right\} \leq e^{-u^2/2}}$ for $u = x - {\cal E}$, we now have 
\begin{eqnarray}  
 {\displaystyle \mathbb{E}  \left( {\sqrt{\delta^2+
      \frac{\alpha^2 \|G\|_2^2}{\beta^2 + \gamma^2
      \|G\|_2^2}}}\right)}  & \leq &  
 {\displaystyle \delta + \int^{{\cal E}}_{0} \widehat{g}'(x) dx + \int^{\infty}_{{\cal E}} \widehat{g}'(x) \mathbb{P} \left\{\|G \|_2 \geq x   \right\} dx } \nonumber \\
& \leq & {\displaystyle  {\sqrt{\delta^2+
      \frac{\alpha^2 {\cal E}^2}{\beta^2 + \gamma^2
      {\cal E}^2}}} + 
\int^{\infty}_{{\cal E}} 
\frac{\alpha^2 \beta^2 x}{\left(\beta^2 + \gamma^2x^2\right)^2
\sqrt{\delta^2 + \frac{\alpha^2  x^2}{\beta^2 +
    \gamma^2x^2}}} e^{-(x - {\cal E})^2/2} dx }  \nonumber \\
& \leq & {\displaystyle   {\sqrt{\delta^2+
      \frac{\alpha^2 {\cal E}^2}{\beta^2 + \gamma^2
      {\cal E}^2}}} + 
\frac{\alpha^2\beta^2}
{\left(\beta^2 + \gamma^2{\cal E}^2\right)^2
\sqrt{\delta^2 + \frac{\alpha^2  {\cal E}^2}{\beta^2 +
    \gamma^2{\cal E}^2}}}\int^{\infty}_{0} \left({\cal E} + u\right) e^{-u^2/2} d u  }  \nonumber \\
& = & {\displaystyle   {\sqrt{\delta^2+
      \frac{\alpha^2 {\cal E}^2}{\beta^2 + \gamma^2
      {\cal E}^2}}} + 
\frac{\alpha^2\beta^2\left({\cal E}\sqrt{\pi/2} + 1\right)}
{\left(\beta^2 + \gamma^2{\cal E}^2\right)^2
\sqrt{\delta^2 + \frac{\alpha^2  {\cal E}^2}{\beta^2 +
    \gamma^2{\cal E}^2}}}.} \label{Eqn:SVB1b}
\end{eqnarray}

Comparing equations~(\ref{Eqn:SVAB2}) and~(\ref{Eqn:SVB1b}),
we now must seek a ${\cal C} >0$ so that 
\begin{equation}\label{Eqn:SVB1c}
{\displaystyle    {\sqrt{\delta^2+
      \frac{\alpha^2 {\cal E}^2}{\beta^2 + \gamma^2
      {\cal E}^2}}} + 
\frac{\alpha^2\beta^2\left({\cal E}\sqrt{\pi/2} + 1\right)}
{\left(\beta^2 + \gamma^2{\cal E}^2\right)^2
\sqrt{\delta^2 + \frac{\alpha^2  {\cal E}^2}{\beta^2 +
    \gamma^2{\cal E}^2}}} \leq    {\sqrt{\delta^2+
      \frac{\alpha^2 {\cal C} ^2}{\beta^2 + \gamma^2
      {\cal C} ^2}}}} 
\end{equation}
for all values of $\alpha > 0$. Equivalently,
\[{\displaystyle   {\sqrt{\delta^2+
      \frac{\alpha^2 {\cal C} ^2}{\beta^2 + \gamma^2
      {\cal C} ^2}}} - {\sqrt{\delta^2+
      \frac{\alpha^2 {\cal E}^2}{\beta^2 + \gamma^2
      {\cal E}^2}}}  \geq 
\frac{\alpha^2\beta^2\left({\cal E}\sqrt{\pi/2} + 1\right)}
{\left(\beta^2 + \gamma^2{\cal E}^2\right)^2
\sqrt{\delta^2 + \frac{\alpha^2  {\cal E}^2}{\beta^2 +
    \gamma^2{\cal E}^2}}}, } \]
or 
\[{\displaystyle 
\frac{{\displaystyle \frac{{\cal C} ^2-{\cal E}^2}{
\beta^2 + \gamma^2
      {\cal C} ^2}}}
{\displaystyle 
  {\sqrt{\delta^2+
      \frac{\alpha^2 {\cal C} ^2}{\beta^2 + \gamma^2
      {\cal C} ^2}}} + {\sqrt{\delta^2+
      \frac{\alpha^2 {\cal E}^2}{\beta^2 + \gamma^2
      {\cal E}^2}}}}   \geq 
\frac{\left({\cal E}\sqrt{\pi/2} + 1\right)}
{\left(\beta^2 + \gamma^2{\cal E}^2\right)
\sqrt{\delta^2 + \frac{\alpha^2  {\cal E}^2}{\beta^2 +
    \gamma^2{\cal E}^2}}}, } \]
which is the same as
\[{\displaystyle  {\cal C} ^2 - {\cal E}^2 \geq 
\left({\cal E}\sqrt{\pi/2} + 1\right) 
\frac{\left(\beta^2 + \gamma^2 {\cal C} ^2\right)}
{\left(\beta^2 + \gamma^2  {\cal E}^2\right)}
\left(1 + {\displaystyle \sqrt{
\frac{\delta^2+
      \frac{\alpha^2 {\cal C} ^2}{\beta^2 + \gamma^2
      {\cal C} ^2}}
{\delta^2+
      \frac{\alpha^2 {\cal E}^2}{\beta^2 + \gamma^2
      {\cal E}^2}}}}\right).}\]
The right hand side approaches the maximum value as $\gamma$
approaches $\infty$. Hence ${\cal C} $ must satisfy
\[{\displaystyle  {\cal C} ^2 - {\cal E}^2 \geq 
\frac{2\left({\cal E}\sqrt{\pi/2} + 1\right){{\cal C} }^2
}{{\cal E}^2} ,}  \]
which solves to
\[{\displaystyle  {\cal C}  \geq \frac{{\cal E}^2}
{\sqrt{{\cal E}^2 - 2\left({\cal E}\sqrt{\pi/2} + 1\right)}}
 .}  \]
Again the choice ${\cal C}  = \sqrt{m}+\sqrt{n} + 7 = {\cal E}
 + 4$ satisfies this equation. \Q 

The Proof for Proposition~\ref{Lem:Omega1Bound} will follow
a similar track. However, due to the complications with $p
\leq 1$, we will seek help from
Lemma~\ref{Lem:LargDev} instead of
Theorem~\ref{Thm:Gaussfun} to shorten the estimation
process.

{\noindent \bf Proof of Proposition~\ref{Lem:Omega1Bound}:}
As in the proof of Proposition~\ref{Lem:Omega2Bound}, we can
write 
\[ {\displaystyle \mathbb{E} \left( \frac{1}{\sqrt{1+\alpha^2 \|G^{\dagger} \|_2^2}}\right)  
= 1 - \int^{\infty}_{0} g'(x) \mathbb{P} \left\{\|G^{\dagger} \|_2 \geq x  \right\} dx . }\]
By Lemma~\ref{Lem:LargDev} we have for any $x > 0$, 
\[{\displaystyle \mathbb{P} \left\{\|G^{\dagger} \|_2 \geq x  \right\}  \leq \left(
\frac{p+1}{e \sqrt{\ell}} x\right)^{-(p+1)} .} \]
Following arguments similar to those in the proof of Proposition~\ref{Lem:Omega2Bound}, 
we have for a constant $C> 0$ to be later determined, 
\begin{eqnarray}  
 {\displaystyle \mathbb{E} \left( \frac{1}{\sqrt{1+\alpha^2 \|G^{\dagger}\|_2^2}}\right) } & \geq &  
 {\displaystyle 1 - \int^{C}_{0} g'(x) dx - \int^{\infty}_{C} g'(x) \mathbb{P} \left\{\|G^{\dagger} \|_2 \geq x
  \right\} dx } \nonumber \\
& \geq & {\displaystyle  \frac{1}{\sqrt{1 + \alpha^2 C^2}} - 
\int^{\infty}_{C} \frac{\alpha^2 x }{\left(\sqrt{1 + \alpha^2 x^2}\right)^3} \left(
\frac{p+1}{e \sqrt{\ell}} x\right)^{-(p+1)} dx . }
 \label{Eqn:bnd1}
\end{eqnarray}
Below we will derive lower bounds on~(\ref{Eqn:bnd1}) for
the three difference cases of $p$ in
Proposition~\ref{Lem:Omega1Bound}. For $p \geq 2$,
equation~(\ref{Eqn:bnd1}) can be simplified as
\begin{eqnarray*}  
 {\displaystyle \mathbb{E} \left( \frac{1}{\sqrt{1+\alpha^2 \|G^{\dagger}\|_2^2}}\right) } 
& \geq & {\displaystyle  \frac{1}{\sqrt{1 + \alpha^2 C^2}} - 
\frac{\alpha^2}{\left(\sqrt{1 + \alpha^2 C^2}\right)^3} \int^{\infty}_{C} x \left(\frac{p+1}{e \sqrt{\ell}} x\right)^{-(p+1)} dx  }  \nonumber \\  
& =  & {\displaystyle  \frac{1}{\sqrt{1 + \alpha^2 C^2}} - 
\frac{\alpha^2C^2}{(p-1)\left(\sqrt{1 + \alpha^2 C^2}\right)^3} \left(\frac{p+1}{e \sqrt{\ell}}C\right)^{-(p+1)} . }
\end{eqnarray*}
We now seek a ${\cal C} >0$ so that
\[{\displaystyle  \frac{1}{\sqrt{1 + \alpha^2 C^2}} - 
\frac{\alpha^2C^2}{(p-1)\left(\sqrt{1 + \alpha^2 C^2}\right)^3} \left(\frac{p+1}{e \sqrt{\ell}}C\right)^{-(p+1)} 
\geq \frac{1} {\sqrt{1 + \alpha^2 {{\cal C} }^2}}} \]
for all values of $\alpha > 0$. This condition is very
similar to equation~(\ref{Eqn:ConDelta}). Arguments
similar to those used to solve~(\ref{Eqn:ConDelta}) lead to
\[{\displaystyle  {\cal C}  \geq \frac{C^3}{C^2 - 
{\displaystyle \frac{C^2}{p-1} \left(\frac{p+1}{e
    \sqrt{\ell}}C\right)^{-(p+1)}}} .} \]
The choice ${\displaystyle {\cal C}  = \frac{4e
    \sqrt{\ell}}{p+1}}$ satisfies this equation for 
\[ {\displaystyle C = \left(\frac{e
    \sqrt{\ell}}{p+1}\right)
    \left(\frac{p}{p-1}\right)^{1/(p+1)} .} \]

Now we consider the case $p = 1$. We rewrite
      equation~(\ref{Eqn:bnd1}) in light of
      equation~(\ref{Eqn:calfac1}) in~\ref{Sec:ApPrelim}: 
\[ {\displaystyle \mathbb{E} \left( \frac{1}{\sqrt{1+\alpha^2 \|G^{\dagger}\|_2^2}}\right) } \geq   
{\displaystyle  \frac{1}{\sqrt{1 + \alpha^2 C^2}} - \alpha^2
\left(\frac{e \sqrt{\ell}}{2}\right)^2\left(\log \frac{1 + \sqrt{1 +
      \alpha^2 C^2}}{\alpha C} - \frac{1}{\sqrt{1 +
      \alpha^2 C^2}}\right) . }\]
To prove Proposition~\ref{Lem:Omega1Bound}, we just need to find a
constant ${\cal C}$ so that  
\begin{equation}\label{Eqn:bnd1a}
{\displaystyle  \frac{1}{\sqrt{1 + \alpha^2 C^2}} - \alpha^2
\left(\frac{e \sqrt{\ell}}{2}\right)^2\left(\log \frac{1 + \sqrt{1 +
      \alpha^2 C^2}}{\alpha C} - \frac{1}{\sqrt{1 +
      \alpha^2 C^2}}\right) \geq \frac{1}{1 + {\cal C}^2
  \alpha^2 \log {\displaystyle  \frac{2\sqrt{1 + \alpha^2 C^2}}{ \alpha C}}},}
\end{equation}
where the asymptotic term $\alpha^2 \log {\displaystyle
  \frac{2\sqrt{1 + \alpha^2 C^2}}{ \alpha C}}$ behaves like
  $ {\displaystyle O\left(\alpha^2 \log {\displaystyle
  \frac{1}{\alpha}}\right)} $ when $\alpha$ is tiny and like
  $ {\displaystyle O\left(\alpha\right)} $ when $\alpha$ is
  very large. Equation~(\ref{Eqn:bnd1a}) is equivalent to 
\begin{eqnarray} 
{\displaystyle {\cal C}^2 \alpha^2 \log {\displaystyle
  \frac{2\sqrt{1 + \alpha^2 C^2}}{ \alpha C}}} & \geq & 
{\displaystyle \frac{{\displaystyle 1- 
\frac{1}{\sqrt{1 + \alpha^2 C^2}}+ \alpha^2
\left(\frac{e \sqrt{\ell}}{2}\right)^2\left(\log \frac{1 + \sqrt{1 +
      \alpha^2 C^2}}{\alpha C} - \frac{1}{\sqrt{1 +
      \alpha^2 C^2}}\right)}}{{\displaystyle \frac{1}{\sqrt{1 + \alpha^2 C^2}} - \alpha^2
\left(\frac{e \sqrt{\ell}}{2}\right)^2\left(\log \frac{1 + \sqrt{1 +
      \alpha^2 C^2}}{\alpha C} - \frac{1}{\sqrt{1 +
      \alpha^2 C^2}}\right)}}} \nonumber \\
& = & {\displaystyle \frac{{\displaystyle \frac{\alpha^2 C^2}{1+\sqrt{1 + \alpha^2 C^2}}+ \alpha^2
\left(\frac{e \sqrt{\ell}}{2}\right)^2\left(\sqrt{1 + \alpha^2 C^2}\log \frac{1 + \sqrt{1 +
      \alpha^2 C^2}}{\alpha C} -1\right)}}{{\displaystyle 1  - \alpha^2
\left(\frac{e \sqrt{\ell}}{2}\right)^2\left(\sqrt{1 + \alpha^2 C^2}\log \frac{1 + \sqrt{1 +
      \alpha^2 C^2}}{\alpha C} - 1\right)}}} \label{Eqn:bnd1b}.
\end{eqnarray}
All the extra terms involving the $\log $ function have
added much complexity to the above expression. We cut it
down with equations~(\ref{Eqn:calfac5})
and~(\ref{Eqn:calfac6}) in {Appendix}~\ref{Sec:ApPrelim} by replacing
all relevant expressions involving $\psi = \alpha C$ by
their corresponding calculus upper bounds. This gives 
\[{\displaystyle  {\cal C}^2  \alpha^2 \log {\displaystyle
  \frac{2\sqrt{1 + \alpha^2 C^2}}{ \alpha C}}} \geq 
{\displaystyle \frac{{\displaystyle  \alpha^2 C^2\log {\displaystyle
  \frac{2\sqrt{1 + \alpha^2 C^2}}{ \alpha C}} + \alpha^2 
\left(\frac{e \sqrt{\ell}}{2}\right)^2\log {\displaystyle
  \frac{2\sqrt{1 + \alpha^2 C^2}}{ \alpha C}}}}{{\displaystyle  1  -
    \left(\frac{e \sqrt{\ell}}{2C}\right)^2}},} \]
or ${\displaystyle  {\cal C}^2  \geq {\displaystyle \left( C^2 + \left(\frac{e \sqrt{\ell}}{2}\right)^2\right)\left/\right.
\left( 1 - \left(\frac{e
  \sqrt{\ell}}{2C}\right)^2\right)}}$, which holds for $C =
e \sqrt{\ell}$ and ${\cal C} = {\displaystyle  2 e \sqrt{\ell} = \frac{4e \sqrt{\ell}}{p+1}}.$

The last case for our lower bound in
Proposition~\ref{Lem:Omega1Bound} is $p = 0$. With equation~(\ref{Eqn:calfac2}) in~\ref{Sec:ApPrelim}: 
and the choice $C = e \sqrt{\ell}$, equation~(\ref{Eqn:bnd1}) reduces to 
\begin{eqnarray*}  
 {\displaystyle \mathbb{E} \left(
  \frac{1}{\sqrt{1+\alpha^2 \|G^{\dagger}\|_2^2}}\right)}
& \geq & {\displaystyle  \frac{1}{{\displaystyle  \sqrt{1 + \alpha^2 C^2}}} -
  \frac{e \sqrt{\ell} \alpha}{\sqrt{1+\alpha^2C^2}\left(\sqrt{1+\alpha^2C^2}+\alpha
 C\right)} =  \frac{1}{\sqrt{1+\alpha^2C^2}+\alpha C}} \\
& \geq & {\displaystyle  \frac{1}{1+\alpha{\cal C}}}
\end{eqnarray*}
for ${\cal C} = {\displaystyle 
4 e \sqrt{\ell} = \frac{4e \sqrt{\ell}}{p+1}}.$

It is now time to prove equation~(\ref{Eqn:SVAB4}). Our
approach for $p \geq 2$ is similar. We rewrite, by way of
function $\widehat{g}(x)$ in~(\ref{Eqn:gghat}),  
\[ {\displaystyle \mathbb{E} \left( {\sqrt{\delta^2+
      \frac{\alpha^2 \|G^{\dagger}\|_2^2}{\beta^2 + \gamma^2
      \|G^{\dagger}\|_2^2}}}\right)} =  {\displaystyle  
\delta + \mathbb{E} \left( \widehat{g}\left(\|G^{\dagger}\|\right)\right) = \delta + \int^{\infty}_{0} \widehat{g}'(x) 
\mathbb{P} \left\{\|G^{\dagger} \|_2 \geq x
  \right\} dx . }\]
Since for any $x > 0$, 
\[{\displaystyle \mathbb{P} \left\{\|G^{\dagger} \|_2 \geq x  \right\}  \leq \left(
\frac{p+1}{e \sqrt{\ell}} x\right)^{-(p+1)} ,} \]
we now have 
\begin{eqnarray}  
 {\displaystyle \mathbb{E}  \left( {\sqrt{\delta^2+
      \frac{\alpha^2 \|G^{\dagger}\|_2^2}{\beta^2 + \gamma^2
      \|G^{\dagger}\|_2^2}}}\right)}  & \leq &  
 {\displaystyle \delta + \int^{C}_{0} \widehat{g}'(x) dx + \int^{\infty}_{C} \widehat{g}'(x) \mathbb{P} \left\{\|G^{\dagger} \|_2 \geq x   \right\} dx } \nonumber \\
& \leq & {\displaystyle  {\sqrt{\delta^2+
      \frac{\alpha^2 C^2}{\beta^2 + \gamma^2
      C^2}}} + 
\int^{\infty}_{C} 
\frac{\alpha^2 \beta^2 x}{\left(\beta^2 + \gamma^2x^2\right)^2
\sqrt{\delta^2 + \frac{\alpha^2  x^2}{\beta^2 +
    \gamma^2x^2}}} \left(
\frac{p+1}{e \sqrt{\ell}} x\right)^{-(p+1)} dx }  \nonumber \\
& \leq & {\displaystyle   {\sqrt{\delta^2+
      \frac{\alpha^2 C^2}{\beta^2 + \gamma^2
      C^2}}} + 
\frac{\alpha^2\beta^2}
{\left(\beta^2 + \gamma^2C^2\right)^2
\sqrt{\delta^2 + \frac{\alpha^2  C^2}{\beta^2 +
    \gamma^2C^2}}}\int^{\infty}_{C} \left(
\frac{p+1}{e \sqrt{\ell}} x\right)^{-(p+1)} d x  }  \nonumber \\
& = & {\displaystyle   {\sqrt{\delta^2+
      \frac{\alpha^2 C^2}{\beta^2 + \gamma^2
      C^2}}} + 
\frac{\alpha^2\beta^2C^2}
{(p-1) \left(\beta^2 + \gamma^2C^2\right)^2
\sqrt{\delta^2 + \frac{\alpha^2  C^2}{\beta^2 +
    \gamma^2C^2}}} \left(
\frac{p+1}{e \sqrt{\ell}} C\right)^{-(p+1)}.} \nonumber 
\end{eqnarray}
Similarly, we seek a ${\cal C}  > 0$ so that
\[{\displaystyle   {\sqrt{\delta^2+
      \frac{\alpha^2 C^2}{\beta^2 + \gamma^2
      C^2}}} + 
\frac{\alpha^2\beta^2C^2}
{(p-1) \left(\beta^2 + \gamma^2C^2\right)^2
\sqrt{\delta^2 + \frac{\alpha^2  C^2}{\beta^2 +
    \gamma^2C^2}}} \left(
\frac{p+1}{e \sqrt{\ell}} C\right)^{-(p+1)} \leq 
   {\sqrt{\delta^2+
      \frac{\alpha^2 {\cal C} ^2}{\beta^2 + \gamma^2
      {\cal C} ^2}}}.} \]
This last equation is very similar to
      equation~(\ref{Eqn:SVB1c}), with the only difference
      being the coefficients in the second term on the left
      hand side. Thus its solution similarly satisfies
\[{\displaystyle  {\cal C}  \geq \frac{C^2}
{\sqrt{C^2 - \frac{2C^2}{p-1} \left(
\frac{p+1}{e \sqrt{\ell}} C\right)^{-(p+1)}}}} .\]
Again, the value ${\displaystyle {\cal C}  = \frac{4e
    \sqrt{\ell}}{p+1}}$ satisfies this equation for ${\displaystyle C = \left(\frac{e
    \sqrt{\ell}}{p+1}\right)
    \left(\frac{2p}{p-1}\right)^{1/(p+1)} .}$

The special cases $p = 0$ and $p=1$ lead to some involved
calculations with Lemma~\ref{Lem:LargDev}. Instead, we will
appeal to Lemma~\ref{Lem:ChenDon}, an upper bound on the
probability density function of smallest eigenvalue of the
Wishart matrix $GG^T$. It is a happy coincidence that this
upper bound is reasonably tight for $p \leq 1$. By
Lemma~\ref{Lem:ChenDon}, 
\begin{eqnarray}
{\displaystyle \mathbb{E} \left({\sqrt{\delta^2+\frac{\alpha^2 \|G^{\dagger\|_2^2}}{\beta^2 + \gamma^2 \|G^{\dagger}\|_2^2}}}\right) }  & =  & 
 {\displaystyle \delta + \mathbb{E}
 \left({\sqrt{\delta^2+\frac{\alpha^2
 \|G^{\dagger\|_2^2}}{\beta^2 + \gamma^2
 \|G^{\dagger}\|_2^2}}}-\delta\right) }  \nonumber \\
&\leq   & 
 {\displaystyle \delta + L_{\ell-p,\ell}\int^{\infty}_{0} \left({\sqrt{\delta^2+{\displaystyle \frac{\alpha^2 /x}{\beta^2 + \gamma^2 /x}}}} - \delta\right)
e^{-x/2} x^{1/2(p-1)} dx } \nonumber \\
& = &   {\displaystyle \delta +  L_{\ell-p,\ell}\int^{\infty}_{0} \frac{{\displaystyle \frac{\alpha^2 /x}{\beta^2 + \gamma^2 /x}}}
{\sqrt{\delta+{\displaystyle \frac{\alpha^2 /x}{\beta^2 + \gamma^2 /x}}} + \delta} e^{-x/2} x^{1/2(p-1)} dx }  \nonumber \\
&\leq&   {\displaystyle \delta +  L_{\ell-p,\ell}\int^{\infty}_{0} \frac{\alpha^2 }
{\sqrt{\delta^2\left(\beta^2x+\gamma^2\right) + \alpha^2}\sqrt{\beta^2x+\gamma^2}} e^{-x/2} x^{1/2(p-1)} dx .} \label{Eqn:ChenDon2}
\end{eqnarray}
The integral in equation~(\ref{Eqn:ChenDon2}) can be bounded as
\begin{eqnarray}
&& {\displaystyle \int^{\infty}_{0} \frac{\alpha^2 }
{\sqrt{\delta^2\left(\beta^2x+\gamma^2\right) +
    \alpha^2}\sqrt{\beta^2x+\gamma^2}} e^{-x/2} x^{1/2(p-1)}
dx}  \nonumber \\
& \leq & 
{\displaystyle \int^{1}_{0} \frac{\alpha^2 x^{1/2(p-1)} dx}
{\sqrt{\delta^2\left(\beta^2x+\gamma^2\right) + \alpha^2}\sqrt{\beta^2x+\gamma^2}} } + 
{\displaystyle \int^{\infty}_{1} \frac{\alpha^2 e^{-x/2} dx}
{\sqrt{\delta^2\left(\beta^2+\gamma^2\right) + \alpha^2}\sqrt{\beta^2+\gamma^2}}} \nonumber \\
& \leq & {\displaystyle \int^{1}_{0} \frac{\alpha^2 x^{1/2(p-1)} dx}
{\sqrt{\delta^2\left(\beta^2x+\gamma^2\right) +
    \alpha^2}\sqrt{\beta^2x+\gamma^2}} +  
\frac{2\alpha^2
}{\sqrt{\delta^2\left(\beta^2+\gamma^2\right) +
    \alpha^2}\sqrt{\beta^2+\gamma^2}},} \nonumber \\
& \leq & {\displaystyle \int^{1}_{0} \frac{\alpha^2 x^{1/2(p-1)} dx}
{\sqrt{\delta^2\left(\beta^2x+\gamma^2\right) +
    \alpha^2}\sqrt{\beta^2x+\gamma^2}} +  \frac{2\alpha^2
}{\delta\beta^2}} \label{Eqn:ChenDon3}
\end{eqnarray}

Below we further simplify equation~(\ref{Eqn:ChenDon3}). For
$p=1$, the integral in~(\ref{Eqn:ChenDon3}) becomes,
according to equation~(\ref{Eqn:Calfac3})
in~\ref{Sec:ApPrelim}: 
\begin{eqnarray*} 
&& {\displaystyle \frac{2\alpha^2}{\delta \beta^2} \log \frac{\sqrt{\beta^2\left(\delta^2\left(\beta^2+\gamma^2\right) + \alpha^2\right)}  
+ \sqrt{\delta^2\beta^2\left(\beta^2+\gamma^2\right)}}{\sqrt{\beta^2\left(\delta^2\gamma^2 + \alpha^2\right)}  
+ \sqrt{\delta^2\beta^2\gamma^2}}} \\
& \leq &  {\displaystyle \frac{2\alpha^2}{\delta
      \beta^2} \log \frac{2
      \left(\sqrt{\delta^2\left(\beta^2+\gamma^2\right) +
	\alpha^2}\right)}{\sqrt{\delta^2\gamma^2 +
	\alpha^2}}} \nonumber \\
&\leq & {\displaystyle \frac{2\alpha^2}{\delta \beta^2} \log
    2 \sqrt{1 + \frac{\delta^2\beta^2}{\delta^2\gamma^2 + \alpha^2}} \leq 
\frac{2\alpha^2}{\delta \beta^2} \log 2 \sqrt{1 + \frac{\delta^2\beta^2}{\alpha^2}}}.
\end{eqnarray*}
Replacing the integral in equation~(\ref{Eqn:ChenDon3}), and
plugging the resulting upper bound
into equation~(\ref{Eqn:ChenDon2}), we obtain the desired
equation~(\ref{Eqn:SVAB4}) for $p = 1$. 

Finally we consider the case $ p = 0$. The integral in
equation~(\ref{Eqn:ChenDon3}) can be rewritten as
\[ {\displaystyle \int^{1}_{0} \frac{\alpha^2 }
{\sqrt{x}\sqrt{\delta^2\left(\beta^2x+\gamma^2\right) +
    \alpha^2}\sqrt{\beta^2x+\gamma^2}} dx  
= \int^{1}_{0} \frac{ 2 \alpha^2 }
  {\sqrt{\delta^2\left(\beta^2y^2+\gamma^2\right) +
      \alpha^2}\sqrt{\beta^2y^2+\gamma^2}}  dy ,} \]
where we have used the substitution $x = y^2$. Applying the inequality 
\[ A^2 y^2 + B^2 \geq \frac{1}{2} \left(Ay+B\right)^2 \]
to both factors in the denominator above, and utilizing the
identity~(\ref{Eqn:Calfac4}) from~\ref{Sec:ApPrelim}, 
we bound the integral from above as
\begin{eqnarray*}
{\displaystyle \int^{1}_{0} \frac{ 4 \alpha^2 }
  {\left(\delta\beta y+\sqrt{\delta^2\gamma^2 +
      \alpha^2}\right) \left(\beta y+\gamma\right)}} dy &= & {\displaystyle \frac{ 4 \alpha^2 }
  {\delta\beta \gamma - \beta \sqrt{\delta^2\gamma^2 +
      \alpha^2}} \log \frac{\gamma\left(\delta\beta+\sqrt{\delta^2\gamma^2 +
      \alpha^2}\right)}{\sqrt{\delta^2\gamma^2 +
      \alpha^2}\left(\beta+\gamma\right)} } \\
&=&  {\displaystyle \frac{ 4 \alpha^2 \left(\delta \gamma + \sqrt{\delta^2\gamma^2 +
      \alpha^2}\right)}{\beta\alpha^2} \log  
\frac{\sqrt{\delta^2\gamma^2 +
      \alpha^2}\left(\beta+\gamma\right)} {\gamma\left(\delta\beta+\sqrt{\delta^2\gamma^2 +
      \alpha^2}\right)}} \\
&\leq & {\displaystyle \frac{ 8 \sqrt{\delta^2\gamma^2 +
      \alpha^2}}{\beta} \log  \frac{\sqrt{\delta^2\gamma^2 +
      \alpha^2}}{\gamma \delta} \leq 
\frac{ 4 \sqrt{\delta^2\gamma^2 +
      \alpha^2}}{\beta} \log  \left(1 + \left(\frac{\alpha}{\delta\gamma}\right)^2\right),} 
\end{eqnarray*}
which leads to the desired equation~(\ref{Eqn:SVAB4}) for $p
= 0$. \Q

\renewcommand{\theequation}{S5.\arabic{equation}}
\renewcommand{\thetheorem}{S5.\arabic{theorem}}
\renewcommand{\thesection}{S5}

\section{Facts from Calculus} \label{Sec:ApPrelim}
Here we list the facts we have used from calculus. Their
proofs have been left out, since they do not provide any
additional insight into our analysis. We start with $4$
definite integrals:
\begin{eqnarray}
{\displaystyle \int^{\infty}_{C} \frac{dx}{x\left(\sqrt{1 +
      \alpha^2 x^2}\right)^3}} & = & {\displaystyle \log \frac{1 + \sqrt{1 +
      \alpha^2 C^2}}{\alpha C} - \frac{1}{\sqrt{1 +
      \alpha^2 C^2}}}, \label{Eqn:calfac1} \\
{\displaystyle \int_{C}^{\infty}
\frac{dx}{\left(1+\alpha^2x^2\right)^{3/2}}} & = & {\displaystyle
\frac{1}{\alpha
  \sqrt{1+\alpha^2C^2}\left(\sqrt{1+\alpha^2C^2}+\alpha
  C\right)},}  \label{Eqn:calfac2} \\
{\displaystyle  \int_0^1
  \frac{dx}{\sqrt{\left(Ax+B\right)\left(Cx+D\right)}}}  &=& {\displaystyle \frac{2}{\sqrt{AC}} \log 
\frac{\sqrt{C\left(A+B\right)} + \sqrt{A\left(C+D\right)}}{
\sqrt{BC} + \sqrt{AD}}} \label{Eqn:Calfac3} \\
{\displaystyle \int_0^1
  \frac{dx}{{\left(Ax+B\right)\left(Cx+D\right)}}} & =& {\displaystyle  \frac{1}{{AD-BC}} \log
\frac{D\left(A+B\right)}{B\left(C+D\right)},}\label{Eqn:Calfac4}
\end{eqnarray}
where $\alpha$, $A, B, C, D$ are all positive constants. We
will also list the following inequalities for any $\psi >0$: 
\begin{eqnarray} 
 {\displaystyle \max\left(\frac{1}{1+ \sqrt{1+\psi^2}},
 {\displaystyle \sqrt{1+\psi^2}\log
 \frac{1+\sqrt{1+\psi^2}}{\psi} 
  - 1}\right)} & \leq & {\displaystyle \log \frac{2\sqrt{1+\psi^2}}{\psi}}     ,\label{Eqn:calfac5} \\
{\displaystyle \psi^2 \left(\sqrt{1+\psi^2}\log \frac{1+\sqrt{1+\psi^2}}{\psi}
  - 1\right)} & \leq & 1 . \label{Eqn:calfac6}
\end{eqnarray}

\bibliographystyle{plain}

\end{document}